\documentclass{article}
\usepackage{amsthm,amsfonts,amssymb,epsfig,graphics,amsmath,fancyheadings}

\def\R{\textrm{Re\ }}

\newcommand{\tyek}{t^{-1}}
\newcommand{\tnim}{t^{-\frac 12}}
\newcommand{\trob}{t^{-\frac 14}}

\newcommand{\ubar}{\bar{u}}

\newcommand{\gcal}{\mathcal{G}}

\newcommand{\fcal}{\mathcal{F}}

\newcommand{\ddp}{\frac{\partial{\bar u}^\delta}{\partial \delta}}
\newcommand{\mathbbR}{\mathbb{R}}

\theoremstyle{plain}
\newtheorem{theorem}{Theorem}
\newtheorem{lemma}{Lemma}
\newtheorem{proposition}{Proposition}

\newtheorem{remark}{Remark}
\newtheorem{definition}{Definition}

\numberwithin{equation}{section}

\title{Sharp pointwise bounds for perturbed viscous shock waves
}\bigskip

\author{\sc \small
Peter Howard\thanks{Texas A\& M University, College Station, TX
77843; phoward@math.tamu.edu} \, Mohammadreza
Raoofi\thanks{Indiana University, Bloomington, IN 47405;
mraoofi@indiana.edu. } \, and Kevin Zumbrun
\thanks{Indiana University, Bloomington, IN 47405; kzumbrun@indiana.edu.
 } }

\begin{document}

\maketitle

\begin{abstract}
Refining previous work in \cite{Z.3, MaZ.3, Ra, HZ, HR}, we derive sharp
pointwise bounds on behavior of perturbed viscous shock profiles
for large-amplitude Lax or overcompressive type shocks and physical viscosity.
These extend well-known results of Liu \cite{Liu97} 
obtained by somewhat different techniques for small-amplitude Lax type 
shocks and artificial viscosity, completing a program set out in \cite{ZH}.
As pointed out in \cite{Liu91, Liu97}, the key to obtaining sharp bounds
is to take account of cancellation associated with the property
that the solution decays faster along characteristic than in other
directions.
Thus, we must here estimate characteristic derivatives for the entire
nonlinear perturbation, rather than judicially chosen parts 
as in \cite{Ra, HR}.
a requirement that greatly complicates the analysis.
\end{abstract}

\section{Introduction}\label{introduction}


In the landmark paper \cite{Liu97}, Liu
established sharp pointwise bounds on the asymptotic
behavior of a perturbed viscous shock profile, for small-amplitude,
Lax type shocks and artificial viscosity.
A long-standing program of the authors, set out in \cite{ZH}, 
has been to extend to
large-amplitudes, general type profiles, and physical (partially parabolic)
viscosities results obtained by Liu and others for small amplitude
Lax type profiles and artificial viscosity.
For various results in this direction, see, e.g., \cite{Z.3, MaZ.3, Ra, HZ, HR}.
In this paper, we achieve a definitive result 
recovering the full bounds of \cite{Liu97} for general Lax or 
overcompressive type shocks and physical or artificial viscosity,
thus completing the program of \cite{ZH}. 
The analysis involves an interesting blend of the spectral
techniques introduced by the authors to deal with large 
amplitudes,  described in Sections \ref{preliminaries} and
\ref{nonlinear},
a delicate type of cancellation estimate introduced by Liu
in the small-amplitude context, described in Section \ref{liu},
and sharp 
convolution estimates of the type developed in
\cite{H.1, H.2, H.3, HZ, HR}, taking account of transversality
but not cancellation of two interacting signals.
%
%

Consider a (possibly) large-amplitude
{\it viscous shock profile}, or traveling-wave solution 
\begin{equation} \label{profile}
\bar{u} (x - st); \quad \lim_{x \to \pm \infty} \bar{u} (x) = u_{\pm}, 
\end{equation}
of systems of partially or fully parabolic conservation laws
\begin{equation} \label{main}
\begin{aligned}
u_t + F(u)_x &= (B(u)u_x)_x, \\
\end{aligned}
\end{equation}
$x \in \mathbb{R}$, $u$, $F \in \mathbb{R}^n$, $B\in\mathbb{R}^{n \times n}$, 
where
\begin{equation}
u = \begin{pmatrix} u^I \\ u^{II} \end{pmatrix},
\quad
B = \begin{pmatrix} 0 & 0 \\ b_1 & b_2 \end{pmatrix},
\end{equation}
$u^I \in \mathbb{R}^{n-r}$, $u^{II} \in \mathbb{R}^r$, $r$ some positive integer, 
possibly $n$ (full regularization), and 
\begin{equation*}
\text{Re }\sigma(b_2) \ge \theta > 0.
\end{equation*} 
Here and elsewhere, $\sigma$ denotes spectrum of a matrix or other
linear operator.  
Working in a coordinate system moving along with the shock, we may without loss
of generality consider a standing profile $\bar{u} (x)$, $s=0$.

\medskip
Following \cite{Z.2}, we assume that, 
by some invertible change of coordinates $u\to w(u)$, 
followed if necessary by multiplication on the left
by a nonsingular matrix function $S(w)$, equations (\ref{main})
may be written in the {\it quasilinear, partially symmetric
hyperbolic-parabolic form}
\begin{equation}
\tilde A^0 w_t +  \tilde A w_{x}= (\tilde B w_{x})_{x} + G, \quad
w=\left(\begin{matrix} w^I
\\w^{II}\end{matrix}\right), \label{symm}
\end{equation}
$w^I\in \mathbb{R}^{n-r}$, $w^{II}\in \mathbb{R}^r$, $x\in
\mathbb{R}$, $t\in \mathbb{R}_+$, where, defining $w_\pm:= w(u_\pm)$:\\
\medskip
(A1)\quad $\tilde A(w_\pm)$, $\tilde A_{11}$, $\tilde
A^0$ are symmetric, $\tilde A^0 >0$.\\
\medskip
(A2)\quad Dissipativity: no eigenvector of $ dF(u_\pm)$ lies in the
kernel of $ B(u_\pm)$. (Equivalently, no eigenvector of $ \tilde A
(\tilde A^0)^{-1} (w_\pm)$
lies in the kernel of $\tilde B (\tilde A^0)^{-1}(w_\pm)$.)\\
\medskip
(A3) \quad $ \tilde B= \left(\begin{matrix} 0 & 0 \\ 0 & \tilde b
\end{matrix}\right) $, $ \tilde G= \left(\begin{matrix}  0 \\ \tilde g\end{matrix}\right) $,
with $ Re  \tilde b(w) \ge \theta $ for some $\theta>0$, for all
$W$, and $\tilde g(w_x,w_x)=\mathbf{O}(|w_x|^2)$.
\medskip
Here, the coefficients of (\ref{symm}) may be expressed in terms of
the original equation (\ref{main}), the coordinate change $u\to
w(u)$, and the approximate symmetrizer $S(w)$, as
\begin{equation}
\begin{aligned} \tilde A^0&:= S(w)(\partial u/\partial w),\quad
\tilde A:= S(w)d F(u(w))(\partial u/\partial w),\\
\tilde B&:= S(w)B(u(w))(\partial u/\partial w), \quad G= -(dS w_{x})
B(u(w))(\partial u/\partial w) w_{x}.
\end{aligned}
\label{coeffs}
\end{equation}
Alternatively, we assume, simply,

\medskip
(B1)\quad Strict parabolicity: $n=r$, or, equivalently, $\Re \sigma(B)>0$.\\
\medskip

Along with the above structural assumptions, we make
the technical
hypotheses:\\
\medskip
(H0)\quad $F$, $B$, $w$, $S\in C^{10}$.\\
\smallskip
(H1)\quad The eigenvalues of $\tilde A_*:=\tilde{A}_{11}(\tilde{A}^0_{11})^{-1}$
are (i) distinct from $0$; (ii) of common sign; and (iii) of
constant multiplicity with respect to
$u$.\\
\smallskip
(H2)\quad The eigenvalues of $dF(u_\pm)$ are real, distinct, and
nonzero.\\
\smallskip
({H3})\quad Nearby $\bar u$, the set of all solutions of
\eqref{profile}--\eqref{main} connecting the same values $u_\pm$
forms a smooth manifold $\{\bar u^\delta\}$, $\delta\in
\mathcal{U}\subset \mathbb{R}^\ell$, $\bar u^0=\bar u$.

\begin{remark}\label{profileegs}
\textup{
Structural assumptions (A1)--(A3) [alt. (B1)] 
and technical hypotheses (H0)--(H2) admit such physical systems
as the compressible Navier--Stokes equations, the equations
of magnetohydrodymics, and Slemrod's model for van der Waal gas
dynamics \cite{Z.1}.   Moreover, existence of waves $\bar{u}$
satisfying (H3) has been established in each of these cases.}
\end{remark}

\bigskip
\begin{definition}\label{type}
An ideal shock
\begin{equation} \label{shock}
u(x,t) =
\begin{cases}
u_- &x < st, \\
u_+ &x > st,
\end{cases}
\end{equation}
is classified as
{\it undercompressive}, {\it Lax},
or {\it overcompressive} type according as $i-n$
is less than, equal to, or greater than $1$,
where $i$, denoting the sum of the dimensions $i_-$ and $i_+$
of the center--unstable subspace of $df(u_-)$ and the center--stable
subspace of $df(u_+)$, represents the total number of characteristics
incoming to the shock.
                                                                                                                                                             
A viscous profile \eqref{profile} is classified
as {\it pure undercompressive} type if the associated ideal shock
is undercompressive and $\ell=1$,
{\it pure Lax} type
if the corresponding ideal shock is Lax type and $\ell=i-n$,
and {\it pure overcompressive} type if
the corresponding ideal shock is overcompressive and $\ell=i-n$,
$\ell$ as in (H3).
Otherwise it is classified as {\it mixed under--overcompressive} type;
see \cite{ZH}.
\end{definition}

Pure Lax type profiles are the most common type, and the only type
arising in standard gas dynamics, while pure over- and undercompressive type
profiles arise in magnetohydrodynamics (MHD) and phase-transitional models.
In this paper, we restrict to the case of pure Lax or 
overcompressive type shocks.

Finally, we assume that the profile satisfies a 
linearized stability criterion based on the 
{\it Evans function}.  As described, e.g., in
\cite{AGJ, E, GZ, J, KS, MaZ.3, Z.3, ZH},
the Evans function $D(\lambda)$, 
a Wronskian constructed from solutions of the associated eigenvalue
equation, serves as a characteristic
function for the linear operator $L$ that arises upon
linearization of (\ref{main}) about $\bar{u} (x)$.  More
precisely, away from essential spectrum, zeros of the
Evans function correspond in location and multiplicity with
eigenvalues of $L$ \cite{AGJ, GZ, ZH}.  It was shown in
\cite{ZH} and \cite{MaZ.3}, respectively for the
strictly parabolic and real viscosity cases, that
$L^1 \cap L^p \to L^p$ linearized orbital stability of the
profile, $p>1$, is equivalent to the Evans function
condition
\medbreak
($\mathcal{D}$) \quad There exist
precisely $\ell$ zeroes of $D(\cdot)$ in the nonstable half-plane
$\R \lambda \ge 0$, necessarily at $\lambda=0$,
where $\ell$ as in (H3) is the dimension of the manifold connecting
$u_-$ and $u_+$.
\medbreak\noindent
Under assumptions (A0)--(A3) [alt. (B1)] 
and (H0)--(H3), condition ($\mathcal{D}$) is equivalent to (i)
{\it strong spectral stability}, $\sigma(L)\subset \{\R \lambda \le
0\}\cup \{0\}$, (ii) {\it hyperbolic stability} of the associated
ideal shock, and (iii) {\it transversality} of $\bar u$ as a
solution of the connection problem in the associated traveling-wave
ODE, where hyperbolic stability is defined for Lax and
undercompressive shocks by the Lopatinski condition of \cite{M.1,
M.2, M.3, Fre} and for Lax and overcompressive shocks by the analogous
long-wave stability condition ($\mathcal{D}$ii), below; see
\cite{ZH, MaZ.2, ZS, Z.2} for further explanation. 

\begin{remark}\label{satisfaction}
\textup{
Stability criterion $(\mathcal{D})$ has been shown to hold for 
general small amplitude Lax shocks \cite{FreS, Go.1, HuZ.2, KM, KMN, MN, PZ}, 
and for large-amplitude shocks in such cases as Lax type waves 
arising in isentropic Navier--Stokes equations for the gamma-law gas 
as $\gamma \to 1$ \cite{MN}, and undercompressive shocks arising in Slemrod's
model for van der Waal gas dynamics \cite{Z.1} (see \cite{S} for
Slemrod's model). 
On the other hand, it has been shown to fail for certain
large-amplitude and or nonclassical type shocks \cite{GZ, FreZ, Z.1, ZS}.
More generally, condition $(\mathcal{D})$ can be
readily checked by numerical calculation 
\cite{B.1, B.2, BZ, BDG, HuZ.1, HSZ}. }
\end{remark}

%
%

\medskip
Setting $A_\pm:=df(u_\pm)$, $\Gamma_\pm:=d^2f(u_\pm)$, and $B_\pm:=B(u_\pm)$,
denote by
\begin{equation}
a_1^-< a_2^-<\dots < a_n^- \quad \text{\rm and }
a_1^+< a_2^+<\dots < a_n^+
\label{a}
\end{equation}
the eigenvalues of $A_-$ and $A_+$, and $l_j^\pm$, $r_j^\pm$ left
and right eigenvectors associated with each $a_j^\pm$, normalized so
that $(l_j^{T}r_k)_\pm=\delta^j_k$, where $\delta^j_k$ is the
Kronecker delta function, returning $1$ for $j=k$ and $0$ for $j\ne
k$.  Define scalar diffusion coefficients
\begin{equation}
\beta_j^\pm:= (l_j^T Br_j)_\pm
\label{beta}
\end{equation}
and scalar coupling coefficients
\begin{equation}
\gamma_j^\pm:= (l_j^T \Gamma (r_j,r_j))_\pm.
\label{gamma}
\end{equation}
Under this notation,  hyperbolic stability (a consequence
of the assumed ($\mathcal{D})$) of a Lax or
overcompressive shock profile $\bar u$ is the condition:
\medbreak ($\mathcal{D}$ii)\quad The set $\{r_j^\pm: a_j^\pm \gtrless
0\} \cup \{\int_{-\infty}^{+\infty}\frac{\partial
\bar{u}^\delta}{\partial \delta_i} dx; i=1, \cdots, \ell \}$ forms a
basis for $\mathbb{R}^n$, with
$\int_{-\infty}^{+\infty}\frac{\partial \bar{u}^\delta}{\partial
\delta_i} dx$  computed at $\delta=0.$ \medskip \\ 

\begin{remark}\label{majda}
\textup{
For Lax profiles, ($\mathcal{D}$ii)
reduces to the Liu--Majda condition
$$
\det \begin{pmatrix}
\{r_j^\pm: a_j^\pm \gtrless 0\} & (u_+-u_-)
\end{pmatrix}\ne 0.
$$
}
\end{remark}

Following \cite{Liu85, Liu97},
define for a given mass $m_j^-$ the
scalar diffusion waves $\varphi_j^-(x,t;m_j^-)$ as
(self-similar) solutions of the Burgers equations
\begin{equation}
\varphi_{j,t}^- + a^-_j\varphi^-_{j,x} - \beta^-_j \varphi_{j,xx}^- = -\gamma^-_j
((\varphi_{j}^-)^2)_x
\label{diffusionwaves}
\end{equation}
with point-source initial data
\begin{equation}
\varphi_j^- (x, -1) = m_j^-\delta_0(x),
\label{burgersdata}
\end{equation}
and similarly for $\varphi_j^+(x,t;m_j^+)$.
Given a collection of masses
$m_j^\pm$
prescribed on outgoing characteristic modes
$a_j^-<0$ and $a_j^+>0$, define
\begin{equation}
\varphi(x,t) =\sum_{a_j^- <0}\varphi_j^-(x,t; m_j^-) r^-_j
+ \sum_{a_j^+>0}\varphi_j^+(x,t; m_j^+) r^+_j.
\label{phi}
\end{equation}


In the setting described above, we will determine estimates 
on perturbed viscous shock profiles in terms of $\phi$
and a refined collection of ~template~ functions (terminology following
\cite{ZH}; notation following [Liu97]).

\begin{definition}[Template functions.] 
Let
\begin{equation} \label{templates}
\begin{aligned}
\psi_1 (x, t) &= \sum_{a_k^\pm \gtrless 0} (1 + |x - a_k^\pm t| + t^{1/2})^{-3/2} \\
\psi_2 (x, t) &= (1+|x|)^{-1/2} (1 + |x| + t)^{-1/2} (1 + |x| + t^{1/2})^{-1/2} \chi \\
\psi_3 (x, t) &= (1 + |x| + t)^{-1} (1 + |x|)^{-1} \chi \\
\psi_4 (x, t) &= (1 + |x| + t)^{-7/4} \chi,
\end{aligned}
\end{equation}
where $\chi$ denotes an indicator function on $x \in [a_1^- t, a_n^+ t]$,
and also  
\begin{equation*}
\begin{aligned}
\psi_1^{j, \pm} (x, t) &= (1 + |x - a_j^\pm t| + t^{1/2})^{-3/2} \\
\bar{\psi}_1^{j,\pm} (x, t) &= \psi_1 (x, t) - \psi_1^{j, \pm}.
\end{aligned}
\end{equation*}
\end{definition}

The goal of our analysis is to establish the following 
sharp pointwise description of asymptotic behavior, generalizing
bounds obtained by Liu \cite{Liu97} for small-amplitude
Lax type profiles with artificial viscosity $B=I$.

\begin{theorem} \label{pointwiseestimates}
Assume (A1)--(A3) [alt. (B1)], (H0)--(H3) and $(\mathcal{D})$ hold, 
and  $\bar{u}$ is a pure Lax or overcompressive shock profile.
Assume also that $\tilde{u}$ solves (\ref{main}) with
initial data $\tilde{u}_0$ and that, for initial perturbation
$u_0:=\tilde{u}_0-\bar{u}$, we have $|u_0|_{H^5} \le E_0$
and
$|u_0(x)|$, $|\partial_x u_0(x)|$, 
and $|\partial_x^2 u_0(x)|\le E_0(1+|x|)^{-\frac32}$, 
for $E_0$ sufficiently small. 
Then, the solution $\tilde u$ continues globally in time,
with $\tilde u-\bar u\in L^1\cap H^5$.
Moreover, there exist a choice of $|m_j^\pm|$, $|\delta_*|\le CE_0$ 
and a function $\delta(t)$ (determined, respectively,
by (\ref{masses} and (\ref{delta})), 
such that, for $v:=\tilde{u}-\bar{u}^{\delta_*}-
\varphi-\frac{\partial \bar{u}^\delta}{\partial \delta}(\delta_*)\delta$, 
\begin{equation} \label{1stestimate}
\begin{aligned}
|v (x,t)| &\le C E_0 \Big[\psi_1 + \psi_2 \Big] (x, t) \\
|\partial_x v (x,t)| &\le C E_0 \Big[ t^{-1/2} (\psi_1+\psi_2) + \psi_3 + \psi_4 \Big] (x,t) \\  
|\delta (t)|&\le C E_0(1+t)^{-\frac 12} \\
|\dot \delta (t)| &\le C E_0(1+t)^{-1}
\end{aligned}
\end{equation}
and, for all $a_j^\pm \gtrless 0$,
\begin{equation}\label{1stest2}
|(\partial_t + a_j^\pm \partial_x) v(x,t)|
\le 
CE_0 \Big[t^{-1} (1+t)^{1/4} \psi_1^{j,\pm} + t^{-1/2} (\bar{\psi}_1^{j,\pm}
+ \psi_2) + \psi_3 + \psi_4 \Big] (x, t),
\end{equation}
for some constant $C$ (independent of $x,t$ and $E_0$).
\end{theorem}

\begin{remark}\label{taylor}
\textup{
\bigskip  By Taylor's Theorem,
\begin{equation}
\partial_x^k(\bar u^{\delta_*+\delta(t)}-\bar u^{\delta_*}-\ddp\delta(t) )=
\mathbf O (|\delta(t)|^2 e^{-k|x|}),\label{ishl}
\end{equation}
for $k=0$ [resp. $k=1$],
since $\ddp$ and its derivatives decay at exponential rate 
(see Lemma \ref{expdecay}). 
Since the righthand side of (\ref{ishl}) for $k=0$
[resp. $k=1$] is
smaller than the right hand side of (\ref{1stestimate})(i) [resp. (ii)], 
we may conclude that
$$\tilde v:=  u- \bar u^{\delta_*+\delta(t)}-\phi$$
decays at the same rate as $v$, i.e.,
$u$ is well approximated by $\bar u^{\delta_*+\delta(t)}+\phi$:
a dynamically changing nearby profile with the same endstates $u_\pm$,
superposed with outgoing diffusion waves.
This is a slight refinement of the asymptotic
description $u\sim \bar u^{\delta_*}+\phi$ introduced 
by Liu \cite{Liu85, Liu97}.
}
\end{remark}

\begin{remark}\label{modeests}
\textup{
In contrast to the modulus bounds given here,
the estimates of \cite{Liu97} are stated in
each separate characteristic field, giving
the further information that along characteristic
rays the perturbation $v$ lies mainly in the
associated characteristic direction.
The difference is rather subtle, however, and so
we have opted for simplicity of exposition to
omit this level of detail.  The additional information
can easily be recovered at the expense of additional
bookkeeping; see Remark \ref{refinedcharbds} for further discussion.
}
\end{remark}

%
Liu's analysis involved essentially two main ingredients,
which were strongly coupled.
The first was to obtain approximate Green function bounds
taking advantage of the weakly coupled nature of the equations
in the small-amplitude case to approximate by a superposition of 
solutions of scalar conservation laws; the second, by delicate
pointwise interaction estimates between the approximate Green
kernel and various algebraically decaying source terms, to close
a nonlinear iteration and obtain the result.
The ``coupling'' we mention refers to the fact that the Green
function estimates blow up as amplitude goes to zero and
the shock becomes more and more characteristic, whereas the
source terms decay with amplitude; thus, the two effects must
be delicately balanced to close the iteration and achieve
a correct result.
In the large-amplitude case, bounded away from the characteristic
limit, the issues are somewhat different; namely, we do not
have this ``characteristic coupling'' problem, but on the other
hand we cannot as in \cite{Liu97} obtain approximate Green
function bounds by asymptotic development in the amplitude.
Likewise, for physical, partially parabolic viscosity, there
are new difficulties associated with regularity and the need
to gain derivatives.

These new difficulties have largely been surmounted in the
authors' previous work.  In particular, (i) sharp Green function
bounds have been obtained by Mascia--Zumbrun \cite{MaZ.2} in great 
generality using Laplace transform/stationary phase estimates, and 
(ii) global existence and sharp $H^s\cap L^p$ estimates 
on the nonlinear residual $v$ have been obtained by Raoofi \cite{Ra}
using the linearized estimates of \cite{MaZ.3}, a key cancellation estimate
of Liu, and a nonstandard energy estimate to control higher derivatives.
These results are described in more detail in Sections \ref{preliminaries}
and \ref{liu}, and we shall use them freely in our analysis.
A major advantage of this ``bootstrap'' approach is that 
we need not close a nonlinear iteration, but only carry out a {\it linear} 
fixed-point argument to obtain our result.
This was used in the previous work \cite{HR} to establish``nearly optimal'' 
pointwise bounds in a relatively uncomplicated manner.
As pointed out in \cite{ZH, HR}, however,
to obtain the full bounds of Liu requires estimates not only
on $v$ and $v_x$ as in \cite{HR}, but also {\it characteristic derivatives}
$(\partial_t +a_j^\pm\partial_x)v$ as in \cite{Liu97}.
To carry out these bounds requires an immense ammount of additional
work, with consideration of numerous different cases, and accounts for
most of the work of this paper.

\begin{remark}\label{sp}
\textup{
In the strictly parabolic case, the regularity requirement
on the initial perturbation $u=\tilde u-\bar u$ may
be significantly relaxed in Theorem \ref{pointwiseestimates},
from $u\in H^5$ to $u\in C^{0+\alpha}$, $\alpha>0$,
and the bound $E+0(1+|x|)^{-3/2}$ imposed on $|u_0|$ alone,
and not derivatives $|\partial_x u_0|$ and $|\partial_x^2 u_0|$,
by using parametrix bounds as described in \cite{ZH, Z.2, HZ, HR} to
establish short-time pointwise smoothing estimates.
(That these bounds are necessary in the real viscosity case,
on the other hand,
may be clearly seen in the proof of Proposition \ref{Hlinear}.)
In the artificial viscosity case $B\equiv \text{\rm constant}$,
which includes the case $B\equiv I$ treated by Liu \cite{Liu97},
the short-time existence/regularity theory becomes trivial 
\cite{ZH, HZ}, and we require
no additional regularity beyond the single pointwise hypothesis
$|u_0|\le E+0(1+|x|)^{-3/2}$.
}
\end{remark}

\begin{remark}\label{energysp}
\textup{
Alternatively, we may subsume the strictly parabolic case (B1) under
the same analysis used to treat the partially symmetric case (A1)--(A3),
substituting for parametrix estimates the more elementary Sobolev
estimates of \cite{MaZ.2, Ra}.
For, in this case, we may drop
the requirement of symmetry of $A_\pm$, which
was used only to obtain a skew-symmetric $K$ such that $\Re (KA+B)_\pm>0$,
whereas $K=0$ suffices when $B$ is full rank.
Since the $A_{11}$ component is empty,
the remaining structural conditions are trivially
satisfied upon multiplying \ref{main} by a smooth symmetric positive
definite $A^0(u)$ such that $\Re A^0 B>0$ (guaranteed by Lyapunov's Theorem).
By a more careful accounting, taking advantage of associated
simplifications, the results of \cite{Z.2, Ra}, with the
exception of derivative bounds,
may be obtained in the strictly parabolic case under regularity
$H^1$, matching the assumption of Goodman in the classic
paper \cite{Go.1} stability with respect to zero-mass perturbations,
and the results of \cite{HR} (for which are required $L^\infty$
bounds on $\partial_x u$) under regularity $H^2$.
In the present case, for which pointwise derivative bounds are crucial to the
argument, we require $H^4$.
}
\end{remark}

{\bf Plan of the paper.} In Section \ref{preliminaries},
we recall the results of \cite{MaZ.2, HR} for our use.
In Section \ref{nonlinear}, we give the brief argument
establishing Theorem \ref{pointwiseestimates}, 
subject to certain integral estimates.
In Section \ref{liu}, we motivate the analysis to follow
by reviewing a key cancellation estimate of Liu \cite{Liu85}
upon which our analysis, and that of \cite{Liu97}, depends.
In Sections \ref{integralestimatesL}, 
\ref{integralestimatesNLI}, 
and \ref{integralestimatesNLII} we carry out the
main work of the paper, establishing the deferred integral estimates.

\section{Preliminaries}\label{preliminaries}

We start by recalling some needed, known results.

\subsection{Profile estimates}\label{profileestimates}

We first recall the profile analysis carried out in \cite{MaZ.2},
generalizing results of \cite{MP} in the strictly parabolic case.
Profile $\bar u (x)$ satisfies
the standing-wave ordinary differential equation (ODE)
\begin{equation}
B(\bar u) \bar u'=F(\bar u)-F(u_-). \label{ODE}
\end{equation}
Considering the block structure of $B$, this can be written as:
\begin{equation}
F^I(u^I, u^{II})\equiv F^I(u_-^I, u_-^{II})\label{eq1}
\end{equation}
and
\begin{equation}
b_1(u^I)' + b_2(u^{II})'= F^{II}(u^I, u^{II}) - F^{II}(u_-^I,
u_-^{II}). \label{eq2}
\end{equation}

\begin{proposition} \label{expdecay}[\cite{MaZ.2}]
Given (H1)--(H3), (\ref{eq2}) determines a smooth $r$-dimensional manifold
on which (\ref{eq2}) determines a nondegenerate ODE.  Moreover, 
endstates $u_\pm$ are hyperbolic restpoints of this ODE,
i.e., the coefficients of the linearized equations
about $u_\pm$, written in local coordinates, have no center
subspace. In particular, under regularity (H0),
\begin{equation}
D_x^j D_\delta^i(\bar
u^\delta(x)-u_{\pm})=\mathbf{O}(e^{-\alpha|x|}) \, \,
\text{\rm as $x\rightarrow\pm\infty$},  \qquad \alpha>0, \,
0\le j\le 10, \,i=0,1.
\end{equation}
\end{proposition}

\subsection{Linearized equations and Green distribution bounds} \label{linearized}

We next recall some linear theory from \cite{MaZ.2, ZH, HZ}.
Linearizing (\ref{main}) about $\ubar^{\delta_*}(\cdot)$,
$\delta_*$ to be determined later, gives 
\begin{equation}
v_t=Lv:=-(Av)_x+(Bv_x)_x, \label{linearov}
\end{equation}
with
\begin{equation}
B(x):= B(\ubar^{\delta_*}(x)), \quad  A(x)v:=
df(\ubar^{\delta_*}(x))v-dB(\ubar^{\delta_*}(x))v\ubar^{\delta_*}_x.
\label{AandBov}
\end{equation}
Denoting $A^\pm := A(\pm \infty)$,  $B^\pm:= B(\pm \infty)$, and
considering  Lemma \ref{expdecay}, it follows that
\begin{equation}
|A(x)-A^-|= \mathbf {O} (e^{-\eta |x|}), \quad |B(x)-B^-|= \mathbf
{O} (e^{-\eta |x|}) \label{ABboundsov}
\end{equation}
as $x\to -\infty,$ for some positive $\eta.$ Similarly for $A^+$ and
$B^+,$ as $x\to +\infty.$ Also $|A(x)-A^\pm|$ and $|B(x)-B^\pm|$ are
bounded for all $x$.

 Define the {\it (scalar) characteristic speeds} $a^\pm_1<
\cdots < a_n^\pm$ (as above) to be the eigenvalues of $A^\pm$, and
the left and right {\it (scalar) characteristic modes} $l_j^\pm$,
$r_j^\pm$ to be corresponding left and right eigenvectors,
respectively (i.e., $A^\pm r_j^\pm = a_j^\pm r_j^\pm,$ etc.),
normalized so that $l^+_j \cdot r^+_k=\delta^j_k$ and $l^-_j \cdot
r^-_k=\delta^j_k$. Following Kawashima \cite{Kaw}, define associated
{\it effective scalar diffusion rates} $\beta^\pm_j:j=1,\cdots,n$ by
relation
\begin{equation}
\left(
\begin{matrix}
\beta_1^\pm &&0\\
&\vdots &\\
0&&\beta_n^\pm
\end{matrix}
\right) \quad = \hbox{diag}\ L^\pm B^\pm R^\pm, \label{betaov}
\end{equation}
where $L^\pm:=(l_1^\pm,\dots,l_n^\pm)^t$,
$R^\pm:=(r_1^\pm,\dots,r_n^\pm)$ diagonalize $A^\pm$.

Assume for  $A$ and $B$ the block structures:
$$A=\left(\begin{matrix}A_{11}\quad A_{12}\\A_{21}\quad A_{22}\end{matrix}\right),
B=\left(\begin{matrix}0& 0\\B_{21}& B_{22}\end{matrix}\right).$$


Also, let $a^{*}_j(x)$, $j=1,\dots,(n-r)$ denote the eigenvalues of
$$
A_{*}:= A_{11}- A_{12} B_{22}^{-1}B_{21}, \label{A*}
$$
with $l^*_j(x)$, $r^*_j(x)\in \mathbbR^{n-r}$ associated left and
right eigenvectors, normalized so that $l^{*t}_jr_j\equiv
\delta^j_k$. More generally, for an $m_j^*$-fold eigenvalue, we
choose $(n-r)\times m_j^* $ blocks $L_j^*$ and $R_j^*$ of
eigenvectors satisfying the {\it dynamical normalization}
$$
L_j^{*t}\partial_x R_j^{*}\equiv 0,
$$
along with the usual static normalization $L^{*t}_jR_j\equiv
\delta^j_kI_{m_j^*}$; as shown in Lemma 4.9, \cite{MaZ.1}, this may
always be achieved with bounded $L_j^*$, $R_j^*$. Associated with
$L_j^*$, $R_j^*$, define extended, $n\times m_j^*$ blocks
$$
\mathcal{L}_j^*:=\left(\begin{matrix} L_j^* \\
0\end{matrix}\right), \quad \mathcal{R}_j^*:=
\left(\begin{matrix} R_j^*\\
-B_{22}^{-1}B_{21} R_j^*\end{matrix}\right). \label{CalLR}
$$
%
Eigenvalues $a_j^*$ and eigenmodes $\mathcal{L}_j^*$,
$\mathcal{R}_j^*$ correspond, respectively, to short-time hyperbolic
characteristic speeds and modes of propagation for the reduced,
hyperbolic part of degenerate system (\ref{main}).

Define local, $m_j\times m_j$ {\it dissipation coefficients}
$$
\eta_j^*(x):= -L_j^{*t} D_* R_j^* (x), \quad j=1,\dots,J\le n-r,
\label{eta}
$$
where
$$
\aligned
&{D_*}(x):= 
 \, A_{12}B_{22}^{-1} \Big[A_{21}-A_{22} B_{22}^{-1} B_{21}+ A_{*}
B_{22}^{-1} B_{21} + B_{22}\partial_x (B_{22}^{-1} B_{21})\Big]
\endaligned
\label{D*}
$$
is an effective dissipation analogous to the effective diffusion
predicted by formal, Chapman--Enskog expansion in the (dual)
relaxation case.


The {\it Green distribution} (fundamental solution)
associated with (\ref{linearov}) is defined by
\begin{equation}
G(x,t;y):= e^{Lt}\delta_y (x). \label{ov3.7}
\end{equation}
Recalling the standard notation
$
\textrm{errfn} (z) := \frac{1}{2\pi} \int_{-\infty}^z e^{-\xi^2}
d\xi,
$
we have the following pointwise description. 

\begin{proposition} \label{greenbounds}\cite{MaZ.2}  Under assumptions 
(A1)--(A3) [alt. (B1)],
(H0)--(H3), and $(\mathcal D)$,
the Green distribution $G(x,t;y)$
associated with the linearized evolution equations  may be
decomposed as
\begin{equation}
G(x,t;y)= H + E+  S + R, \label{ourdecomp}
\end{equation}
where, for $y\le 0$:
\begin{equation}
\begin{aligned} H(x,t;y)&:= \sum_{j=1}^{J} a_j^{*-1}(x) a_j^{*}(y)
\mathcal{R}_j^*(x) \zeta_j^*(y,t) \delta_{x-\bar a_j^* t}(-y)
\mathcal{L}_j^{*t}(y)\\
&= \sum_{j=1}^{J} \mathcal{R}_j^*(x) \mathcal{O}(e^{-\eta_0 t})
\delta_{x-\bar a_j^* t}(-y) \mathcal{L}_j^{*t}(y),
\end{aligned}
\label{multH}
\end{equation}
where the averaged convection rates $\bar a_j^*= \bar a_j^*(x,t)$ in
(\ref{multH}) denote the time-averages over $[0,t]$ of $a_j^*(x)$
along backward characteristic paths $z_j^*=z_j^*(x,t)$ defined by
$$
dz_j^*/dt= a_j^*(z_j^*), \quad z_j^*(t)=x, 
$$
and the dissipation matrix $\zeta_j^*=\zeta_j^*(x,t)\in
\mathbb{R}^{m_j^*\times m_j^*}$ is defined by the {\it dissipative
flow}
$$
d\zeta_j^*/dt= -\eta_j^*(z_j^*)\zeta_j^*, \quad
\zeta_j^*(0)=I_{m_j};
$$
and $\delta_{x-\bar a_j^* t}$ denotes  Dirac distribution centered
at $x-\bar a_j^* t$.
\begin{equation}\label{E}
E(x,t;y):=\sum_{j=1}^\ell \frac{\partial \bar
u^\delta(x)}{\partial \delta_j}_{|\delta=\delta_*}e_j(y,t),
\end{equation}
\begin{equation}\label{e}
  e_j(y,t):=\sum_{a_k^{-}>0}
  \left(\textrm{errfn }\left(\frac{y+a_k^{-}t}{\sqrt{4\beta_k^{-}t}}\right)
  -\textrm{errfn }\left(\frac{y-a_k^{-}t}{\sqrt{4\beta_k^{-}t}}\right)\right)
  l_{jk}^{-}
\end{equation}
\begin{equation}
\begin{aligned} S(x,t;y)&:= \chi_{\{t\ge 1\}}\sum_{a_k^-<0}r_k^-  {l_k^-}^t
(4\pi \beta_k^-t)^{-1/2} e^{-(x-y-a_k^-t)^2 / 4\beta_k^-t} \\
&+\chi_{\{t\ge 1\}} \sum_{a_k^- > 0} r_k^-  {l_k^-}^t (4\pi
\beta_k^-t)^{-1/2} e^{-(x-y-a_k^-t)^2 / 4\beta_k^-t}
\left({\frac {e^{-x}}{e^x+e^{-x}}}\right)\\
&+ \chi_{\{t\ge 1\}}\sum_{a_k^- > 0, \,  a_j^- < 0}
[c^{j,-}_{k,-}]r_j^-  {l_k^-}^t (4\pi \bar\beta_{jk}^- t)^{-1/2}
e^{-(x-z_{jk}^-)^2 / 4\bar\beta_{jk}^- t}
\left({\frac{e^{ -x}}{e^x+e^{-x}}}\right),\\
&+\chi_{\{t\ge 1\}} \sum_{a_k^- > 0, \,  a_j^+ > 0}
[c^{j,+}_{k,-}]r_j^+  {l_k^-}^t (4\pi \bar\beta_{jk}^+ t)^{-1/2}
e^{-(x-z_{jk}^+)^2 / 4\bar\beta_{jk}^+ t}
\left({\frac{e^{ x}}{e^x+e^{-x}}}\right),\\
\end{aligned}
\label{Sov}
\end{equation}
with
\begin{equation}
z_{jk}^\pm(y,t):=a_j^\pm\left(t-\frac{|y|}{|a_k^-|}\right)
\label{zjkov}
\end{equation}
and
\begin{equation}
\bar \beta^\pm_{jk}(x,t;y):= \frac{x^\pm}{a_j^\pm t} \beta_j^\pm +
\frac{|y|}{|a_k^- t|} \left( \frac{a_j^\pm}{a_k^-}\right)^2
\beta_k^-, \label{betaaverageov}
\end{equation}
The remainder $R$ and its derivatives have the following bounds.
\begin{equation}
\begin{aligned}
R(x,t;y)&=
\mathbf{O}(e^{-\eta(|x-y|+t)})\\
&+\sum_{k=1}^n \mathbf{O} \left( (t+1)^{-1/2} e^{-\eta x^+}
+e^{-\eta|x|} \right)
t^{-1/2}e^{-(x-y-a_k^{-} t)^2/Mt} \\
&+ \sum_{a_k^{-} > 0, \, a_j^{-} < 0} \chi_{\{ |a_k^{-} t|\ge |y|
\}} \mathbf{O} ((t+1)^{-1/2} t^{-1/2})
e^{-(x-a_j^{-}(t-|y/a_k^{-}|))^2/Mt}
e^{-\eta x^+}, \\
&+ \sum_{a_k^{-} > 0, \, a_j^{+}> 0} \chi_{\{ |a_k^{-} t|\ge |y| \}}
\mathbf{O} ((t+1)^{-1/2} t^{-1/2}) e^{-(x-a_j^{+}
(t-|y/a_k^{-}|))^2/Mt}
e^{-\eta x^-}, \\
\end{aligned}
\label{Rbounds}
\end{equation}
\begin{equation}
\begin{aligned}
R_y(x,t;y)&= \sum_{j=1}^J \mathbf{O}(e^{-\eta t})\delta_{x-\bar
a_j^* t}(-y) +
\mathbf{O}(e^{-\eta(|x-y|+t)})\\
&+\sum_{k=1}^n \mathbf{O} \left( (t+1)^{-1/2} e^{-\eta x^+}
+e^{-\eta|x|}+e^{-\eta|y|}  \right)
t^{-1}
e^{-(x-y-a_k^{-} t)^2/Mt} \\
&+ \sum_{a_k^{-} > 0, \, a_j^{-} < 0} \chi_{\{ |a_k^{-} t|\ge |y|
\}} \mathbf{O} ((t+1)^{-1/2} t^{-1})
e^{-(x-a_j^{-}(t-|y/a_k^{-}|))^2/Mt}
e^{-\eta x^+} \\
&+ \sum_{a_k^{-} > 0, \, a_j^{+} > 0} \chi_{\{ |a_k^{-} t|\ge |y|
\}} \mathbf{O} ((t+1)^{-1/2} t^{-1})
e^{-(x-a_j^{+}(t-|y/a_k^{-}|))^2/Mt}
e^{-\eta x^-}, \\
\end{aligned}
\label{Rybounds}
\end{equation}
\begin{equation}
\begin{aligned}
R_x(x,t;y)&= \sum_{j=1}^J \mathbf{O}(e^{-\eta t})\delta_{x-\bar
a_j^* t}(-y) +
\mathbf{O}(e^{-\eta(|x-y|+t)})\\
&+\sum_{k=1}^n \mathbf {O} \left( (t+1)^{-1} e^{-\eta x^+}
+e^{-\eta|x|} \right)
t^{-1} (t+1)^{1/2}
e^{-(x-y-a_k^{-} t)^2/Mt} \\
&+ \sum_{a_k^{-} > 0, \, a_j^{-} < 0} \chi_{\{ |a_k^{-} t|\ge |y|
\}} \mathbf{O}((t+1)^{-1/2} t^{-1})
e^{-(x-a_j^{-}(t-|y/a_k^-|))^2/Mt}
e^{-\eta x^+} \\
&+ \sum_{a_k^{-} > 0, \, a_j^{+} > 0} \chi_{\{ |a_k^{-} t|\ge |y|
\}} \mathbf{O}((t+1)^{-1/2} t^{-1})
e^{-(x-a_j^{+}(t-|y/a_k^{-}|))^2/Mt}
e^{-\eta x^-}. \\
\end{aligned}
\label{Rxbounds}
\end{equation}
Moreover, for $|x-y|/t$ sufficiently large, $|G|\le Ce^{-\eta
t}e^{-|x-y|^2/Mt)}$ as in the strictly parabolic case.
\end{proposition}

Setting $\tilde{G} := S+R$, so that $G = H + E + \tilde{G}$,  we have the
following useful alternative bounds for $\tilde{G}$.

\begin{proposition}[\cite{ZH}, \cite{MaZ.2, HZ}]\label{altGFbounds}
Under the assumptions of Proposition \ref{greenbounds}, $\tilde{G}$ has
the following bounds.
\begin{equation}\label{Gbounds}
\begin{aligned}
|\partial_{x,y}^\alpha &\tilde G(x,t;y)|\le 
\sum_{j=1}^J\sum_{\beta=0}^{\max \{0, |\alpha|-1\}} \mathbf{O}(e^{-\eta t})\partial_y^\beta \delta_{x-\bar
a_j^* t}(-y) + \mathbf{O}(e^{-\eta(|x-y|+t)})\\
&+ C(t^{-|\alpha|/2} + |\alpha_x| e^{-\eta|x|}) \Big( \sum_{k=1}^n
t^{-1/2}e^{-(x-y-a_k^{-} t)^2/Mt} e^{-\eta x^+} \\
&+ \sum_{a_k^{-} > 0, \, a_j^{-} < 0} \chi_{\{ |a_k^{-} t|\ge |y|
\}} t^{-1/2} e^{-(x-a_j^{-}(t-|y/a_k^{-}|))^2/Mt}
e^{-\eta x^+}, \\
&+ \sum_{a_k^{-} > 0, \, a_j^{+}> 0} \chi_{\{ |a_k^{-} t|\ge |y| \}}
t^{-1/2} e^{-(x-a_j^{+} (t-|y/a_k^{-}|))^2/Mt}
e^{-\eta x^-}\Big), \\
|(\partial_t + a_l^- \partial_x) &\tilde G(x,t;y)|\le 
\sum_{j=1}^J \mathbf{O}(e^{-\eta t})\delta_{x-\bar
a_j^* t}(-y) + \mathbf{O}(e^{-\eta(|x-y|+t)})\\
&+C t^{-3/2}
\Big( e^{-(x-y-a_l^{-} t)^2/Mt} e^{-\eta x^+} 
+ \sum_{a_k^- > 0} e^{-(x-a_l^{-}(t-|y/a_k^{-}|))^2/Mt} e^{-\eta x^+} \Big) \\
&+
Ct^{-1/2} 
\Big( e^{-(x-y-a_l^{-} t)^2/Mt} e^{-\eta |x|} I_{\{x \ge 0\}} 
+ \sum_{a_k^- > 0} e^{-(x-a_l^{-}(t-|y/a_k^{-}|))^2/Mt} e^{-\eta |x|} I_{\{x \ge 0\}} \Big) \\
& + C(t^{-1} + e^{-\eta|x|}) \Big( \sum_{k \ne l}
t^{-1/2}e^{-(x-y-a_k^{-} t)^2/Mt} e^{-\eta x^+} \\
&+ \sum_{a_k^{-} > 0, \, a_j^{-} < 0, j \ne l} \chi_{\{ |a_k^{-} t|\ge |y|
\}} t^{-1/2} e^{-(x-a_j^{-}(t-|y/a_k^{-}|))^2/Mt}
e^{-\eta x^+}, \\
&+ \sum_{a_k^{-} > 0, \, a_j^{+}> 0} \chi_{\{ |a_k^{-} t|\ge |y| \}}
t^{-1/2} e^{-(x-a_j^{+} (t-|y/a_k^{-}|))^2/Mt}
e^{-\eta x^-}\Big), \\
\end{aligned}
\end{equation}
$0\le |\alpha|\le 2$,
for $y\le 0$, and symmetrically for $y\ge 0$, for some $\eta$, $C$,
$M>0$, where $a_j^\pm$ are as in Proposition \ref{greenbounds},
$\beta_k^\pm>0$, $x^\pm$ denotes the positive/negative part of $x$,
and indicator function $\chi_{\{ |a_k^{-}t|\ge |y| \}}$ is $1$ for
$|a_k^{-}t|\ge |y|$ and $0$ otherwise. Moreover, all estimates are
uniform in the supressed parameter $\delta_*$.
\end{proposition}

\begin{remark}
\textup{ In the strictly parabolic case 
$r=n$, the hyperbolic part $H$ is
absent in the decomposition of $G$ given in (\ref{ourdecomp}). }
\end{remark}
\begin{remark}\label{eboundsrmk}
\textup{ {}From \eqref{e}, we obtain by straightforward calculation
(see \cite{MaZ.2}) the bounds
\begin{equation}\label{ebounds}
\begin{aligned}
&|e_j(y,t)|\le C\sum_{a_k^->0}
  \left(\textrm{errfn }\left(\frac{y+a_k^{-}t}{\sqrt{4\beta_k^{-}t}}\right)
  -\textrm{errfn }\left(\frac{y-a_k^{-}t}{\sqrt{4\beta_k^{-}t}}\right)\right),\\
&|\partial_t  e_j(y,t)|\le C t^{-1/2} \sum_{a_k^->0} e^{-|y+a_k^-t|^2/Mt},\\
&|\partial_y  e_j(y,t)|\le C t^{-1/2} \sum_{a_k^->0} e^{-|y+a_k^-t|^2/Mt}\\
&|\partial_{yt}  e_j(y,t)|\le C
t^{-1} \sum_{a_k^->0} e^{-|y+a_k^-t|^2/Mt}\\
\end{aligned}
\end{equation}
for $y\le 0$, and symmetrically for $y\ge 0$. }
\end{remark}

\subsection{$L^p$ decay} \label{Lp}

Finally, we recall the $H^s\cap L^p$ result of \cite{Ra}
established using the linearized bounds of \cite{MaZ.2} and
Hausdorff--Young's inequality together
with an $H^s$ energy estimate and the cancellation estimate 
described in Section \ref{liu}, below.

\begin{proposition}[\cite{Ra, HR}]\label{propra}
Under the conditions of Theorem \ref{pointwiseestimates},
there exists a unique, global solution 
$\tilde u$ of (\ref{main}), $\tilde u\in L^1\cap H^5$.
Moreover, for $v:=\tilde{u}-\bar{u}^{\delta_*}-
\varphi-\frac{\partial \bar{u}^\delta}{\partial \delta}(\delta_*)\delta$, 
with $|m_j^\pm|$, $|\delta_*|\le CE_0$ 
and $\delta(t)$ defined as in (\ref{masses} and (\ref{delta}), 
\begin{equation}\label{vlp}
|v(\cdot,t)|_{L^p}\le
 C E_0
(1+t)^{-\frac{1}{2}(1-\frac1p)-\frac 14},
\qquad
1 \leq p \leq \infty,
\end{equation}
and
\begin{equation}\label{vhs}
|v(\cdot,t)|_{H^5}\le
 C E_0
(1+t)^{-\frac{1}{2}},
\end{equation}
for some constant $C>0$,
\end{proposition}

\section{Nonlinear analysis: Proof of Theorem 
\ref{pointwiseestimates}}\label{nonlinear}

We now carry out the proof of Theorem \ref{pointwiseestimates}
assuming certain integral estimates to be established in 
Sections \ref{integralestimatesL}, \ref{integralestimatesNLI},
and \ref{integralestimatesNLII}. 
Let $\tilde{u}$ be the solution of (\ref{main}) guaranteed by
Proposition \ref{propra}. 
Following \cite{Liu85, Liu97, Ra}, let
$|m_i|$, $|\delta_*|\le CE_0$ be the unique solutions guaranteed
by ($\cal{D}$ii) and the Implicit Function Theorem of
\begin{equation}\label{masses}
\int_{-\infty}^{+\infty} (\tilde{u} (x, 0) -\bar{u}^{\delta_*})(x)\, dx
= \sum_{a_j^- <0}m_j r_j^- + \sum_{a_j^+ >0}m_j r_j^+,
\end{equation}
or, equivalently,
$\int_{-\infty}^{+\infty} \tilde{u} (x, 0) 
=\int_{-\infty}^{+\infty} (\bar{u}^{\delta_*}+\phi)(x)\, dx.
$
This determines the asymptotic state, by conservation of mass.

\begin{remark}[\cite{Liu85, ZH}]\textup{
In the case of Lax--type shock waves,
$\bar{u}^\delta(x)=\bar{u}(x+\delta),$ 
hence $\int\frac{\partial\bar{u}^\delta}{\partial\delta}\, dx=
(u_+-u_-)$ and
$\delta_*$ can be explicitly computed as the solution of a linear equation.}
\end{remark}

Setting $u(x,t):=\tilde{u}(x,t)-\bar{u}^{\delta_*}(x)$, use Taylor's
expansion around $\bar{u}^{\delta_*}(x)$ to find
\begin{equation}
u_t + (A(x)u)_x - (B(x)u_x)_x = -(\Gamma(x)(u,u))_x + Q(u, u_x)_x,
\label{utaylor}
\end{equation}
where $\Gamma(x)(u,u) = d^2f(\bar{u}^{\delta_*})(u,u)-
d^2B(\bar{u}^{\delta_*})(u,u)\bar{u}^{\delta_*}_x$ and $$Q(u, u_x) =
\mathbf {O} (|u||u_x|+|u|^3),$$  
with $\Gamma_\pm = \Gamma(\pm
\infty).$  Define constant coefficients $b^{\pm}_{ij}$ and
$\Gamma_{ijk}^{\pm}$ to satisfy
\begin{equation}
\Gamma_\pm (r^\pm_j, r^\pm_k) = \sum_{i=1}^n \Gamma_{ijk}^{\pm}
r^\pm_i,\quad B_\pm r^\pm_j= \sum_{i=1}^n b^{\pm}_{ij} r^\pm_i.
\label{coefshock}
\end{equation}
Then, of course, $\beta^\pm_i = b^{\pm}_{ii}$ and $\gamma^\pm_i :=
\Gamma_{iii}^{\pm}.$

Continuing, set 
$v:= u -\varphi-\frac{\partial\bar{u}^\delta}{\partial\delta}\delta(t)$, 
with $\delta(t)$ to be defined
later, and assuming $\delta(0)=0$. Notice that, by our choice of
$\delta_*$ and diffusion waves $\varphi_i^\pm$'s, we have zero initial
mass of $v$, i.e.,
\begin{equation} \label{initial0mass}
\int_{-\infty}^{+\infty} v(x,0)dx = 0.
\end{equation}
Replacing $u$ with $v + \varphi+\frac{\partial\bar{u}^\delta}{\partial\delta}\delta(t)$ in (\ref{utaylor}) (
$\frac{\partial \bar{u}^\delta}{\partial \delta_i} $  computed at
$\delta=\delta_*$), and using the fact that
$\frac{\partial\bar{u}^\delta}{\partial \delta_i}$ satisfies the
linear time independent equation $Lv=0$, we obtain
\begin{equation}
v_t - Lv = \Phi(x,t)+ \mathcal{F}(\varphi, v,  \frac{\partial\bar{u}^\delta}{\partial\delta}\delta(t))_x 
\frac{\partial\bar{u}^\delta}{\partial\delta}+ \dot
\delta (t), \label{khati}
\end{equation}
where
\begin{equation}\label{fcalt}
\begin{split}
\mathcal{F}(\varphi, v, \frac{\partial\bar{u}^\delta}{\partial\delta}\delta &) = \mathbf {O}  ( |v|^2 + |\varphi||v|
+ |v||\frac{\partial\bar{u}^\delta}{\partial\delta}\delta|
+ |\varphi|| \frac{\partial\bar{u}^\delta}{\partial\delta} \delta|
+|\frac{\partial\bar{u}^\delta}{\partial\delta}\delta|^2 \\
&+ |(\varphi + v+\frac{\partial\bar{u}^\delta}{\partial\delta}\delta)
(\varphi + v + \frac{\partial\bar{u}^\delta}{\partial\delta}\delta)_x| + |\varphi
+v+\frac{\partial\bar{u}^\delta}{\partial\delta}\delta|^3 ).
\end{split}
\end{equation}

Furthermore,
\begin{equation}
\begin{aligned}
\mathcal{F}(v, \varphi, \frac{\partial\bar{u}^\delta}{\partial\delta}\delta(t))_x& =
\mathbf{O}\big(\mathcal{F}(v,\varphi,\frac{\partial\bar{u}^\delta}{\partial\delta}\delta)\\
&+
|(v+\varphi+\frac{\partial\bar{u}^\delta}{\partial\delta}\delta)_x|
|(v^{II}+\varphi+\frac{\partial\bar{u}^\delta}{\partial\delta}\delta)_x|\\
&+|v+\varphi+\frac{\partial\bar{u}^\delta}{\partial\delta}\delta|
|(v^{II}+\varphi+\frac{\partial\bar{u}^\delta}{\partial\delta}\delta)_{xx}|\big),
\label{fcalx}
\end{aligned}
\end{equation}
and $\Phi(x,t) := -\varphi_t - (A(x)\varphi)_x  + (B(x)\varphi_x)_x
-(\Gamma(x)(\varphi , \varphi))_x$. For $\Phi$  we write
\begin{equation}
\begin{split}
\Phi(x,t) = &-(\varphi_t +A \varphi_x - B \varphi_{xx} + \Gamma(\varphi, \varphi)_x) \\
= &-\sum_{a_i^- < 0} \varphi_t^i r_i^- + (A(x)\varphi^i r_i^-)_x -
(B(x)\varphi_x^i r_i^-)_x + (\Gamma(x)(\varphi^ir_i^-,\varphi^ir_i^-))_x\\
&- \sum_{a_i^+ > 0} \varphi_t^i r_i^+ + (A(x)\varphi^i r_i^+)_x -
(B(x)\varphi_x^i r_i^+)_x + (\Gamma(x)(\varphi^ir_i^+,\varphi^ir_i^+))_x\\
& - \sum_{i\neq j}(\varphi_i \varphi_j \Gamma(x)(r_i^\pm, r_j^\pm))_x.
\end{split} \label{ghati}
\end{equation}
Let us write a typical term of the first summation ($a_i^- < 0$)
in the following form:
\begin{equation}
\begin{split}
\varphi_t^i & r_i^- + (A(x)\varphi^i r_i^-)_x -
(B(x)\varphi_x^i r_i^-)_x + (\Gamma(x)(\varphi^ir_i^-,\varphi^ir_i^-))_x\\
&= \left[(A(x)-A^-)\varphi^i r_i^- -
(B(x)-B^-)\varphi_x^i r_i^- + (\Gamma(x)-\Gamma^-)(\varphi^ir_i^-,\varphi^ir_i^-)\right]_x\\
&+ \varphi_t^i r_i^- + (\varphi^i_xA^- r_i^-) - (\varphi_{xx}^i B^- r_i^-) +
((\varphi^{i})^2_x\Gamma^-(r_i^-,r_i^-)).
\end{split} \label{ghatitar}
\end{equation}
Now we use the definition of $\varphi^i$ in (\ref{diffusionwaves}) and the
definition of coefficients $b_{ij}$ and $\Gamma_{ijk}$ in
(\ref{coefshock}) to write the last part of (\ref{ghatitar}) in the
following form:
\begin{equation}
\begin{split}
\varphi_t^i r_i^- + (\varphi^i_xA^- r_i^-) - &(\varphi_{xx}^i B^- r_i^-) +
((\varphi^i)^2_x\Gamma^-(r_i^-,r_i^-))\\
&= -\varphi^i_{xx}\sum_{j\ne i} b^-_{ij}r^-_j - (\varphi^i)^2_{x}\sum_{j\ne
i} \Gamma^-_{jii}r^-_j.
\end{split} \label{ajabaa}
\end{equation}
Similar statements hold for $a_i^+ > 0$ with minus signs replaced
with plus signs.
\medskip

Applying Duhamel's principle, we obtain from (\ref{khati})
\begin{equation}
\begin{aligned} v(x,t)
&=\int^{+\infty}_{-\infty} G(x,t;y)v_0(y)dy\\
&+\int^t_0
\int^{+\infty}_{-\infty}G(x,t-s;y)\mathcal{F}(\varphi, v, \frac{\partial\bar{u}^\delta}{\partial\delta}\delta)_y(y,s) dy \, ds,\\
&+ \int^t_0\int^{+\infty}_{-\infty}G(x,t-s;y)\Phi(y,s) dy
\, ds 
+ \frac{\partial\bar{u}^\delta}{\partial\delta} \delta(t),\\
 \label{shswgreen}
\end{aligned}
\end{equation}
where we have used the identity
$\int_{-\infty}^{+\infty}
G(x,t;y)\frac{\partial\bar{u}^\delta}{\partial
\delta_i}(y)dy=e^{Lt}\frac{\partial\bar{u}^\delta}{\partial
\delta_i}=\frac{\partial\bar{u}^\delta}{\partial \delta_i}$ and
$\delta(0)=0$.
Assuming
\begin{equation}
\begin{aligned}
 \delta_i (t)
&=-\int^\infty_{-\infty}e_{i}(y,t) v_0(y)dy\\
&-\int^t_0\int^{+\infty}_{-\infty} e_{i}(y,t-s)\mathcal{F}(\varphi,
v,\frac{\partial\bar{u}^\delta}{\partial\delta}\delta)_y(y,s)
 dy ds\\
&-\int^t_0\int^{+\infty}_{-\infty}e_i(x,t-s;y)\Phi(y,s)dy \, ds,
\end{aligned} \label{delta}
\end{equation}
and using (\ref{shswgreen}),  (\ref{delta}) and $G=H+E+\tilde{G}$ we
obtain:
\begin{equation}
\begin{aligned} v(x,t)
&=\int^{+\infty}_{-\infty} (H+ \tilde{G})(x,t;y)v_0(y)dy\\
&+\int^t_0
\int^{+\infty}_{-\infty}(H+\tilde{G})(x,t-s;y)\mathcal{F}(\varphi, v,\frac{\partial\bar{u}^\delta}{\partial\delta}\delta)_y(y,s) dy \, ds\\
&+ \int^t_0\int^{+\infty}_{-\infty}(H+\tilde{G})(x,t-s;y)\Phi(y,s) dy
\, ds,  \label{tilde{G}}
\end{aligned}
\end{equation}
which we can clearly augment with derivative estimates through differentiation 
on both sides.
In addition to $v(x, t)$ and $\delta (t)$, we will keep track in our
argument of $v_x (x, t)$, $(v_t+a^\pm_jv_x)$, and $\dot \delta (t)$, 
the latter of which satisfies
\begin{equation}
\begin{aligned}
 \dot \delta_i (t)
&=-\int^\infty_{-\infty}\partial_t e_{i}(y,t) v_0(y)dy\\
&-\int^t_0\int^{+\infty}_{-\infty} \partial_t e_{i}(y,t-s)\mathcal{F}(\varphi,
v,\frac{\partial\bar{u}^\delta}{\partial\delta}\delta)_y(y,s)
 dy ds\\
&-\int^t_0\int^{+\infty}_{-\infty} \partial_t e_i(y,t-s)\Phi(y,s)dy \, ds,
\end{aligned} \label{dotdelta}
\end{equation}
where we have taken advantage of the observation, apparent
from (\ref{ebounds}), that $e(y,0)=0$.

Define
\begin{equation} \label{zetadefined}
\begin{aligned}
\zeta(t) &:= \sup_{y, 0 \le s \le t} \frac{|v (y, s)|}{\psi_1 (y, s) + \psi_2 (y,s)} \\
&+
\sup_{y, 0 \le s \le t}
\frac{|v_y (y, s)|}{s^{-1/2} (\psi_1 (y, s) + \psi_2 (y,s)) + \psi_3 (y, s) + \psi_4 (y, s)} \\
&+
\sum_{a_j^- < 0, a_j^+ > 0}
\sup_{y, 0 \le s \le t}
\frac{|(\partial_s + a_j^\pm \partial_y) v (y, s)|}
{s^{-1} (1+s)^{1/4} \psi_1^{j, \pm} (y, s) + s^{-1/2} (\bar{\psi}_1^{j,\pm} (y, s) + \psi_2 (y, s))
+ \psi_3 (y, s) + \psi_4 (y, s)} \\
&+
\sup_{0 \le s \le t} \frac{|\delta (s)|}{(1+s)^{-1/2}}
+
\sup_{0 \le s \le t} \frac{|\dot{\delta} (s)|}{(1+s)^{-1}}.
\end{aligned}
\end{equation}
The desired estimates follow easily 
using the following integral estimates, to be established
in Sections \ref{integralestimatesL}, \ref{integralestimatesNLI},
and \ref{integralestimatesNLII}.

\begin{lemma}[Linear estimates I.] \label{linearintegralestimates}
If $\int_{-\infty}^{+\infty} v_0 (y) dy = 0$ and
$|v_0 (y)| \le E_0 (1 + |y|)^{-3/2}$, $E_0 > 0$, then
\begin{equation}\label{Ginit}
\begin{aligned}
\Big|\int_{-\infty}^{+\infty} \tilde{G} (x, t; y) v_0 (y) dy \Big|
&\le 
C E_0 \psi_1 (x, t) \\
\Big|\int_{-\infty}^{+\infty} \partial_x \tilde{G} (x, t; y) v_0 (y) dy \Big|
&\le
CE_0 (t^{-1/2} + e^{-\eta |x|}) \psi_1 (x, t) \\
\Big|\int_{-\infty}^{+\infty}(\partial_t + a_j^\pm \partial_x) \tilde{G} (x, t; y) v_0 (y) dy \Big|
&\le
CE_0 \Big(t^{-1} \psi_1^{j, \pm} (x,t) + t^{-1/2} \bar{\psi}_1^{j, \pm}(x,t)
+ e^{-\eta |x|} \psi_1 (x, t) \Big) \\
\Big|\int_{-\infty}^{+\infty} e_i (y,t) v_0 (y) dy \Big|
&\le 
C E_0 (1 + t)^{-1/2} \\
\Big|\int_{-\infty}^{+\infty} \partial_t e_i(y, t) v_0 (y) dy \Big|
&\le 
C E_0 (1+t)^{-3/2}.
\end{aligned}
\end{equation}
\end{lemma}

\begin{lemma}[Linear estimates II.]\label{Hlinear}
If $|v_0(x)|$, $|\partial_x v_0(x)|$,
$|\partial_x^2 v_0(x)|\le E_0 (1+|x|)^{-\frac32}$, 
$E_0>0$,  then, for some $\theta>0$,
\begin{equation}\label{Hinit}
\begin{aligned}
\int_{-\infty}^{+\infty} H(x, t; y) v_0(y) dy&\le C E_0 e^{-\theta
t}(1+|x|)^{-\frac32}
\le CE_0\psi_1(x,t), \\
\int_{-\infty}^{+\infty} H_x(x, t; y) v_0(y) dy&\le C E_0 e^{-\theta
t}(1+|x|)^{-\frac12}
\le CE_0 (1+t)^{-1/2}\psi_1(x,t), \\
\int_{-\infty}^{+\infty}(\partial_t+a_j^\pm\partial_x) 
H_x(x, t; y) v_0(y) dy&\le C E_0 e^{-\theta
t}(1+|x|)^{-\frac12}
\le CE_0(1+t)^{-3/4}\psi_1(x,t). \\
\end{aligned}
\end{equation}
\end{lemma}

\begin{lemma} [Nonlinear estimates I.] \label{nonlinearintegralestimates}
If $\|v(y, s)\|_{L^\infty} \le C E_0 (1 + s)^{-3/4}$,  
$\|v(y, s)\|_{H^4} \le C E_0 (1 + s)^{-1/2}$,  
and
$\zeta (t)<+\infty$ (see \ref{zetadefined}), then
\begin{equation}\label{Ginteract}
\begin{aligned}
\Big| \int_0^t &\int_{-\infty}^{+\infty} \tilde{G} (x, t-s; y) 
\mathcal{F} (\varphi, v, \frac{\partial \bar{u}^\delta}{\partial \delta} \delta)_y dy ds \Big|
\le
CE_0 \zeta(t) [\psi_1 (x, t) + \psi_2 (x, t)] \\
\Big|\int_0^t &\int_{-\infty}^{+\infty} \tilde{G} (x, t-s; y) \Phi (y, s) dy ds \Big|
\le
CE_0 \psi_1 (x, t) \\ 
\Big|\int_0^{t} &\int_{-\infty}^{+\infty} \tilde{G}_{x} (x, t-s; y) 
\mathcal{F} (\varphi, v, \frac{\partial \bar{u}^\delta}{\partial \delta} \delta)_y dy ds \Big| \\
&\le
CE_0 \zeta (t) \Big[t^{-1/2}(\psi_1 (x, t) + \psi_2 (x, t)) + \psi_3 (x, t) + \psi_4 (x, t)]  \\
\Big|\int_0^{t} &\int_{-\infty}^{+\infty} \tilde{G}_{x} (x, t-s; y) \Phi (y, s) dy ds
\le
CE_0 t^{-1/2} \psi_1 (x, t) \\
\Big| (\partial_t + a_j^\pm \partial_x)  
\int_0^t &\int_{-\infty}^{+\infty} \tilde{G} (x, t-s; y) 
\mathcal{F} (\varphi, v, \frac{\partial \bar{u}^\delta}{\partial \delta} \delta)_y dy ds \Big| \\
&\le
CE_0 \zeta (t) \Big[t^{-1} (1+t)^{1/4} \psi_1^{j, \pm} (x, t) 
+
t^{-1/2}(\bar{\psi}_1^{j, \pm} (x, t) + \psi_2 (x, t)) + \psi_3 (x, t) + \psi_4 (x, t)] \\
\Big|(\partial_t + a_j^\pm \partial_x) 
\int_0^t &\int_{-\infty}^{+\infty} \tilde{G} (x, t-s; y) \Phi (y, s) dy ds \Big| \\
&\le
CE_0 \Big[t^{-1} (1+t)^{1/4} \psi_1^{j,\pm} (x, t) + t^{-1/2} (\bar{\psi}_1^{j,\pm}(x,t) + \psi_2 (x,t)) \Big] \\
\Big|\int_0^t &\int_{-\infty}^{+\infty} \partial_y e_i (y, t-s) 
\mathcal{F} (\varphi, v, \frac{\partial \bar{u}^\delta}{\partial \delta} \delta)_y dy ds \Big|
\le C E_0 \zeta (t) (1 + t)^{-3/4} \\
\Big|\int_0^t &\int_{-\infty}^{+\infty} \partial_{yt} e_i (y, t-s) 
\mathcal{F} (\varphi, v, \frac{\partial \bar{u}^\delta}{\partial \delta} \delta)_y dy ds \Big|
\le C E_0 \zeta (t) (1 + t)^{-1} \\
\Big|\int_0^t &\int_{-\infty}^{+\infty} e_i (y, t-s) \Phi (y, s) dy ds \Big| \le C E_0 (1 + t)^{-1/2} \\
\Big|\int_0^t &\int_{-\infty}^{+\infty} \partial_t e_i (y, t-s) 
\Phi (y, s) dy ds \Big| \le C E_0 (1 + t)^{-1}. \\
\end{aligned}
\end{equation}
\end{lemma}

\begin{lemma}[Nonlinear estimates II.]\label{Hnonlinear} 
If $\|v(y, s)\|_{L^\infty} \le C E_0 (1 + s)^{-3/4}$,  
$\|v(y, s)\|_{H^5} \le C E_0 (1 + s)^{-1/2}$,  
and
$\zeta (t)<+\infty$, then 
\begin{equation} \label{Hinteract}
\begin{aligned}
\Big| \int_0^t &\int_{-\infty}^{+\infty} H (x, t-s; y) 
\mathcal{F} (\varphi, v, \frac{\partial \bar{u}^\delta}{\partial \delta} \delta)_y dy ds \Big|
\le
CE_0 \zeta(t) [\psi_1 (x, t) + \psi_2 (x, t)] \\
\Big|\int_0^t &\int_{-\infty}^{+\infty} H (x, t-s; y) \Phi (y, s) dy ds \Big|
\le
CE_0 \psi_1 (x, t) \\ 
\Big|\int_0^{t} &\int_{-\infty}^{+\infty} H_{x} (x, t-s; y) 
\mathcal{F} (\varphi, v, \frac{\partial \bar{u}^\delta}{\partial \delta} \delta)_y dy ds \Big| \\
&\le
CE_0 \zeta (t) \Big[t^{-1/2}(\psi_1 (x, t) + \psi_2 (x, t)) + \psi_3 (x, t) + \psi_4 (x, t)]  \\
\Big|\int_0^{t} &\int_{-\infty}^{+\infty} H_{x} (x, t-s; y) \Phi (y, s) dy ds
\le
CE_0 t^{-1/2} (\psi_1+ \psi_2) (x, t) \\
\Big| (\partial_t + a_j^\pm \partial_x)  
\int_0^t &\int_{-\infty}^{+\infty} {H} (x, t-s; y) 
\mathcal{F} (\varphi, v, \frac{\partial \bar{u}^\delta}{\partial \delta} \delta)_y dy ds \Big| \\
&\le
CE_0 \zeta (t) \Big[t^{-1} (1+t)^{1/4} \psi_1^{j, \pm} (x, t) 
+
t^{-1/2}(\bar{\psi}_1^{j, \pm} (x, t) + \psi_2 (x, t)) + \psi_3 (x, t) + \psi_4 (x, t)] \\
\Big|(\partial_t + a_j^\pm \partial_x) 
\int_0^t &\int_{-\infty}^{+\infty} H (x, t-s; y) \Phi (y, s) dy ds \Big| \\
&\le
CE_0 \Big[t^{-1} (1+t)^{1/4} \psi_1^{j,\pm} (x, t) + t^{-1/2} (\bar{\psi}_1^{j,\pm}(x,t) + \psi_2 (x,t)) \Big] .\\
\end{aligned}
\end{equation}
\end{lemma}

\medskip
\noindent {\bf Proof of Theorem \ref{pointwiseestimates}.} 
We prove Theorem \ref{pointwiseestimates} directly from the integral
representations (\ref{delta}) and (\ref{tilde{G}}), augmented with 
similar representations for $\dot{\delta}_i (t)$, $v_x (x, t)$, 
and $(\partial_t + a_j^\pm \partial_x) v(x, t)$, obtained through 
direct differentiation of (\ref{delta}) and (\ref{tilde{G}}).  Recalling 
the definition of our iteration variable $\zeta(t)$ in (\ref{zetadefined}),
our goal will be to employ the estimates of Lemmas 
\ref{linearintegralestimates}--\ref{Hnonlinear} to establish the 
inequality
\begin{equation} \label{continuousinduction1}
\zeta(t) \le CE_0 (1 + \zeta (t)),
\end{equation}
where $E_0$ is precisely as in Theorem \ref{pointwiseestimates}.
Choosing, then, $E_0 \le \frac{1}{2C}$, and noting that 
the possibility of $\zeta (t)$ jumping discontinuously from a 
finite value to an infinite value is precluded by 
short-time theory (see Remark 3.6 in \cite{HR})\footnote{
Note also that the pointwise bounds of \cite{HR} 
directly imply that $\zeta$ is finite for each given $t$.
},
we will be able to conclude
\begin{equation} \label{continuousinduction2}
\zeta(t) \le \frac{C E_0}{1 - C E_0} \le 1.
\end{equation}
The estimates of Theorem \ref{pointwiseestimates} follow
immediately from (\ref{continuousinduction2}) and 
(\ref{zetadefined}).  

We begin with a careful consideration of the integral 
represention for $v(t, x)$, (\ref{tilde{G}}), which we
will separate for clarity into terms that arise in the 
case of strict parabolicity---involving 
$\tilde{G}$---and the additional terms arising from the 
relaxation of strict parabolicity---involving $H$.  
For the terms involving $\tilde{G}$, we estimate 
\begin{equation*}
\begin{aligned}
\Big|\int_{-\infty}^{+\infty} \tilde{G} (x, t; y) v_0 (y) dy \Big|
& +
\Big|\int_0^t \int_{-\infty}^{+\infty} \tilde{G} (x, t - s; y) 
\mathcal{F} (\varphi, v, \frac{\partial \bar{u}^\delta}{\partial \delta} \delta)_y (y, s) dy ds \Big| \\
&+
\Big|\int_0^t \int_{-\infty}^{+\infty} \tilde{G} (x, t - s; y) 
\Phi (y, s) dy ds \Big| \\
&\le
CE_0 \psi_1 (x, t) + CE_0 \zeta(t) [\psi_1 (x, t) + \psi_2(x, t)] + CE_0 \psi_1 (x, t),
\end{aligned}
\end{equation*}   
where the estimates follow respectively from the first estimate of Lemma \ref{linearintegralestimates},
the first estimate of Lemma \ref{nonlinearintegralestimates}, and the second estimate of 
Lemma \ref{nonlinearintegralestimates} (though the the first and last estimates are the same, 
we write both for clarity).  For the terms in (\ref{tilde{G}}) involving $H$, we have 
similarly
\begin{equation*}
\begin{aligned}
\Big|\int_{-\infty}^{+\infty} H (x, t; y) v_0 (y) dy \Big|
& +
\Big|\int_0^t \int_{-\infty}^{+\infty} H (x, t - s; y) 
\mathcal{F} (\varphi, v, \frac{\partial \bar{u}^\delta}{\partial \delta} \delta)_y (y, s) dy ds \Big| \\
&+
\Big|\int_0^t \int_{-\infty}^{+\infty} H (x, t - s; y) 
\Phi (y, s) dy ds \Big| \\
&\le
CE_0 \psi_2 (x, t) + CE_0 \zeta(t) [\psi_1 (x, t) + \psi_2(x, t)] + CE_0 \psi_1 (x, t),
\end{aligned}
\end{equation*}
for which the estimates follow respectively from the first estimate of Lemma \ref{Hlinear},
the first estimate of Lemma \ref{Hnonlinear}, and the second estimate of Lemma 
\ref{Hnonlinear}.  Combining these estimates, we conclude 
\begin{equation*}   
|v(x, t)| \le CE_0 [\psi_1 (x, t) + \psi_2 (x, t)]
+
CE_0 \zeta (t) [\psi_1 (x, t) + \psi_2 (x, t)],
\end{equation*}
which can be rearranged as 
\begin{equation*}
\frac{|v(t, x)|}{\psi_1 (x, t) + \psi_2 (x, t)} = C E_0 (1 + \zeta (t)).
\end{equation*}
Keeping in mind that $\zeta (t)$ is a nondecreasing function of $t$, we 
have at last
\begin{equation} \label{quotientbound}
\sup_{y, 0 \le s \le t} \frac{|v(y, s)|}{\psi_1 (y, s) + \psi_2 (y, s)}
\le C E_0 (1 + \zeta (t)).
\end{equation}
Proceeding similarly for each of the expressions $\delta_i (t)$,
$\dot{\delta}_i (t)$, $v_x (x, t)$, 
and $(\partial_t + a_j^\pm \partial_x) v(x, t)$, we can bound
each summand in the definition of $\zeta (t)$ in precisely the 
same way by $C E_0 (1 + \zeta (t))$.  We conclude the sought
estimate
\begin{equation*}
\zeta (t) \le C E_0 (1 + \zeta (t)).
\end{equation*}
As discussed in (\ref{continuousinduction1}) and  (\ref{continuousinduction2}),
we conclude $\zeta (t) \le 1$, and from this the estimates of 
Theorem \ref{pointwiseestimates}.
\hfill $\square$

\section{Liu's cancellation estimate} \label{liu}

%
%
Before carrying out the deferred integral estimates, we revisit
a key estimate of Liu \cite{Liu85} that 
at once determines the ultimate rate of decay and 
motivates the analysis to follow.
Consider the illustrative convolution 
\begin{equation}
\begin{aligned} \label{dovomi}
u(x,t)&= \int_0^t \int_{-\infty}^{+\infty}g(x-y, t-s)(K(y-s,
s)^2)_y \, dy \, ds\\
&= \int_0^t \int_{-\infty}^{+\infty}g_y(x-y, t-s) K(y-s, s)^2 \,
dy \, ds,
\end{aligned}
\end{equation}
$g(x,t)= K(x,t)=(4\pi t)^{-\frac 12}e^\frac{-x^2}{4t}$,
similar to quadratic interaction integrals arising through the
integration of scattering terms against diffusion waves
(see Remark \ref{charbds}).
If we replace the integrands in (\ref{dovomi}) by their absolute value, 
we obtain (see \cite{HZ}) the sharp estimate
\begin{equation}\label{abs}
|u(x,t)| \leq C  \left( g(x, 4t) +
g(x-t, 4t) \right)\\
+ C \chi_{\{\sqrt{t} \leq x \leq t-\sqrt{t}\}} (  x^{-\frac 12}
(t-x)^{-\frac12}).
\end{equation}

By taking account of cancellation, however,
we may obtain the following stronger bound (also sharp)
pointed out by Liu \cite{Liu85}.
We follow the notation of Raoofi \cite{Ra}, and also the proof,
based on integration by parts in the characteristic direction,
which abstracts from the more concrete calculation of Liu the central
ideas that will be of use here.

\begin{proposition}[\cite{Liu85, Ra}]\label{sect1main}
For $u(x,t)$ defined in (\ref{dovomi}) and $t \geq 1$, 
\begin{equation}
\begin{split}
|u(x,t)| &\leq C \trob \left( g(x, 8t) +
g(x-t, 8t) \right)
+ C \chi_{\{\sqrt{t} \leq x \leq t-\sqrt{t}\}} \left(t^{-\frac 12}
(t-x)^{-1}+ x^{- \frac 32} \right),
\end{split}\label{liubounds}
\end{equation}
 where $\chi$ stands for the
indicator function, and $C$ is a constant independent of $t$
and $x$.\\
The same result holds if $K^2$ in (\ref{dovomi}) is replaced with $K_x$.
\end{proposition}

This was first pointed out in
 
\begin{proof}[Proof of Proposition \ref{sect1main}]
As $K^2 \sim K_x$, the proof will be stated only for $K^2$. It
is straightforward to verify that the same argument works
for $K_x$ at every step.
 We first state a simple lemma.
\begin{lemma} \label{stlem}
If \, $0\leq s\leq \sqrt{t}$, \, then $e^{\frac{-(x\pm s)^2}{4t}}
\leq Ce^{\frac{-x^2}{8t}}$ with $C$ independent of $t,s$ and $x$.
\end{lemma}
\begin{proof}
The statement of the lemma is equivalent to
\begin{equation}
\frac{-(x\pm s)^2}{4t} \leq \frac{-x^2}{8t}+D \notag
\end{equation}
for some $D$, which (after some calculation) in its turn is
equivalent to \linebreak $(x\pm 2s)^2-2s^2 \geq -8Dt$, which holds
for $D > \frac 14$, since $s^2 < t$.
\end{proof}
The argument relies on the following simple properties of the
heat kernel $g$.
\begin{align}
 \int_{-\infty}^{+\infty}g(x-y, t)g(y,
t') dy = g(x,t+t') \label{gg1}\\
 |g_x(x,t)|\leq C\tnim g(x,2t)
\label{gg2}\\
 |g_t(x,t)|\leq C\tyek
g(x,2t) \label{gg3}\\ |g(x,t)| \leq C\, \tnim. \label{gg4}
\end{align}
Rewriting (\ref{dovomi}), we have
\begin{equation}
\begin{aligned}
u(x,t)&= \int_0^t \int_{-\infty}^{+\infty}g(x-y, t-s)(g(y-s,s)^2)_y \, dy \, ds\\
&=\int_0^{\sqrt{t}} \int_{-\infty}^{+\infty}g(x-y, t-s)(g(y-s,
s)^2)_y \, dy \, ds\\
&+ \int_{\sqrt{t}}^{t-\sqrt{t}} \int_{-\infty}^{+\infty}g(x-y,
t-s)(g(y-s, s)^2)_y \, dy \, ds\\
&+ \int_{t-\sqrt{t}}^t \int_{-\infty}^{+\infty}g(x-y, t-s)(g(y-s,
s)^2)_y \, dy \, ds\\
&=: I + II + III.
\end{aligned} \label{shekaf}
\end{equation}
$(I)$ and $(III)$ are easy to estimate:
\begin{align}
|I| &= \left| \int_0^{\sqrt{t}} \int_{-\infty}^{+\infty}g(x-y,
t-s)(g(y-s, s)^2)_y \, dy \, ds\right| \notag\\
&= \left| \int_0^{\sqrt{t}} \int_{-\infty}^{+\infty}g_y(x-y, t-s)
g(y-s, s)^2 \, dy \, ds\right|.\notag
\end{align}
By (\ref{gg2}) and (\ref{gg4}), the above is less than or equal to
$$   C\int_0^{\sqrt{t}}
\int_{-\infty}^{+\infty}(t-s)^{-\frac12} s^{-\frac12} g(x-y,
2(t-s))\,g(y-s, 2s) \, dy \, ds $$
which, by (\ref{gg1}), is less than or equal to
$
  C\int_0^{\sqrt{t}}
(t-s)^{-\frac12} s^{-\frac12} g(x-s, 2t) \, ds. \notag
$
Now, using Lemma \ref{stlem}, the above is
\begin{equation}
\begin{aligned}
 \leq  C\,g(x, 4t) \int_0^{\sqrt{t}} 
(t-s)^{-\frac12} s^{-\frac12} \, ds &\leq  C\tnim g(x,
4t) \int_0^{\sqrt{t}} s^{-\frac12} \, ds\\ 
 &\leq C\trob g(x,4t).
\end{aligned}
\end{equation}
 Part $(III)$
in (\ref{shekaf}) can be handled similarly.

The more difficult part is part $(II)$ of (\ref{shekaf}):
\begin{equation}
II= \int_{\sqrt{t}}^{t-\sqrt{t}} \int_{-\infty}^{+\infty}g(x-y,
t-s)(g(y-s, s)^2)_y \, dy \, ds.
\end{equation}
In order to estimate $(II)$, let us write $g=g(x, \tau)$. Then we
have
\begin{equation}
\begin{aligned}
g(x-y, t-s)(g(y-s, s)^2)_y \,&= (g(x-y, t-s)g(y-s, s)^2)_s \\ & -
\,g_\tau (x-y, t-s)g(y-s, s)^2\\ &+\, g(x-y, t-s)(g^2)_\tau(y-s,
s) .
\end{aligned} \label{shekaf2}
\end{equation}
We will do the estimates piece by piece.
The first part of (\ref{shekaf2}) can be decomposed as
\begin{align}
\int_{\sqrt{t}}^{t-\sqrt{t}} \int_{-\infty}^{+\infty} (g(x-y,
t-s)g(y-s, s)^2)_s\, dy\, ds\\
= \int_{-\infty}^{+\infty} g(x-y, \sqrt{t})\,g(y-t+\sqrt{t},
t-\sqrt{t})^2\,dy\\
- \int_{-\infty}^{+\infty} g(x-y, t-\sqrt{t})\,g(y-\sqrt{t},
\sqrt{t})^2\,dy.
\end{align}
Using (\ref{gg1}) and (\ref{gg4}), it follows that
$$ \int_{-\infty}^{+\infty} g(x-y, \sqrt{t})\,g(y-t+\sqrt{t},
t-\sqrt{t})^2\,dy \leq C\tnim g(x-t+\sqrt{t}, t), $$ and
$$\int_{-\infty}^{+\infty} g(x-y, t-\sqrt{t})\,g(y-\sqrt{t},
\sqrt{t})^2\,dy \leq C\trob g(x-\sqrt{t}, t),$$ but, by Lemma
\ref{stlem},
$g(x-\sqrt{t}, t) \leq g(x, 2t)$ and
$g(x-t+\sqrt{t}, t) \leq g(x-t, 2t).$
These terms fit in the right hand side of (\ref{liubounds}).

For the other parts in (\ref{shekaf2}), we use (\ref{gg3}) to obtain
\begin{equation}
\begin{aligned}
&\left| \int_{\sqrt{t}}^{t-\sqrt{t}} \int_{-\infty}^{+\infty}
g_\tau (x-y, t-s) g(y-s, s)^2 \, dy\, ds
\right|  \\
\leq &\int_{\sqrt{t}}^{t-\sqrt{t}} s^{-\frac 12} (t-s)^{-1} g(x-s,
2t) ds  \\
=&\int_{\sqrt{t}}^{\frac t2} s^{-\frac 12} (t-s)^{-1} g(x-s, 2t)
ds \\
+ &\, \int_{\frac t2}^{t-\sqrt{t}} s^{-\frac 12} (t-s)^{-1}
g(x-s, 2t) ds \\
=: \,&\mathcal{A} \,+\, \mathcal{B}.
\end{aligned} \label{AB}
\end{equation}

If $x \leq \sqrt{t}$, then
\begin{equation}
\begin{aligned}
\mathcal{A} &\leq C\,\tyek g(x-\sqrt{t}, 2t)\int_{\sqrt{t}}^{\frac t2} s^{-\frac 12}\, ds\\
&\leq  C \, t^{-\frac 12} \, g(x-\sqrt{t}, 2t)\\
&\leq C \tnim g(x, 4t).
\end{aligned}
\end{equation}
Similarly when $x\ge t-\sqrt{t}$, we obtain $\mathcal A \le C\tnim
g(x-t, 4t).$

Now for $\sqrt{t} \leq x \leq t-\sqrt t$ we have
\begin{equation}
\begin{aligned}
\mathcal{A} &= \int_{\sqrt{t}}^{\frac t2} s^{-\frac 12}
(t-s)^{-1} g(x-s, 2t) ds \\
&\leq C t^{-1} \int_{\sqrt{t}}^{\frac t2} s^{-\frac 12}
 g(x-s, 2t) ds \\
&= C t^{-1} \int_{\sqrt{t}}^{\frac x2} s^{-\frac 12}
 g(x-s, 2t) ds
 +  t^{-1} \int_{\frac x2}^{\frac t2} s^{-\frac 12}
 g(x-s, 2t) ds \\
 &\le C t^{-1} g(\frac x2, 2t)\int_{\sqrt{t}}^{\frac x2} s^{-\frac
 12} ds + t^{-1} x^{-\frac 12} \int_{\frac x2}^{\frac t2}
 g(x-s, 2t) ds\\
 &\le C\left( t^{-\frac 12} g(x, 4t) + t^{-1} x^{-\frac 12}\right)
\end{aligned}
\end{equation}
also acceptable, as $ t^{-1} x^{-\frac 12}\le x^{-\frac32}$ for
$\sqrt t \le x \le t-\sqrt t.$
%

Part $\mathcal{B} $ in (\ref{AB}) can be estimated similarly. We
carry it out briefly  only for   $\frac t2\le x\le t-\sqrt t.$ Let
$\xi:= t-\sqrt t -x.$ Then
\begin{equation}\begin{aligned}
\mathcal{B} &\le C t^{-\frac 12}\left(\int_{\frac
t2}^{x+\xi}+\int_{x+\xi}^{t-\sqrt t}\right)(t-s)^{-1} g(x-s, 2t)ds\\
&\le C t^{-\frac 12}(t-x)^{-1} + C t^{-\frac 14}g(\frac {t-\sqrt
t-x}2, 2t),
\end{aligned}\end{equation}
to which we apply Lemma \ref{stlem}.

There remains the last part of (\ref{shekaf2}), i.e.,
\begin{equation}\left|
\int_{\sqrt{t}}^{t-\sqrt{t}}\int_{-\infty}^{+\infty} g(x-y,
t-s)(g^2)_\tau(y-s, s) \, dy\,
ds\right|,\label{lastone}\end{equation} which can easily be shown to
be less than or equal to
\begin{equation}
C \int_{\sqrt{t}}^{t-\sqrt{t}}s^{-\frac 32} g(x-s, 2t)\, ds.
\label{aan}
\end{equation}

If $x\,\leq\, \sqrt{t}$, then (\ref{aan}) is of the order
$$ C\, g(x-\sqrt{t}, 2t)\int_{\sqrt{t}}^{t-\sqrt{t}}s^{-\frac 32}\, ds$$
$$\leq C\,\trob \, g(x-\sqrt{t}, 2t)$$
$$\leq C\,\trob \, g(x,4t). $$
For $x\, \geq \, t-\sqrt{t}$ we use a similar method.\\

For $\sqrt{t} \, \leq \,x\, \leq \, t-\sqrt{t}$, we use a similar
method to those used for the previous cases to get  $$(\ref{aan})
\, \leq \,
C \left(\trob g(x,4t) + x^{-\frac 32}\right).$$\\
 This completes the proof.
\end{proof}

\begin{remark}\label{charbds}
\textup{
Replacing $K(y-at, t)$ by $\kappa(y,t)$ in (\ref{dovomi}),
and $g(x-y, t)$ by $\gcal(x,t;y)$ we obtain bounds
similar to (\ref{liubounds})
(with appropriate modifications due to different speeds) 
provided $\kappa$ and $\gcal$ satisfy 
\begin{align}
|\gcal(x,t;y)| \leq C g(x-y-at, \beta t),\\
|\gcal_y(x,t;y)| \leq C \tnim g(x-y-at, 2\beta t),\\
|\gcal_t(x,t;y)| \leq C \tyek g(x-y-at, 2\beta t),\\
|\kappa(x,t)| \leq C g(x-bt, \beta t), \label{bd3}\\
|\kappa_y(x,t)| \leq C \tnim g(x-bt,2\beta t),\\
|\kappa_t(x,t;y)| \leq C \tyek g(x-bt, 2\beta t), \label{bd6}
\end{align}
for some $a\neq b,$ and some constants $C, \beta >0$, as hold in particular
for $\kappa$ a single diffusion wave and $\gcal$ a single
component of the scattering term $S$ in the decomposition $G=E+H+S+R$.
Comparing with (\ref{1stestimate}), we see that 
{these are the rate-determing terms}. }

\textup{
More generally \cite{ZH}, we may expect analogous cancellation whenever
$\kappa$ and $\gcal$ have the property that they decay
faster along characteristic directions $\partial_s+a_j^\pm \partial_y$
[resp.  $\partial_t+a_j^\pm \partial_x$] than $\partial_y$ or
$\partial_s$ [resp.  $\partial_x$ or $\partial_t$].
{This is the underlying principle guiding our analysis, and that
of \cite{Liu97}.}
}
\end{remark}

\begin{remark}\label{refinedcharbds}
\textup{
Reviewing (\ref{liubounds}) more carefully, we see that
the signal $u$ is slightly more concentrated near the
characteristic direction $dx/dt=0$ of the propagator $g$
than in the direction $dx/dt=1$ of the source $K^2_y$,
decaying as
$(x+ t^{\frac 12})^{- \frac 32} $
rather than
$t^{-\frac 12}(|t-x|+t^{\frac 14})^{-1}$.
Systematic, characteristic-by-characteristic bookkeeping 
taking account of this difference leads to the 
slightly refined bounds of \cite{Liu97}; these sum, of course,
to the simpler modulus bounds presented here.
}
\end{remark}

\section{Linear integral Estimates}\label{integralestimatesL}

It remains to establish the deferred integral estimates used
in Section \ref{nonlinear}.
We begin, in this section, with
the linear integral estimates of Lemmas
\ref{linearintegralestimates} and \ref{Hlinear}.

\medskip
\noindent
{\bf Proof of Lemma \ref{linearintegralestimates}.}  The first, fourth, and fifth 
estimates of Lemma \ref{linearintegralestimates} have been established in \cite{HR}.  The 
second and third follow in a similar fashion from the estimates of Proposition
\ref{altGFbounds}.  \hfill $\square$

\noindent{\bf Proof of Lemma \ref{linearintegralestimates}.}  
Looking at (\ref{multH}), we notice that in order to estimate
$\int_{-\infty}^{+\infty}H(x,t;y)v_0(y) dy$ it suffices to estimate
\begin{equation}\label{randomcalc}
\begin{aligned}
\int_{-\infty}^{+\infty}\mathcal{R}_j^*(x) \mathcal{O}(e^{-\eta_0 t}) \delta_{x-\bar a_j^*
t}(-y) \mathcal{L}_j^{*t}(y)v_0(y) dy&
\le CE_0 e^{-\eta_0 t}v_0(\bar a_j^*t-x)\\
&\le CE_0 e^{-\eta_0 t} (1+|\bar a_j^*t-x|)^{-\frac32} \\
&\le CE_0 e^{-\eta_0 t} (1+|x|)^{-\frac32}(1+|\bar a_j^*t|)^{\frac32}\\
&\le CE_0 e^{-\frac{\eta_0 t}2} (1+|x|)^{-\frac32}.\\
\end{aligned}
\end{equation}
Here we used the crude inequality
\begin{equation}\label{abineq}
\frac {1}{1+|a+b|}\le \frac{1+|b|}{1+|a|}
\end{equation}
and the fact that $\bar a_j^*$ $\mathcal{R}_j^*$ and
$\mathcal{L}_j^{*t}$ are bounded. 
Observing that $e^{-
\frac {eta_0t}{4}}\le C(N)(1+t)^{-N}$ for any $N$, and
$(1+|x-at|)\le C(1+t)(1+|x|)$, we can bound the righthand side
of (\ref{randomcalc}) in turn by $CE_0 e^{-\frac{\eta_0 t}4} \psi_1(x,t)$,
giving (\ref{Hinit})(i).
Estimates (\ref{Hinit})(ii) and (\ref{Hinit})(iii) are obtained similarly.
 \hfill $\square$

\hfill $\square$

\section{Nonlinear integral Estimates I}\label{integralestimatesNLI}

In this section, we carry out the main work of the paper,
establishing the nonlinear integral estimates of 
Proposition \ref{nonlinearintegralestimates}.

\medskip
\noindent
{\bf Proof of Lemma \ref{nonlinearintegralestimates}.}
Under the assumption that $\zeta (t)$ is bounded, we have the estimates 
\begin{equation*}
\begin{aligned}
|v (x,t)| &\le \zeta(t) \Big[\psi_1 + \psi_2 \Big] (x, t), \\
|\partial_x v (x,t)| &\le \zeta(t) \Big[ t^{-1/2} (\psi_1+\psi_2) + \psi_3 + \psi_4 \Big] (x,t) \\  
|(\partial_t + a_j^\pm \partial_x) v(x,t)|
&\le 
\zeta(t) \Big[t^{-1} (1+t)^{1/4} \psi_1^{j,\pm} + t^{-1/2} (\bar{\psi}_1^{j,\pm}
+ \psi_2) + \psi_3 + \psi_4 \Big] (x, t).
\end{aligned}
\end{equation*}
In the analysis that follows, we will omit $\zeta(t)$ from our estimates
in most cases and focus only on the terms with a given form, the 
{\it template}.
   
We observe at the outset that in proving estimates of form $\psi_1 (x ,t)$,
we will frequently make use of the inequality
\begin{equation} \label{kerneltopsi1}
(1+t)^{-3/4} e^{-\frac{(x - a_j^\pm t)^2}{Lt}}
\le
C (1 + |x - a_j^\pm t| + t^{-1/2})^{-3/2}. 
\end{equation}
In the case $|x - a_k^\pm t| \le K \sqrt{t}$, for some constant $K$, 
this inequality is immediate.  On the other hand, for  
$|x - a_k^\pm t| \ge K \sqrt{t}$, we observe
\begin{equation*}
\begin{aligned}
(1+t)^{-3/4} e^{-\frac{(x - a_j^\pm t)^2}{Lt}}
&=
(1+t)^{-3/4} |x - a_j^\pm t|^{-3/2} t^{3/4} 
\frac{|x - a_j^\pm t|^{3/2}}{t^{3/4}} e^{-\frac{(x - a_j^\pm t)^2}{Lt}} \\
&\le
C_1 (1+t)^{-3/4} t^{3/4} |x - a_j^\pm t|^{-3/2}
e^{-\frac{(x - a_j^\pm t)^2}{2Lt}}
\le
C (1 + |x - a_j^\pm t| + t^{1/2})^{-3/2},
\end{aligned}
\end{equation*}
where the seeming blow-up as $|x - a_j^\pm t| \to 0$ is controlled by 
the size of $\sqrt{t}$, which must be smaller than $|x - a_j^\pm t|$.

\medskip
{\bf Proof of (\ref{Ginteract}(i)).}
For the first estimate in Lemma \ref{nonlinearintegralestimates},
we begin by considering the nonlinearities $(v(y,s)^2)_y$ and 
$(v(y,s)v_y(y,s))_y$.  
Here and in the remaining cases, the 
analyses of the convection, reflection, and transmission 
contributions to the Green's kernel $\partial_y \tilde{G} (x, t; y)$
are similar, and we provide full details only for the case 
of convection.  We observe that the contribution 
\begin{equation*}
\sum_{j=1}^J\sum_{\beta=0}^{\max \{0, |\alpha|-1\}} \mathbf{O}(e^{-\eta t})\partial_y^\beta \delta_{x-\bar
a_j^* t}(-y)  
\end{equation*}
is similar to, though less singular than, terms arising in 
$\partial_y H(x, t; y)$ (see \ref{multH}), and can be analyzed
as in the proof of Lemma \ref{Hnonlinear}.  Finally, we remark that
the contribution 
\begin{equation*} 
\mathbf{O}(e^{-\eta(|x-y|+t)})
\end{equation*}
has no effect on the iteration.

{\it (\ref{Ginteract}(i)), term one.}
We first consider integration of our 
convecting Green's kernel against 
the nonlinearity $s^{-1/2} (1 + s)^{-1/4} \psi_1 (y, s)$.  In this case, we 
have integrals of the form 
\begin{equation} \label{gf1nl1}
\int_0^t \int_{-\infty}^0 (t-s)^{-1} e^{-\frac{(x - y - a_j^- (t-s))^2}{M(t-s)}}
(1 + |y - a_k^- s| + s^{1/2})^{-3/2} s^{-1/2} (1+s)^{-1/4} dy ds,
\end{equation}
with a similar integral for $y \ge 0$.
Proceeding as in \cite{HR}, we write 
\begin{equation} \label{maindecomp}
x - y - a_j^- (t - s)
= (x - a_j^- (t-s) -a_k^- s) - (y - a_k^- s),
\end{equation}
from which we have the estimate
\begin{equation} \label{gf1nl1balance1}
\begin{aligned}
&e^{-\frac{(x - y - a_j^- (t-s))^2}{M(t-s)}} (1 + |y - a_k^- s| + s^{1/2})^{-\gamma} \\
&\le
C\Big[ e^{-\epsilon \frac{(x - y - a_j^- (t-s))^2}{M(t-s)}}
e^{-\frac{(x - a_j^- (t-s) - a_k^- s)^2}{\bar{M}(t-s)}}
(1 + |y - a_k^- s| + s^{1/2})^{-\gamma} \\
&+
e^{-\frac{(x - y - a_j^- (t-s))^2}{M(t-s)}} 
(1 + |y - a_k^- s| + |x - a_j^- (t-s) - a_k^- s| + s^{1/2})^{-\gamma} \Big],
\end{aligned}
\end{equation}
where $\gamma = 3/2$.  (Here and in future cases, we will state estimates 
in terms of a parameter $\gamma$ for general applicability.)
For the first estimate in (\ref{gf1nl1balance1}), we have integrals
\begin{equation} \label{gf1nl1balance1first}
\int_0^t \int_{-\infty}^0 (t-s)^{-1}
e^{-\epsilon \frac{(x - y - a_j^- (t-s))^2}{M(t-s)}}
e^{-\frac{(x - a_j^- (t-s) - a_k^- s)^2}{\bar{M}(t-s)}}
(1 + |y - a_k^- s| + s^{1/2})^{-3/2} s^{-1/2} (1 + s)^{-1/4} dy ds.
\end{equation}
We have three cases to consider: $a_k^- < 0 < a_j^-$, $a_k^- \le a_j^- < 0$,
and $a_j^- < a_k^- < 0$.  We will proceed by analyzing the second of these
in detail and observing that no qualitatively new calculations arise 
in the remaining two.  For this second case and for $|x| \ge |a_k^-| t$, we write
\begin{equation} \label{tminussdecomp}
x - a_j^- (t-s) - a_k^- s = (x - a_k^- t) - (a_j^- - a_k^-) (t-s),
\end{equation}
for which we observe that there is no cancellation between 
summands.  We immediately obtain an estimate on (\ref{gf1nl1balance1first})
by 
\begin{equation} \label{gf1nl1nocanc1}
\begin{aligned}
C_1 &t^{-1} e^{-\frac{(x - a_k^- t)^2}{Lt}} 
\int_0^{t/2} (1 + s^{1/2})^{-1/2} s^{-1/2} (1 + s)^{-1/4} ds \\
&+
C_2 (1 + t^{1/2})^{-3/2} t^{-1/2} (1+t)^{-1/4} e^{-\frac{(x - a_k^- t)^2}{Lt}} 
\int_{t/2}^t (t-s)^{-1/2} ds \\
&\le
C (1+t)^{-1} \ln (e+t) e^{-\frac{(x - a_k^- t)^2}{Lt}},
\end{aligned}
\end{equation}
which is sufficient by (\ref{kerneltopsi1}).
For $|x| \le |a_j^-| t$, we write 
\begin{equation} \label{sdecomp}
x - a_j^- (t-s) - a_k^- s = (x - a_j^- t) - (a_k^- - a_j^-) s,
\end{equation}
for which we observe that there is no cancellation between summands,
and proceeding as in (\ref{gf1nl1nocanc1}), we obtain an estimate by 
\begin{equation*}
C (1+t)^{-1} \ln (e+t) e^{-\frac{(x - a_j^- t)^2}{Lt}},
\end{equation*}
which again is sufficient.  We remark that this concludes the 
analysis in the case $a_j^- = a_k^-$, so from here on we may take 
$a_k^- < a_j^-$ and the case $|a_j^-| t \le |x| \le |a_k^-| t$.
In this case, and for $s \in [0, t/2]$, we observe through (\ref{sdecomp}) the 
inequality
\begin{equation} \label{gf1nl1balance2} 
\begin{aligned}
&e^{-\frac{(x - a_j^- (t-s) - a_k^- s)^2}{\bar{M}(t-s)}}
(1 + |y - a_k^- s| + s^{1/2})^{-3/2} (1 + s)^{-\gamma} \\
&\le 
C \Big[ e^{-\frac{(x - a_j^- t)^2}{Lt}}
(1 + |y - a_k^- s| + s^{1/2})^{-3/2} (1 + s)^{-\gamma} \\
&+
e^{-\frac{(x - a_j^- (t-s) - a_k^- s)^2}{\bar{M}(t-s)}}
(1 + |y - a_k^- s| + |x - a_j^- t|^{1/2})^{-3/2} 
(1 + |x-a_j^- t|)^{-\gamma} \Big],
\end{aligned}
\end{equation}
with $\gamma = 3/4$.
For the first estimate in (\ref{gf1nl1balance2}), we proceed similarly
as in (\ref{gf1nl1nocanc1}), while for the second, we have, upon integration 
of the $\epsilon$-kernel in $y$, an estimate by 
\begin{equation*}
\begin{aligned}
C_1 &t^{-1/2} (1+|x - a_j^- t|^{1/2})^{-3/2}  (1 + |x - a_j^- t|)^{-3/4}
\int_0^{t/2} s^{-1/2} (1+s)^{1/2} e^{-\frac{(x - a_j^- (t-s) - a_k^- s)^2}{\bar{M}(t-s)}} ds \\
&\le
C (1 + |x - a_j^- t|)^{-3/2},
\end{aligned}
\end{equation*} 
where we have used in this last inequality that $a_j^- \ne a_k^-$.
We observe that this is clearly sufficient in the case $|x - a_j^- t| \ge \sqrt{t}$,
whereas in the case $|x - a_j^- t| \le \sqrt{t}$, decay of kernel type
$\exp((x-a_j^-t)^2/(Lt))$ is immediate, for $L$ sufficiently large, and so we only require decay at rate
$t^{-3/4}$, which is straightforward.
For $s \in [t/2, t]$, we observe through (\ref{tminussdecomp}) the inequality 
\begin{equation} \label{gf1nl1balance3}
\begin{aligned}
(t-s)^{-\gamma} & e^{-\frac{(x - a_j^- (t-s) - a_k^- s)^2}{\bar{M}(t-s)}} \\
&\le
C \Big[ |x - a_k^- t|^{-\gamma} e^{-\frac{(x - a_j^- (t-s) - a_k^- s)^2}{\bar{M}(t-s)}}
+ (t - s)^{-\gamma} e^{-\frac{(x - a_j^- t)^2}{Lt}} \Big],
\end{aligned}
\end{equation}
where $\gamma = 1$.
For the second estimate in (\ref{gf1nl1balance3}), we proceed similarly as in 
(\ref{gf1nl1nocanc1}), while for the first we have, upon integration in $y$
of the nonlinearity, an estimate by 
\begin{equation*}
\begin{aligned}
C_2 &(1+t^{1/2})^{-1/2} t^{-1/2} (1+t)^{-1/4} |x - a_k^- t|^{-1}
\int_{t/2}^t e^{-\frac{(x - a_j^- (t-s) - a_k^- s)^2}{\bar{M}(t-s)}} ds \\
&\le
C (1 + t)^{-1/2} |x - a_k^- t|^{-1}.
\end{aligned}
\end{equation*}
Recalling that we are currently working in the case $t \ge |x|/|a_k^-|$, we see that
this final estimate is sufficient for the case $|x - a_k^- t| \ge \sqrt{t}$.
For the second term in (\ref{gf1nl1balance1}), we have integrals of the form 
\begin{equation} \label{gf1nl1balance1second}
\int_0^t \int_{-\infty}^0 (t - s)^{-1} e^{-\frac{(x - y - a_j^- (t-s))^2}{M(t-s)}} 
(1 + |y - a_k^- s| + |x - a_j^- (t-s) - a_k^- s| + s^{1/2})^{-3/2} 
s^{-1/2} (1+s)^{-1/4} ds.
\end{equation}
We continue to focus on the case $a_k^- \le a_j^-$.  For the cases 
$|x| \ge |a_k^-| t$ and $|x| \le |a_j^-| t$, respectively, we have no 
cancellation between the summands in (\ref{tminussdecomp}) and 
(\ref{sdecomp}) and consequently we obtain estimates
\begin{equation} \label{gf1nl1nocanc2}
C (1+t)^{-1/4} \Big[(1+|x - a_k^- t|)^{-3/2} + (1+|x - a_j^- t|)^{-3/2} \Big].
\end{equation}
For the case $|a_j^-| t \le |x| \le |a_k^-| t$, we divide the analysis into
sub-intervals $s \in [0, t/2]$ and $s \in [t/2, t]$.  For $s \in [0, t/2]$,
we observe through (\ref{sdecomp}) the inequality
\begin{equation} \label{gf1nl1balance4}
\begin{aligned}
&\Big(1 + |y - a_k^- s| + |x - a_j^- (t-s) - a_k^- s| + s^{1/2}\Big)^{-3/2} 
(1+s)^{-\gamma} \\
&\le
C\Big[ (1 + |x - a_j^- t| + s^{1/2})^{-3/2} (1+s)^{-\gamma} \\
&+
(1 + |y - a_k^- s| + |x - a_j^- (t-s) - a_k^- s| + |x - a_j^- t|^{1/2})^{-3/2} 
(1+|x-a_j^- t|)^{-\gamma} \Big],
\end{aligned}
\end{equation}
for $\gamma = 3/4$.
For the first estimate in (\ref{gf1nl1balance4}), we proceed as in 
(\ref{gf1nl1nocanc2}), while for the second we have, upon integration of
the kernel, an estimate on (\ref{gf1nl1balance1second}) by 
\begin{equation*}
\begin{aligned}
C_1 &t^{-1/2} (1+|x-a_j^- t|)^{-3/4} 
\int_0^{t/2} s^{-1/2} (1+s)^{1/2} 
(1 + |x - a_j^- (t-s) - a_k^- s| + |x - a_j^- t|^{1/2})^{-3/2} ds \\
&\le
C (1+t)^{-1/2} (1+|x-a_j^- t|)^{-1}.
\end{aligned}
\end{equation*}
For $s \in [t/2, t]$, we observe through (\ref{tminussdecomp}) the inequality 
\begin{equation} \label{gf1nl1balance5}
\begin{aligned}
(t - s)^{-1} &(1 + |x - a_j^- (t-s) - a_k^- s| + s^{1/2})^{-3/2} \\
&\le
C \Big[|x - a_k^- t|^{-1} (1 + |x - a_j^- (t-s) - a_k^- s| + s^{1/2})^{-3/2} \\
&+
(t - s)^{-1} (1 + |x - a_j^- (t-s) - a_k^- s| + |x - a_j^- t| + s^{1/2})^{-3/2} \Big].
\end{aligned}
\end{equation}
For the second estimate in (\ref{gf1nl1balance5}), we proceed similarly 
as in (\ref{gf1nl1nocanc2}), while for the first we have, upon 
integration of the kernel, an estimate by 
\begin{equation*}
\begin{aligned}
C_2 &|x - a_k^- t|^{-1} t^{-1/2} (1+t)^{-1/4} 
\int_{t/2}^t (t-s)^{1/2} (1 + |x - a_j^- (t-s) - a_k^- s| + t^{1/2})^{-3/2} ds \\
&\le
C |x - a_k^- t|^{-1} (1+t)^{-1/4} (1+t^{1/2})^{-1/2},
\end{aligned}
\end{equation*}
which is sufficient for $t \ge |x|/|a_j^-|$.


{\it (\ref{Ginteract}(i)), term two.}
We next consider integration against the nonlinearity $s^{-1/2} (1+s)^{-1/4} \psi_2 (y, t)$,
for which we have integrals of the form 
\begin{equation} \label{gf1nl2}
\int_0^t \int_{-|a_1^-| s}^0 (t-s)^{-1} e^{-\frac{(x - y - a_j^- (t-s))^2}{M(t-s)}}
(1+|y|)^{-1/2} (|y| + s)^{-1/2} (1 + |y| + s)^{-3/4} (1 + |y| + s^{1/2})^{-1/2} dy ds,
\end{equation}
wherein we have observed that for $y \in [-|a_1^-| s, 0]$, $s$ decay yields also
decay in $y$.   We first observe an immediate time decay estimate by
\begin{equation} \label{timedecayestimate2}
\begin{aligned}
C_1 &t^{-1} \int_0^{t/2} (1+s)^{-3/4} (1+s^{1/2})^{-1/2} ds
+
C_2 t^{-1/2} (1+t)^{-3/4} (1+t^{1/2})^{-1/2} \int_{t/2}^t (t-s)^{-1/2} ds \\
&\le
C (1 + t)^{-1} \ln (e+t).
\end{aligned}
\end{equation}
In order to determine estimates in space as well, we observe the inequality
\begin{equation} \label{gf1nl2balance1}
\begin{aligned}
&e^{-\frac{(x - y - a_j^- (t-s))^2}{M(t-s)}}
(1+|y|)^{-1/2} (|y|+s)^{-1/2} (1 + |y| + s)^{-3/4} (1 + |y| + s^{1/2})^{-1/2} \\
&\le
C \Big[e^{-\epsilon \frac{(x - y - a_j^- (t-s))^2}{M(t-s)}}
e^{-\frac{(x - a_j^- (t-s))^2}{\bar{M}(t-s)}}
(1+|y|)^{-1/2} (|y|+s)^{-1/2} (1 + |y| + s)^{-3/4} (1 + |y| + s^{1/2})^{-1/2} \\
&+
e^{-\frac{(x - y - a_j^- (t-s))^2}{M(t-s)}}
(1+|x-a_j^- (t-s)|)^{-1/2} (|x-a_j^- (t-s)| + s)^{-1/2} \\
&\times
(1 + |x-a_j^- (t-s)| + s)^{-3/4} 
(1 + |x-a_j^- (t-s)| + s^{1/2})^{-1/2} \Big].
\end{aligned}
\end{equation}
For the first estimate in (\ref{gf1nl2balance1}), we have integrals
\begin{equation} \label{gf1nl2balance1first}
\begin{aligned}
\int_0^t &\int_{-|a_1^-| s}^0 (t-s)^{-1} e^{-\epsilon \frac{(x - y - a_j^- (t-s))^2}{M(t-s)}}
e^{-\frac{(x - a_j^- (t-s))^2}{\bar{M} (t-s)}} \\
&\times (1+|y|)^{-1/2} (|y| + s)^{-1/2} (1 + |y| + s)^{-3/4} (1 + |y| + s^{1/2})^{-1/2} dy ds.
\end{aligned}
\end{equation}
We have two cases to consider here, $a_j^- < 0$ and $a_j^- > 0$, of which we
focus on the former.  In this case, for $|x| \ge |a_j^-| t$, there is no
cancellation between $x - a_j^- t$ and $a_j^- s$, and so proceeding as in 
(\ref{timedecayestimate2}), we obtain an estimate by 
\begin{equation*}
C (1+t)^{-1} \ln (e+t) e^{-\frac{(x - a_j^- t)^2}{Lt}}.
\end{equation*}
For $|x| \le |a_j^-| t$, we divide the analysis into two cases, 
$s \in [0, t/2]$ and $s \in [t/2, t]$.  For $s \in [0, t/2]$, we
observe the inequality
\begin{equation} \label{gf1nl2balance2}
\begin{aligned}
&e^{-\frac{(x - a_j^- (t-s))^2}{\bar{M} (t-s)}}
(|y| + s)^{-1/2} (1 + |y| + s)^{-3/4} (1 + |y| + s^{1/2})^{-1/2} \\
&\le
C \Big[e^{-\frac{(x - a_j^- t)^2}{L t}}
(|y|+s)^{-1/2} (1 + |y| + s)^{-3/4} (1 + |y| + s^{1/2})^{-1/2} \\
&+
e^{-\frac{(x - a_j^- (t-s))^2}{\bar{M} (t-s)}} (s + |x - a_j^- t|)^{-1/2}
(1 + s + |x - a_j^- t|)^{-3/4} (1 + s^{1/2} + |x - a_j^- t|^{1/2})^{-1/2}
\Big].
\end{aligned}
\end{equation}
For the first estimate in (\ref{gf1nl2balance2}), we proceed similarly 
as in (\ref{timedecayestimate2}), while for the second we have, upon 
integration of the $\epsilon$ kernel, an estimate by 
\begin{equation*}
\begin{aligned}
C_1 &t^{-1/2} (|x - a_j^- t|)^{-1/2}  (1 + |x - a_j^- t|)^{-1} 
\int_0^{t/2} e^{-\frac{(x - a_j^- (t-s))^2}{\bar{M} (t-s)}} ds \\
&\le
C (|x - a_j^- t|)^{-1/2} (1 + |x - a_j^- t|)^{-1},
\end{aligned}
\end{equation*}
where the seeming blow-up as $|x - a_j^- t| \to 0$ can be eliminated
by proceeding alternatively for $|x - a_j^- t|$ bounded.
For $s \in [t/2, t]$, we have
\begin{equation} \label{gf1nl2balance3}
(t-s)^{-\gamma} e^{-\frac{(x - a_j^- (t-s))^2}{\bar{M} (t-s)}}
\le 
C |x|^{-\gamma} e^{-\frac{(x - a_j^- (t-s))^2}{2 \bar{M} (t-s)}},
\end{equation}
$\gamma = 1/2$, for which we have an estimate by 
\begin{equation*}
C_2 |x|^{-1/2} t^{-1/2} (1+t)^{-3/4} (1+t^{1/2})^{-1/2} 
\int_{t/2}^t e^{-\frac{(x - a_j^- (t-s))^2}{2 \bar{M} (t-s)}} ds
\le
C |x|^{-1/2} (1+t)^{-1}, 
\end{equation*}
where the seeming blow-up as $|x| \to 0$ can be eliminated by an 
alternative calculation in the case of $|x|$ bounded.  In the 
current case of $|x| \le |a_j^-| t$, this final estimate is bounded
by $\psi_2 (x, t)$.
For the second estimate in (\ref{gf1nl2balance1}), we have integrals
\begin{equation} \label{gf1nl2balance1second}
\begin{aligned}
\int_0^t &\int_{-|a_1^-| s}^0 (t-s)^{-1} e^{-\frac{(x - y - a_j^- (t-s))^2}{M(t-s)}} 
(1+|x - a_j^- (t-s)|)^{-1/2} \\
&\times (|x - a_j^- (t-s)| + s)^{-1/2}  (1 + |x - a_j^- (t-s)| + s)^{-3/4} 
(1 + |x - a_j^- (t-s)| + s^{1/2})^{-1/2} dy ds.
\end{aligned}
\end{equation}
Focusing again on the case $a_j^- < 0$, we observe that for $|x| \ge |a_j^-| t$,
there is no cancellation between $x - a_j^- t$ and $a_j^- s$, and consequently
that we have an estimate, upon integration of the kernel, by 
\begin{equation} \label{gf1nl2nocanc2}
\begin{aligned}
C_1 &t^{-1/2} (1 + |x - a_j^- t|)^{-7/4} 
\int_0^{t/2} (|x - a_j^- (t-s)| + s)^{-1/2} ds \\
&+ 
C_2 (1 + |x - a_j^- t|)^{-7/4} t^{-1/2} 
\int_{t/2}^t (t-s)^{-1/2} ds \\
&\le
C (1 + |x - a_j^- t|)^{-7/4},
\end{aligned}
\end{equation}
which, for $|x - a_j^- t| \ge \sqrt{t}$ decays faster than the claimed estimates.
For the case $|x - a_j^- t| \le \sqrt{t}$, we require only $t^{-3/4}$ decay, 
which is clear from (\ref{timedecayestimate2}).
For $|x| \le |a_j^-| t$, we divide the analysis into cases, $s \in [0, t/2]$
and $s \in [t/2, t]$.  For $s \in [0, t/2]$, we observe the inequality
\begin{equation} \label{gf1nl2balance4}
\begin{aligned}
&(1+|x - a_j^- (t-s)|)^{-1/2} (|x - a_j^- (t-s)| + s)^{-1/2} \\
&\times \quad  
(1 + |x - a_j^- (t-s)| + s)^{-3/4} 
(1 + |x - a_j^- (t-s)| + s^{1/2})^{-1/2} \\
&\le C\Big[
(1+|x - a_j^- t|)^{-1/2} (|x - a_j^- t| + s)^{-1/2} 
(1 + |x - a_j^- t| + s)^{-3/4} 
(1 + |x - a_j^- t| + s^{1/2})^{-1/2} \\
&+
(1+|x - a_j^- (t-s)|)^{-1/2} (|x - a_j^- t| + s)^{-1/2} 
(1 + |x - a_j^- t| + s)^{-3/4} 
(1 + |x - a_j^- t|^{1/2} + s^{1/2})^{-1/2} \Big].
\end{aligned}
\end{equation}
For the first estimate in (\ref{gf1nl2balance4}), we obtain an 
estimate, upon integration of the kernel, by 
\begin{equation*}
C_1 t^{-1/2} (1+|x - a_j^- t|)^{-1} \int_0^{t/2} 
(|x - a_j^- t| + s)^{-1/2} 
(1 + |x - a_j^- t| +s)^{-3/4} ds 
\le
C t^{-1/2} (1 + |x - a_j^- t|)^{-5/4},
\end{equation*}
while for the second we obtain an estimate by 
\begin{equation*}
C_1 t^{-1/2} (1 + |x - a_j^- t|)^{-3/2} 
\int_0^{t/2} (|x - a_j^- (t-s)|)^{-1/2} ds
\le
C (1 + |x - a_j^- t|)^{-3/2}.
\end{equation*}
For $s \in [t/2, t]$, we write 
\begin{equation} \label{gf1nl2balance5}
\begin{aligned}
(t-s)^{-1/2} &(1+|x - a_j^- (t-s)|)^{-1/2} \\
&\le C\Big[|x|^{-1/2} (1+|x - a_j^- (t-s)|)^{-1/2}
+
(t-s)^{-1/2} (1+|x - a_j^- (t-s)| + |x|)^{-1/2}
\Big].
\end{aligned}
\end{equation}
For the first estimate in (\ref{gf1nl2balance5}), we obtain an estimate
by 
\begin{equation*}
C_2 |x|^{-1/2} t^{-1/2} (1+t)^{-3/4} (1+t^{1/2})^{-1/2} 
\int_{t/2}^t (1 + |x - a_j^- (t-s)|)^{-1/2} ds
\le
C (1+t)^{-1} |x|^{-1/2},
\end{equation*}
which for $|x|$ bounded away from 0 is bounded by $\psi_2 (x, t)$,
while for the second estimate in (\ref{gf1nl1balance5}), we 
obtain an estimate by 
\begin{equation*}
C_2 (1+|x|)^{-1/2} t^{-1/2} (1+t)^{-3/4} (1 + t^{1/2})^{-1/2} 
\int_{t/2}^t (t-s)^{-1/2} ds
\le
C (1 + |x|)^{-1/2} (1+t)^{-1}.
\end{equation*}


{\it (\ref{Ginteract}(i)), term three.}
We next consider integration against the nonlinearity $(1+s)^{-1} e^{-\eta |y|}$
(which arises, for example, from the term 
$(\frac{\partial \bar{u}^\delta}{\partial \delta} \delta)^2$),
for which, in the case of our convection Green's kernel estimate, we have 
integrals of the form 
\begin{equation} \label{gf1nl3}
\int_0^t \int_{-\infty}^0 (t-s)^{-1} e^{-\frac{(x - y - a_j^- (t-s))^2}{M(t-s)}}
(1+s)^{-1} e^{-\eta |y|} dy ds.
\end{equation}
First, we observe a straightforward time decay estimate by 
\begin{equation} \label{gf1nl3timedecay}
\begin{aligned}
C_1 &t^{-1} \int_0^{t/2} (1+s)^{-1} ds
+
C_2 (1 + t)^{-1} \int_{t/2}^{t-1} (t-s)^{-1} ds
+
C_3 (1 + t)^{-1} \int_{t-1}^t (t-s)^{-1/2} ds \\
&\le
C (1+t)^{-1} \ln (e+t).
\end{aligned}
\end{equation}
In order to determine estimates in space as well, we observe the 
inequality
\begin{equation} \label{gf1nl3balance1}
\begin{aligned}
&e^{-\frac{(x - y - a_j^- (t-s))^2}{M(t-s)}} e^{-\eta |y|} \\
&\le
\Big[e^{-\frac{(x - a_j^- (t-s))^2}{\bar{M}(t-s)}} e^{-\eta |y|}
+
e^{-\frac{(x - y - a_j^- (t-s))^2}{M(t-s)}}
e^{-\eta_1 |x - a_j^- (t-s)|} e^{-\eta_2 |y|} \Big].
\end{aligned}
\end{equation}
For the first estimate in (\ref{gf1nl3balance1}), we have 
integrals 
\begin{equation} \label{gf1nl3balance1first}
\int_0^t \int_{-\infty}^0 
(t-s)^{-1} e^{-\frac{(x - a_j^- (t-s))^2}{\bar{M}(t-s)}} (1+s)^{-1} e^{-\eta |y|} dy ds.
\end{equation}
We have two cases to consider, $a_j^- < 0$ and $a_j^- > 0$, of which we focus
on the former.  For $|x| \ge |a_j^-| t$, there is no cancellation between
$x - a_j^- t$ and $a_j^- s$, and proceeding almost precisely as in 
(\ref{gf1nl3timedecay}), we obtain an estimate by 
\begin{equation*}
C (1+t)^{-1} \ln (e+t) e^{-\frac{(x - a_j^- t)^2}{Lt}}.
\end{equation*}
For $|x| \le |a_j^-| t$, we divide the analysis into cases, $s \in [0, t/2]$ 
and $s \in [t/2, t]$.  For $s \in [0, t/2]$, we observe the estimate
\begin{equation} \label{gf1nl3balance2}
\begin{aligned}
&e^{-\frac{(x - a_j^- (t-s))^2}{\bar{M}(t-s)}} (1+s)^{-1} \\
&\le
C \Big[ e^{-\frac{(x - a_j^- t)^2}{Lt}} (1+s)^{-1}
+
e^{-\frac{(x - a_j^- (t-s))^2}{\bar{M}(t-s)}} (1 + |x - a_j^- t|)^{-1}
\Big].
\end{aligned}
\end{equation}
For the first estimate in (\ref{gf1nl3balance2}), we proceed similarly 
as in (\ref{gf1nl3timedecay}), while for the second we have an estimate 
by 
\begin{equation*}
C_1 t^{-1} (1+|x-a_j^- t|)^{-1} 
\int_0^{t/2} e^{-\frac{(x - a_j^- (t-s))^2}{\bar{M}(t-s)}} ds
\le
C (1 + t)^{-1/2} (1 + |x - a_j^- t|)^{-1},
\end{equation*}
which is sufficient for $t \ge |x|/|a_j^-|$.  For $s \in [t/2, t]$,
we have
\begin{equation*}
(t-s)^{-1/2} e^{-\frac{(x - a_j^- (t-s))^2}{\bar{M}(t-s)}}
\le
C |x|^{-1/2} e^{-\frac{(x - a_j^- (t-s))^2}{2\bar{M}(t-s)}},
\end{equation*}
from which we obtain an estimate on (\ref{gf1nl3balance1first})
by 
\begin{equation*}
C_2 |x|^{-1/2} (1+t)^{-1} \int_{t/2}^t (t-s)^{-1/2}
e^{-\frac{(x - a_j^- (t-s))^2}{2\bar{M}(t-s)}} ds
\le
C |x|^{-1} (1+t)^{-1},
\end{equation*}
for which we have oberved the bound
\begin{equation*}
\int_{t/2}^t (t-s)^{-1/2}
e^{-\frac{(x - a_j^- (t-s))^2}{2\bar{M}(t-s)}} ds
\le C.
\end{equation*} 
The second estimate in (\ref{gf1nl3balance1}) can be analyzed similarly.



{\it (\ref{Ginteract}(i)), nonlinearity $v(y, s) \varphi (y, s)$.}
We next consider integration against the critical nonlinearity $v(y, s) \varphi (y, s)$,
which constituted the limiting estimate of the analysis in \cite{HR}.  Here, we 
refine the analysis of \cite{HR} both through refined estimates on 
$v(y, s)$ and through the application of the approach of Liu \cite{Liu85, Liu97}
described in Section \ref{liu}, which takes
advantage of improved decay for derivatives along the characteristic direction. 
It is precisely for this analysis that we must 
keep track of estimates on the characteristic derivatives 
$(\partial_t + a_j^- \partial_x) v$.  We proceed by dividing the integration over $s$
as,
\begin{equation} \label{threeintegrals1}
\begin{aligned}
\int_0^t &\int_{-\infty}^{+\infty} \tilde{G} (x, t-s; y) 
\Big(v(y, s) \varphi (y, s) \Big)_y dy ds \\
&=
\int_0^{\sqrt{t}} \int_{-\infty}^{+\infty} \tilde{G} (x, t-s; y) 
\Big(v(y, s) \varphi (y, s) \Big)_y dy ds \\
&+
\int_{\sqrt{t}}^{t - \sqrt{t}} \int_{-\infty}^{+\infty} \tilde{G} (x, t-s; y) 
\Big(v(y, s) \varphi (y, s) \Big)_y dy ds \\
&+
\int_{t - \sqrt{t}}^t \int_{-\infty}^{+\infty} \tilde{G} (x, t-s; y) 
\Big(v(y, s) \varphi (y, s) \Big)_y dy ds.
\end{aligned}
\end{equation} 
For the first integral on the right-hand side of  
(\ref{threeintegrals1}), we integrate by parts in $y$
and use the supremum 
estimate $\|v (y, s)\|_{L^\infty} \le C (1+s)^{-3/4}$ (valid for either
of our estimates on $v(y,s)$) to obtain integrals of 
the form, 
\begin{equation*}
\int_0^{\sqrt{t}} \int_{-\infty}^0 
(t-s)^{-1} e^{-\frac{(x - y - a_j^- (t-s))^2}{M(t-s)}}
(1 + s)^{-5/4} e^{-\frac{(y - a_k^- s)^2}{Ms}} dy ds,
\end{equation*}
for which we observe the equality
\begin{equation} \label{completedsquare}
e^{-\frac{(x - y - a_j^- (t-s))^2}{M(t-s)}} e^{-\frac{(y - a_k^- s)^2}{Ms}}
=
e^{-\frac{(x - a_j^- (t-s) - a_k^- s)^2}{Mt}}
e^{-\frac{t}{Ms(t-s)}(y - \frac{xs - (a_j^- + a_k^-)s(t-s)}{t})^2}
\end{equation}
(see Lemma 6 of \cite{HZ}).
Integrating over $y$, we immediately obtain an estimate by 
\begin{equation*}
C t^{-1/2} \int_0^{\sqrt{t}} (t-s)^{-1/2} (1+s)^{-3/4}
e^{-\frac{(x - a_j^- (t-s) - a_k^- s)^2}{Mt}} ds.
\end{equation*}
For $|x - a_j^- t| \le K \sqrt{t}$, any fixed $K$, we have
kernel decay $\exp(-(x-a_j^- t)^2/(Lt))$ by boundedness, while
for $|x - a_j^- t| \ge K \sqrt{t}$, with $K$ sufficiently large, 
we have 
\begin{equation} \label{kerneljuggle1}
|x - a_j^- (t-s) - a_k^- s| = |(x - a_j^- t) + (a_j^- - a_k^-) s|
\ge (1 - \frac{a_j^- - a_k^-}{K}) |x - a_j^- t|,
\end{equation}
and we again have kernel decay $\exp(-(x-a_j^- t)^2/(Lt))$.  In 
either case, we obtain a final estimate by 
\begin{equation*}
C (1 + t)^{-3/4} e^{-\frac{(x - a_j^- t)^2}{Lt}},
\end{equation*}
which is sufficient by (\ref{kerneltopsi1}).  Similarly, for the
third integral in (\ref{threeintegrals1}), we obtain, upon integration
in $y$ precisely as above, an estimate by 
\begin{equation*}
C t^{-1/2} \int_{t-\sqrt{t}}^t (t-s)^{-1/2} (1+s)^{-3/4}
e^{-\frac{(x - a_j^- (t-s) - a_k^- s)^2}{Mt}} ds.
\end{equation*}
Proceeding similarly as above, we obtain an estimate in this case
of the form 
\begin{equation*}
(1+t)^{-3/4} e^{-\frac{(x - a_k^- t)^2}{Lt}}.
\end{equation*}
For the second integral in (\ref{threeintegrals1}), we first consider the 
case $j = k$, for which, proceedingly similarly as above, 
we have integrals of the form 
\begin{equation*}
C t^{-1/2} \int_{\sqrt{t}}^{t - \sqrt{t}} (t-s)^{-1/2} (1+s)^{-3/4}
e^{-\frac{(x - a_j^- t)^2}{Mt}} ds,
\end{equation*}
which has been shown sufficient above.  In the case $j \ne k$, the 
claimed estimate does not follow from such a direct method, and 
we employ a Liu-type cancellation estimate, 
as described in Section \ref{liu}, based on integration by parts
in the characteristic direction. 
In order to clarify our analysis, we define the non-convecting variables
\begin{equation} \label{nonconvectingvariables}
\begin{aligned}
g (x, t; y) &= \tilde{G}^j (x, t, y-a_j^- t) \\
\phi (y, s) &= \varphi_k^- (y + a_k^- s, s) \\
V(y, s) &= v(y + a_k^- s, s),
\end{aligned}
\end{equation}
where 
\begin{equation} \label{convectionkernel}
\tilde{G}^j (x, t; y)
=
ct^{-1/2} e^{-\frac{(x - y - a_j^- (t-s))^2}{4\beta_j^- t}}
\quad c = r_j^- (l_j^-)^\text{tr} / \sqrt{4\pi\beta_j^-},
\end{equation} 
and $\varphi_k^-$ is as in (\ref{diffusionwaves}).  (We will consider 
corrections to $\tilde{G}^j$ at the end of the analysis.)  In this
notation, the second integral in (\ref{threeintegrals1}) becomes
\begin{equation} \label{thirdintegral1}
\int_{\sqrt{t}}^{t - \sqrt{t}} \int_{-\infty}^{+\infty} 
g(x, t-s; y + a_j^- (t-s)) \Big(V(y - a_k^- s, s) \phi(y - a_k^- s, s) \Big)_y
dy ds.
\end{equation}
Setting $\xi = y + a_j^- (t-s)$, (\ref{thirdintegral1}) becomes
\begin{equation} \label{thirdintegral1shifted}
\int_{\sqrt{t}}^{t - \sqrt{t}} \int_{-\infty}^{+\infty} 
g(x, t-s; \xi) \Big(V(\xi - a_j^- (t-s) - a_k^- s, s) 
\phi(\xi - a_j^- (t-s) - a_k^- s, s) \Big)_\xi d\xi ds.
\end{equation}
Denoting by $\partial_\tau$ differentiation with respect to 
the second dependent variable of each function (effectively, a 
characteristic derivative $\partial_t + a_k^- \partial_x$ on 
our original variables), we observe the differential relationship
\begin{equation} \label{partsins}
\begin{aligned}
&\Big(g(x,t-s;\xi) V (\xi - a_j^- (t-s) - a_k^- s, s) \phi (\xi - a_j^- (t-s) - a_k^- s, s) \Big)_s \\
&=
- g_\tau (x,t-s;\xi) V (\xi - a_j^- (t-s) - a_k^- s, s) \phi (\xi - a_j^- (t-s) - a_k^- s, s) \\
&+ (a_j^- - a_k^-) g(x,t-s;\xi) 
\Big(V (\xi - a_j^- (t-s) - a_k^- s, s) \phi (\xi - a_j^- (t-s) - a_k^- s, s)\Big)_\xi \\
&+
g(x,t-s;\xi) \Big(V (\xi - a_j^- (t-s) - a_k^- s, s) \phi (\xi - a_j^- (t-s) - a_k^- s, s) \Big)_\tau.
\end{aligned}
\end{equation}
Recalling that the case $j = k$ has already been considered, we can rearrange 
(\ref{partsins}) so that the integrand of (\ref{thirdintegral1shifted}) can be 
written as 
\begin{equation} \label{partsinsrearranged}
\begin{aligned}
g(x,t-s;\xi) 
&\Big(V (\xi - a_j^- (t-s) - a_k^- s, s) \phi (\xi - a_j^- (t-s) - a_k^- s, s)\Big)_\xi \\
&=
(a_j^- - a_k^-)^{-1}
\Big(g(x,t-s;\xi) V (\xi - a_j^- (t-s) - a_k^- s, s) \phi (\xi - a_j^- (t-s) - a_k^- s, s) \Big)_s \\
&+(a_j^- - a_k^-)^{-1}
g_\tau (x,t-s;\xi) V (\xi - a_j^- (t-s) - a_k^- s, s) \phi (\xi - a_j^- (t-s) - a_k^- s, s) \\
&-(a_j^- - a_k^-)^{-1}
g(x,t-s;\xi) \Big(V (\xi - a_j^- (t-s) - a_k^- s, s) \phi (\xi - a_j^- (t-s) - a_k^- s, s) \Big)_\tau.
\end{aligned}
\end{equation}
For the integration over the first expression on the right-hand side of 
(\ref{partsinsrearranged}), we change the order of integration, and evaluate integration over
$s$ to obtain 
\begin{equation} \label{partsinsrearrangedfirst}
\begin{aligned}
&\int_{\sqrt{t}}^{t-\sqrt{t}} \int_{-\infty}^{+\infty} 
\Big(g(x,t-s;\xi) V (\xi - a_j^- (t-s) - a_k^- s, s) \phi (\xi - a_j^- (t-s) - a_k^- s, s) \Big)_s 
d\xi ds \\
&=
\int_{-\infty}^{+\infty} 
g(x,\sqrt{t};\xi) V (\xi - a_j^- \sqrt{t} - a_k^- (t - \sqrt{t}), t-\sqrt{t}) 
\phi (\xi - a_j^- \sqrt{t} - a_k^- (t-\sqrt{t}), t-\sqrt{t}) d\xi \\
&-
\int_{-\infty}^{+\infty} 
g(x,t-\sqrt{t};\xi) V (\xi - a_j^- (t-\sqrt{t}) - a_k^- \sqrt{t}, \sqrt{t}) 
\phi (\xi - a_j^- (t-\sqrt{t}) - a_k^- \sqrt{t}, \sqrt{t}) d\xi.
\end{aligned}
\end{equation}
In each integral on the right hand side of (\ref{partsinsrearrangedfirst}),
we employ the estimate
\begin{equation} \label{vinfty}
\|V(\cdot, t)\|_{L^\infty} \le C (1 + t)^{-3/4},
\end{equation}
and proceed similarly as in (\ref{completedsquare}).
For the first, we obtain an estimate by 
\begin{equation*}
\begin{aligned}
\int_{-\infty}^{+\infty} &(\sqrt{t})^{-1/2} e^{-\frac{(x - \xi)^2}{M \sqrt{t}}}
(1 + (t-\sqrt{t}))^{-5/4} e^{-\frac{(\xi - a_j^- \sqrt{t} - a_k^- (t-\sqrt{t})^2}{M(t-\sqrt{t})}}
d\xi \\
&\le
C t^{-1/2} (1 + t)^{-3/4} 
e^{-\frac{(x - a_j^- \sqrt{t} - a_k^- (t-\sqrt{t})^2}{Mt}},
\end{aligned}
\end{equation*}
which gives an estimate by 
\begin{equation*}
C t^{-1/2} (1 + t)^{-3/4} 
e^{-\frac{(x - a_k^- t)^2}{Mt}}
\end{equation*}
(see (\ref{kerneljuggle1})).  For the second
integral on the right hand side of (\ref{partsinsrearrangedfirst}),
we obtain an estimate by 
\begin{equation*}
\begin{aligned}
\int_{-\infty}^{+\infty} &(t-\sqrt{t})^{-1/2} e^{-\frac{(x - \xi)^2}{M (t-\sqrt{t})}}
(1 + \sqrt{t})^{-5/4} e^{-\frac{(\xi - a_j^- (t-\sqrt{t}) - a_k^- \sqrt{t})^2}{M\sqrt{t}}} d\xi \\
&\le
C t^{-1/2} (1 + \sqrt{t})^{-3/4} 
e^{-\frac{(x - a_j^- (t-\sqrt{t}) - a_k^- \sqrt{t})^2}{M\sqrt{t}}},
\end{aligned}
\end{equation*}
which gives an estimate by 
\begin{equation*}
C (1 + t)^{-7/8} 
e^{-\frac{(x - a_j^- t)^2}{Mt}}.
\end{equation*}
For integration over the second integrand on the right hand side of 
\ref{partsinsrearranged}, 
we have integrals
\begin{equation} \label{partsinsrearrangedsecond}
\begin{aligned}
\Big|
\int_{\sqrt{t}}^{t-\sqrt{t}} &\int_{-\infty}^{+\infty} 
g_\tau (x,t-s;\xi) V (\xi - a_j^- (t-s) - a_k^- s, s) 
\phi (\xi - a_j^- (t-s) - a_k^- s, s) d\xi ds \Big| \\
&\le
C_1 \int_{\sqrt{t}}^{t-\sqrt{t}} \int_{-\infty}^{+\infty} 
(t-s)^{-3/2} e^{-\frac{(x-\xi)^2}{M(t-s)}} (1+s)^{-5/4} 
e^{-\frac{(\xi - a_j^- (t-s) -a_k^- s)^2}{Ms}} d\xi ds \\
&\le
C_2 t^{-1/2} \int_{\sqrt{t}}^{t-\sqrt{t}}
(t-s)^{-1} (1+s)^{-3/4} e^{-\frac{(x - a_j^- (t-s) - a_k^- s)^2}{Ms}} ds,
\end{aligned}
\end{equation}
where for the first inequality in (\ref{partsinsrearrangedsecond}), we 
have employed the estimate (\ref{vinfty}), while for the second we have 
used a calculation similar to (\ref{completedsquare}).  
In this last integral, we have three cases to consider, $a_k^- < 0 < a_j^-$,
$a_k < a_j^- < 0$, and $a_j^- < a_k^- < 0$, for which we focus on the 
second.  (We recall that the case $a_k^- = a_j^-$ has already been considered
above.)  For $|x| \ge |a_k^-| t$, there is no cancellation between 
summands in (\ref{tminussdecomp}), and we obtain kernel decay, 
\begin{equation*}
C t^{-1/2} (1+t)^{-3/4} \ln (e+t) e^{-\frac{(x - a_k^- t)^2}{Lt}},
\end{equation*}
while for $|x| \le |a_j^-| t$, there is no cancellation between 
summands in (\ref{sdecomp}), and we obtain kernel decay
\begin{equation*}
C t^{-1/2} (1+t)^{-3/4} \ln (e+t) e^{-\frac{(x - a_j^- t)^2}{Lt}}.
\end{equation*}
In either case, the seeming blow-up as $t \to 0$ can be eliminated by 
an alternative analysis in the case of $t$ bounded.
For $|a_j^-| t \le |x| \le |a_k^-| t$, we divide the analysis into cases, 
$s \in [\sqrt{t}, t/2]$ and $s \in [t/2, t-\sqrt{t}]$.  For $s \in [\sqrt{t}, t/2]$,
we observe through (\ref{sdecomp}) the inequality
\begin{equation} \label{gf1nl4balance1}
\begin{aligned}
&(1+s)^{-\gamma} e^{-\frac{(x - a_j^- (t-s) - a_k^- s)^2}{Ms}} \\
&\le
C\Big[ (1+s)^{-\gamma} e^{-\frac{(x - a_j^- t)^2}{Lt}}
+ (1+|x - a_j^- t|)^{-\gamma} e^{-\frac{(x - a_j^- (t-s) - a_k^- s)^2}{Ms}} \Big],
\end{aligned}
\end{equation}
where $\gamma = 3/4$.
For the first estimate in (\ref{gf1nl4balance1}), we proceed as in 
the case $|x| \le |a_j^-| t$, while for the second we obtain an estimate
by 
\begin{equation*}
\begin{aligned}
C_1 &t^{-3/2} (1+|x - a_j^- t|)^{-3/4} 
\int_{\sqrt{t}}^{t/2} e^{-\frac{(x - a_j^- (t-s) - a_k^- s)^2}{Ms}} ds \\
&\le
C t^{-1/2} (1+t)^{-1/2} (1+|x - a_j^- t|)^{-3/4},
\end{aligned}
\end{equation*}
which is sufficient.  For $s \in [t/2, t - \sqrt{t}]$, we observe through 
(\ref{tminussdecomp}) the inequality (\ref{gf1nl1balance3}) with $\gamma = 1$.
For the second estimate in (\ref{gf1nl1balance3}), we proceed as in the
case $|x| \ge |a_k^-|t$, while for the first we have an estimate by 
\begin{equation*}
\begin{aligned}
C_2 &t^{-1/2} (1+t)^{-3/4} |x - a_k^- t|^{-1} 
\int_{t/2}^{t - \sqrt{t}} e^{-\frac{(x - a_j^- (t-s) - a_k^- s)^2}{Ms}} ds \\
&\le
C (1+t)^{-3/4} |x - a_k^- t|^{-1},
\end{aligned}
\end{equation*}
which is sufficient for $|x - a_k^- t|$ bounded away from 0.  In the case of 
$|x - a_k^- t|$ bounded, we proceed alternatively.

The third term in (\ref{partsinsrearranged}) can be regarded as the crucial
piece, since it is here that we must keep track not only of estimates on 
$v(y,s)$, but also on estimates of characteristic derivatives on $v$.
We begin by expanding the characteristic derivative, 
\begin{equation} \label{partsinsrearrangedthird}
\begin{aligned}
\int_{\sqrt{t}}^{t-\sqrt{t}} &\int_{-\infty}^{+\infty} 
g(x, t-s; \xi) \Big(V(\xi - a_j^- (t-s) - a_k^- s, s)
\phi(\xi - a_j^- (t-s) - a_k^- s, s)\Big)_\tau d\xi ds \\
&=
\int_{\sqrt{t}}^{t-\sqrt{t}} \int_{-\infty}^{+\infty} 
g(x, t-s; \xi) V(\xi - a_j^- (t-s) - a_k^- s, s)
\phi_\tau (\xi - a_j^- (t-s) - a_k^- s, s) d\xi ds \\
&+
\int_{\sqrt{t}}^{t-\sqrt{t}} \int_{-\infty}^{+\infty} 
g(x, t-s; \xi) V_\tau (\xi - a_j^- (t-s) - a_k^- s, s)
\phi(\xi - a_j^- (t-s) - a_k^- s, s) d\xi ds.
\end{aligned}
\end{equation}
For the first integral on the right-hand side of 
(\ref{partsinsrearrangedthird}), we employ the estimate 
(\ref{vinfty}) to obtain an estimate by 
\begin{equation*}
\begin{aligned}
\int_{\sqrt{t}}^{t-\sqrt{t}} &\int_{-\infty}^{+\infty} 
(t-s)^{-1/2} e^{-\frac{(x-\xi)^2}{M(t-s)}} (1+s)^{-3/4} 
(1+s)^{-3/2} e^{-\frac{(\xi - a_j^- (t-s) - a_k^- s)^2}{Ms}}
d\xi ds \\
&\le
C t^{-1/2} \int_{\sqrt{t}}^{t-\sqrt{t}} (1 + s)^{-7/4}
e^{-\frac{(x - a_j^- (t-s) - a_k^- s)^2}{Mt}} ds,
\end{aligned}
\end{equation*}
which can be analyzed similarly as was (\ref{partsinsrearrangedsecond}).
For the second integral on the right-hand side of 
(\ref{partsinsrearrangedthird}), we note the following relation, 
also useful in calculations below, 
\begin{equation} \label{charestonv}
\begin{aligned}
&|V_\tau (\xi - a_j^- (t-s) - a_k^- s, s)| = |(\partial_t + a_k^- \partial_x) v| (\xi - a_j^- (t-s), s) \\
&\le
C s^{-1} s^{1/4} (1 + |\xi - a_j^- (t-s) - a_k^- s| + s^{1/2})^{-3/2} \\
&+
C s^{-1/2} \sum_{a_l^\pm \gtrless 0, l\ne k} (1 + |\xi - a_j^- (t-s) - a_l^- s| + s^{1/2})^{-3/2} \\
&+C s^{-1/2} (1 + |\xi - a_j^- (t-s)|)^{-1/2}
(1 + |\xi - a_j^- (t-s)| + s)^{-1/2} \\
&\times (1 + |\xi - a_j^- (t-s)| + s^{1/2})^{-1/2} 
I_{\{a_1^- s \le (\xi - a_j^- (t-s)) \le a_n^+ s \}} \\
&+
(1 + |\xi - a_j^- (t-s)|)^{-1} (1 + s)^{-1}
+
(1 + |\xi - a_j^- (t-s)| + s)^{-7/4}.
\end{aligned}
\end{equation}
For the first estimate in (\ref{charestonv}), we use a supremum norm 
to obtain integrals 
\begin{equation*}
\begin{aligned}
\int_{\sqrt{t}}^{t-\sqrt{t}} &\int_{-\infty}^{+\infty} 
(t-s)^{-1/2} e^{-\frac{(x-\xi)^2}{M(t-s)}} s^{-1} (1+s)^{-1} 
e^{-\frac{(\xi - a_j^- (t-s) - a_k^- s)^2}{Ms}}
d\xi ds \\
&\le
C t^{-1/2} \int_{\sqrt{t}}^{t-\sqrt{t}} s^{-1/2} (1 + s)^{-1}
e^{-\frac{(x - a_j^- (t-s) - a_k^- s)^2}{Ms}} ds,
\end{aligned}
\end{equation*}
which can be analyzed similarly as was (\ref{partsinsrearrangedsecond}).
For the remaining estimates in (\ref{charestonv}), the critical observation 
is that when integrated against the convecting diffusion kernel,
\begin{equation*}
e^{-\frac{(\xi - a_j^- (t-s) - a_k^- s)^2}{Ms}},
\end{equation*}
they give increased decay in $s$, from which the claimed estimates can 
readily be observed.
This completes the proof of (\ref{Ginteract})(i).
\medskip

{\bf Proof of (\ref{Ginteract}(ii)).}
Integration against the nonlinearity $\Phi(y,s)$ has been considered in 
\cite{HR}, and we need only refine one estimate from that calculation 
in order to conclude our claim.  The critical calculation regards 
integrals
\begin{equation*}
t^{-1/2} \int_{t/2}^{t-\sqrt{t}} (t-s)^{-1} (1+s)^{-1/2} 
e^{-\frac{(x - a_j^- (t-s) - a_k^- s)^2}{Mt}} ds,
\end{equation*} 
which was estimated in \cite{HR} by the expression,
\begin{equation*}
(1 + t)^{-3/4} (1 + |x - a_k^- t|)^{-1/2}.
\end{equation*}
Here, we observe that alternatively, we may proceed by observing
through (\ref{tminussdecomp}) the inequality (\ref{gf1nl1balance3})
with $\gamma = 1$.
For the first estimate in (\ref{gf1nl1balance3}),
we have an estimate by 
\begin{equation*}
\begin{aligned}
C_2 &t^{-1/2} (1+t)^{-1/2} |x - a_k^- t|^{-1}
\int_{t/2}^{t-\sqrt{t}} 
e^{-\frac{(x - a_j^- (t-s) - a_k^- s)^2}{Mt}} ds \\
&\le
C (1+t)^{-1/2} |x - a_k^- t|^{-1},
\end{aligned}
\end{equation*}
while for the second we have an estimate by 
\begin{equation*}
\begin{aligned}
C_2 &t^{-1/2} (1+t)^{-1} e^{-\frac{(x-a_k^- t)^2}{Lt}}
\int_{t/2}^{t-\sqrt{t}} (t-s)^{-1} ds \\
&\le
C (1 + t)^{-3/2} \ln (e+t) e^{-\frac{(x-a_k^- t)^2}{Lt}},
\end{aligned}
\end{equation*}
either of which is sufficient.
This completes the proof of (\ref{Ginteract}(ii)).
\medskip

{\bf Proof of (\ref{Ginteract}(iii)--(iv)), x-derivatives.}
As will be clear from our analysis of characteristic derivatives
just below, analysis
of $x$-derivatives is almost precisely the same as the case of 
characteristic derivatives for which the direction of differentiation
matches neither the convection rate of the kernel or the convection 
rate of the nonlinearity.  Since we consider the case of 
characteristic derivatives in great detail, we will omit the
case of $x$-derivatives. 


{\bf Proof of (\ref{Ginteract}(v)--(vi)), characteristic derivatives.}  
We next develop estimates on 
characteristic derivatives $(\partial_t + a_l^\pm \partial_x)$,
$a_l^\pm \gtrless 0$,
on the nonlinear interaction integrals.  These estimates are the primary
new contribution of the current analysis.

Observing that 
$\tilde{G} (x,0;y) = \delta_y (x)I$ by construction, we have, for $a_l^\pm \gtrless 0$, 
\begin{equation} \label{charderstart}
\begin{aligned}
(\partial_t + a_l^\pm \partial_x) 
\int_0^t &\int_{-\infty}^{+\infty} \tilde{G} (x,t-s;y)
\Big[\Phi (y,s) 
+ \mathcal{F} (\varphi, v, \frac{\partial \bar{u}^\delta}{\partial \delta}\delta)_y (y,s) \Big]
dy ds \\
&=
\Phi (x,t) 
+ \mathcal{F} (\varphi, v, \frac{\partial \bar{u}^\delta}{\partial \delta}\delta)_x (x,t) \\
&+
\int_0^t \int_{-\infty}^{+\infty} (\partial_t + a_l^\pm \partial_x) \tilde{G} (x,t-s;y)
\Big[\Phi (y,s) 
+ \mathcal{F} (\varphi, v, \frac{\partial \bar{u}^\delta}{\partial \delta}\delta)_y (y,s) \Big]
dy ds.
\end{aligned}
\end{equation}
For $\zeta(t) > 0$ so that 
\begin{equation}
\begin{aligned}
|v(x,t)| &\le \zeta (t) (\psi_1 (x, t) + \psi_2 (x, t)) \\
|\partial_x v (x,t)| 
&\le \zeta(t) \Big[ t^{-1/2}  (\psi_1 (x, t) + \psi_2 (x, t)) + \psi_3 (x,t) + \psi_4 (x,t) \Big],
\end{aligned}
\end{equation}
the first two terms on the right-hand side of (\ref{charderstart}) can 
be seen to be bounded by 
\begin{equation*}
\begin{aligned}
C &E_0 (1+t)^{-3/4} \sum_{a_j^\pm \gtrless 0} (1 + |x - a_j^\pm t| + t^{1/2})^{-3/2} \\
&+\zeta (t) \Big[t^{-1} (1+t)^{1/4} \psi_1^{j, \pm} (x, t) 
+ t^{-1/2} ( \bar{\psi}_1^{j, \pm} + \psi_2 (x,t) ) + \psi_3 (x, t) + \psi_4 (x, t) \Big].
\end{aligned}
\end{equation*}

{\bf Proof of (\ref{Ginteract}(vi)), nonlinearity $\Phi (y, s)$.} The basic elements of 
our proof are more easily seen in the case of estimate (\ref{Ginteract}(vi))---involving
the diffusion wave nonlinearity $\Phi (y, s)$---and so our approach will be to 
establish (\ref{Ginteract}(vi)) first and return to (\ref{Ginteract}(v)).   
Proceeding from the right-hand side of (\ref{charderstart}), with  
nonlinearity $\Phi (y, s)$, we have integrals of the form  
\begin{equation} \label{charderdiffusionwaves}
\int_0^t \int_{-\infty}^{+\infty}
(\partial_t + a_l^\pm \partial_x) \tilde{G} (x, t-s; y) \Phi (y, s) dy ds.
\end{equation}
In particuar, we focus on the nonlinearity $({\varphi_k^-}^2)_y$, with 
$\varphi_k^- (y, s)$ as defined in (\ref{diffusionwaves}), and the leading
order Green's kernel, 
\begin{equation*}
\tilde{G}^j (x, t; y) = ct^{-1/2} e^{-\frac{(x - y - a_j^- t)^2}{4 \beta_j t}},
\end{equation*}
with $c = r_j^- (l_j^-)^\text{tr}/\sqrt{4\pi\beta_j^-}$ (we will discuss corrections 
to $\tilde{G}^j (x, t; y)$ at the end of the analysis).
We divide the analysis into three parts, 
\begin{equation} \label{threeparts}
\begin{aligned}
\int_0^t &\int_{-\infty}^{+\infty}
(\partial_t + a_l^- \partial_x) \tilde{G}^j (x, t-s; y) ({\varphi_k^- (y, s)}^2)_y  dy ds \\
&=
\int_0^{\sqrt{t}} \int_{-\infty}^{+\infty}
(\partial_t + a_l^- \partial_x) \tilde{G}^j (x, t-s; y) ({\varphi_k^- (y, s)}^2)_y  dy ds \\
&+
\int_{\sqrt{t}}^{t-\sqrt{t}} \int_{-\infty}^{+\infty}
(\partial_t + a_l^- \partial_x) \tilde{G}^j (x, t-s; y) ({\varphi_k^- (y, s)}^2)_y  dy ds \\
&+
\int_{t-\sqrt{t}}^t \int_{-\infty}^{+\infty}
(\partial_t + a_l^- \partial_x) \tilde{G}^j (x, t-s; y) ({\varphi_k^- (y, s)}^2)_y  dy ds.
\end{aligned}
\end{equation}
According to the estimates of Proposition \ref{greenbounds}, the case $j=k$ does 
not occur (the diffusion waves have been chosen to eliminate precisely this case),
leaving three cases to consider, $l = j \ne k$, $l = k \ne j$, and 
$l \ne j, l \ne k, j\ne k$.  We remark that it is in this third case that the 
characteristic derivative acts similarly as a derivative with respect to $x$
only.

{\it Case 1: $l = j \ne k$.}
For the case $l = j \ne k$, the characteristic derivative of $\tilde{G}^j$ behaves like a time
derivative of the heat kernel, and we have 
\begin{equation*}
|(\partial_t + a_l^- \partial_x) \tilde{G}^j (x, t-s; y)|
\le C (t-s)^{-3/2} e^{-\frac{(x - y - a_j^- (t-s))^2}{M (t-s)}}.
\end{equation*}
For the first estimate in (\ref{threeparts}), upon integration by parts in $y$, we have
integrals of the form 
\begin{equation*}
\int_0^{\sqrt{t}} \int_{-\infty}^{+\infty} 
(t-s)^{-2} e^{-\frac{(x - y - a_j^- (t-s))^2}{M (t-s)}} (1+s)^{-1} 
e^{-\frac{(y - a_k^- s)^2}{Ms}} dy ds,
\end{equation*}
for which we observe the equality (\ref{completedsquare}).
Integration over $y$ leads immediately to an estimate by 
\begin{equation*}
C t^{-1/2} \int_0^{\sqrt{t}} (t-s)^{-3/2} (1+s)^{-1/2} 
e^{-\frac{(x - a_j^- (t-s) - a_k^- s)^2}{Mt}} ds.
\end{equation*}
In the event that $|x - a_j^- t| \le K \sqrt{t}$, we immediately have decay with 
scaling $\exp((x-a_j^-t)^2/(Lt))$, while for $|x - a_j^- t| \ge K \sqrt{t}$
we have
\begin{equation} \label{largechar}
|x - a_j^- (t-s) - a_k^- s| = |(x - a_j^- t) - (a_k^- - a_j^-) s|
\ge (1 - \frac{a_k^- - a_j^-}{K}) |x - a_j^- t|,
\end{equation} 
for which, with $K$ taken sufficiently large, we again have decay with 
scaling $\exp((x-a_j^-t)^2/(Lt))$.  We have, then, an estimate on this term 
of the form 
\begin{equation}
C_1 t^{-2} (1 + t^{1/2})^{-1/2} \int_0^{\sqrt{t}} 
e^{-\frac{(x - a_j^- (t-s) - a_k^- s)^2}{Mt}} ds
\le
C t^{-3/2} (1 + t^{1/2})^{-1/2} e^{-\frac{(x - a_j^- t)^2}{Lt}},  
\end{equation}
where in obtaining this last inequality we have reserved a small
part of the kernel $\exp(-(x - a_j^- (t-s) - a_k^- s)^2/(Mt))$
for integration.  Finally, the blow-up as $t \to 0$ can be reduced
by putting derivatives on the nonlinearity for small time.
For the third integral in (\ref{threeparts}), we integrate
the charactericstic derivative by parts to avoid blow-up
near $s = t$.  Observing the relation
\begin{equation*}
(\partial_t + a_j^- \partial_x) \tilde{G}^j (x, t-s; y)
= - (\partial_s + a_j^- \partial_y) \tilde{G}^j (x, t-s; y),
\end{equation*}
we have
\begin{equation} \label{charparts}
\begin{aligned}
\int_{t-\sqrt{t}}^t &\int_{-\infty}^{+\infty} 
(\partial_t + a_l^- \partial_x) \tilde{G}^j (x, t- s; y)
({\varphi_k^- (y, s)}^2)_y dy ds \\
&=
\int_{t-\sqrt{t}}^t \int_{-\infty}^{+\infty} 
(- \partial_s - a_l^- \partial_y) \tilde{G}^j (x, t- s; y)
({\varphi_k^- (y, s)}^2)_y dy ds \\
&=
\int_{-\infty}^{+\infty} \tilde{G}^j (x, \sqrt{t}; y)
({\varphi_k^- (y, t - \sqrt{t})}^2)_y dy \\
&-
\int_{-\infty}^{+\infty} \tilde{G}^j (x, 0; y)
({\varphi_k^- (y, t)}^2)_y dy \\
&+
\int_{t-\sqrt{t}}^t \int_{-\infty}^{+\infty} 
\tilde{G}^j (x, t- s; y)
(\partial_s + a_l^- \partial_y) ({\varphi_k^- (y, s)}^2)_y dy ds.  
\end{aligned}
\end{equation}
For the first integral in (\ref{charparts}), we have 
\begin{equation*}
\begin{aligned}
\int_{-\infty}^{+\infty} &(\sqrt{t})^{-1/2} e^{-\frac{(x - y - a_j^- \sqrt{t})}{M \sqrt{t}}}
(1+(t-\sqrt{t}))^{-3/2} e^{-\frac{(y - a_k^- (t-\sqrt{t}))^2}{M(t-\sqrt{t})}} dy \\
&\le
C t^{-1/2} (1+ (t - \sqrt{t}))^{-1} 
e^{-\frac{(x - a_j^- \sqrt{t} - a_k^- (t - \sqrt{t}))^2}{Mt}},
\end{aligned}
\end{equation*}
which gives an estimate by 
\begin{equation*}
t^{-1/2} (1 + t)^{-1} e^{-\frac{(x - a_k^- t)^2}{Lt}},
\end{equation*}
sufficient by an argument similar to (\ref{largechar}).
For the second integral in (\ref{charparts}), $\tilde{G}^j (x, 0; y)$
is a delta function, over which integration yields an estimate by 
\begin{equation*}
C (1+t)^{-3/2} e^{-\frac{(x - a_k^- t)^2}{Mt}}.
\end{equation*}
For the third integral in (\ref{charparts}), we have, upon noting
that the characteristic derivative is not along the direction of
propagation of $\varphi_k^-$, integrals of the form
\begin{equation*}
\begin{aligned}
\int_{t - \sqrt{t}}^t &\int_{-\infty}^{+\infty} (t - s)^{-1/2}
e^{-\frac{(x - y - a_j^- (t-s))^2}{M (t-s)}} 
(1+s)^{-3/2} e^{-\frac{(y - a_k^- s)^2}{Ms}} dy ds \\
&\le
C t^{-1/2} \int_{t - \sqrt{t}}^t (1+s)^{-3/2} 
e^{-\frac{(x - a_j^- (t-s) - a_k^- s)^2}{Mt}} ds.
\end{aligned}
\end{equation*}
Proceeding similarly as in (\ref{largechar}) and the surrounding 
estimates, we obtain an estimate by 
\begin{equation*}
C (1+t)^{-3/2} e^{-\frac{(x - a_j^- t)^2}{Lt}},
\end{equation*}
which has been shown sufficient in (\ref{kerneltopsi1}).
For the second integral in (\ref{threeparts}), we must proceed
by taking advantage of increased decay for derivatives along
characteristic directions.  Recalling our definitions 
(\ref{nonconvectingvariables}), we can write this integral 
in the form 
\begin{equation*}
\int_{\sqrt{t}}^{t-\sqrt{t}} \int_{-\infty}^{+\infty}
g_\tau (x, t-s; y+a_j^- (t-s)) (\phi (y - a_k^- s, s)^2)_y dy ds,
\end{equation*}
which upon the substitution $\xi = y + a_j^- (t-s)$ becomes
\begin{equation*}
\int_{\sqrt{t}}^{t-\sqrt{t}} \int_{-\infty}^{+\infty}
g_\tau (x, t-s; \xi) (\phi (\xi - a_j^- (t-s) - a_k^- s, s)^2)_\xi d\xi ds.
\end{equation*}
Proceeding similarly as in (\ref{partsins}), we write 
\begin{equation} \label{charderpartsins}
\begin{aligned}
&\Big(g_\tau (x,t-s;\xi) \phi (\xi - a_j^- (t-s) - a_k^- s, s)^2 \Big)_s \\
&=
- g_{\tau \tau} (x,t-s;\xi) \phi (\xi - a_j^- (t-s) - a_k^- s, s)^2 \\
&+ (a_j^- - a_k^-) g_\tau (x,t-s;\xi) 
(\phi (\xi - a_j^- (t-s) - a_k^- s, s)^2)_\xi \\
&+
g_\tau (x,t-s;\xi) \Big(\phi (\xi - a_j^- (t-s) - a_k^- s, s)^2 \Big)_\tau,
\end{aligned}
\end{equation}
which can be rearranged as
\begin{equation} \label{charderpartsinsrearranged}
\begin{aligned}
&g_\tau (x,t-s;\xi) (\phi (\xi - a_j^- (t-s) - a_k^- s, s)^2)_\xi \\
&=(a_j^- - a_k^-)^{-1}
\Big(g_\tau (x,t-s;\xi) \phi (\xi - a_j^- (t-s) - a_k^- s, s)^2 \Big)_s \\
&+(a_j^- - a_k^-)^{-1} g_{\tau \tau} (x,t-s;\xi) \phi (\xi - a_j^- (t-s) - a_k^- s, s)^2 \\
&-(a_j^- - a_k^-)^{-1}
g_\tau (x,t-s;\xi) \Big(\phi (\xi - a_j^- (t-s) - a_k^- s, s)^2 \Big)_\tau,
\end{aligned}
\end{equation}
where we recall again that the case $j = k$ does not occur here.
For integration over the first summand on the right-hand side of 
(\ref{charderpartsinsrearranged}), we exchange order of integration to 
obtain 
\begin{equation} \label{charderpartsinsrearrangedfirst}
\begin{aligned}
&\int_{\sqrt{t}}^{t-\sqrt{t}} \int_{-\infty}^{+\infty} 
\Big(g_\tau (x,t-s;\xi) \phi (\xi - a_j^- (t-s) - a_k^- s, s)^2 \Big)_s d\xi ds \\
&=
\int_{-\infty}^{+\infty} 
g_\tau (x,\sqrt{t};\xi) \phi (\xi - a_j^- \sqrt{t} - a_k^- (t-\sqrt{t}), t-\sqrt{t})^2 d\xi \\
&-
\int_{-\infty}^{+\infty} 
g_\tau (x,t-\sqrt{t};\xi) \phi (\xi - a_j^- (t-\sqrt{t}) - a_k^- \sqrt{t}, \sqrt{t})^2 d\xi.
\end{aligned}
\end{equation}
For the first integral in (\ref{charderpartsinsrearrangedfirst}), proceeding
similarly as in (\ref{completedsquare}), we estimate
\begin{equation*}
\begin{aligned}
\int_{-\infty}^{+\infty} &(\sqrt{t})^{-3/2} e^{-\frac{(x-\xi)^2}{M\sqrt{t}}}
(1 + (t-\sqrt{t}))^{-1} e^{-\frac{(\xi - a_j^- \sqrt{t} - a_k^- (t-\sqrt{t}))^2}{M(t-\sqrt{t})}}
d\xi \\
&\le
C t^{-1/2} (\sqrt{t})^{-1} (1 + (t-\sqrt{t}))^{-1/2} 
e^{-\frac{(x - a_j^- \sqrt{t} - a_k^- (t-\sqrt{t}))^2}{Mt}},
\end{aligned}
\end{equation*}
for which we conclude in a manner similar to (\ref{kerneljuggle1}) an estimate by 
\begin{equation*}
C t^{-1} (1+t)^{-1/2} e^{-\frac{(x - a_k^- t)^2}{Lt}}.
\end{equation*}
Proceeding similarly for the second integral in (\ref{charderpartsinsrearrangedfirst}),
we estimate
\begin{equation*}
\begin{aligned}
\int_{-\infty}^{+\infty} &(t - \sqrt{t})^{-3/2} e^{-\frac{(x-\xi)^2}{M(t-\sqrt{t})}}
(1 + \sqrt{t})^{-1} e^{-\frac{(\xi - a_j^- (t-\sqrt{t}) - a_k^- \sqrt{t})^2}{M\sqrt{t}}}
d\xi \\
&\le
C t^{-1/2} (t - \sqrt{t})^{-1} (1 + \sqrt{t})^{-1/2} 
e^{-\frac{(x - a_j^- (t-\sqrt{t}) - a_k^- \sqrt{t})^2}{Mt}},
\end{aligned}
\end{equation*}
for which we conclude an estimate by 
\begin{equation*}
C t^{-3/2} (1+t)^{-1/4} e^{-\frac{(x - a_j^- t)^2}{Lt}}.
\end{equation*}
For integration over the second summand on the right-hand side of 
(\ref{charderpartsinsrearranged}), we estimate
\begin{equation*}
\begin{aligned}
\int_{\sqrt{t}}^{t-\sqrt{t}} &\int_{-\infty}^{+\infty} 
(t-s)^{-5/2} e^{-\frac{(x-\xi)^2}{M(t-s)}}
(1 + s)^{-1} e^{-\frac{(\xi - a_j^- (t-s) - a_k^- s)^2}{Ms}} d\xi ds \\
&\le
Ct^{-1/2} \int_{\sqrt{t}}^{t-\sqrt{t}} (t-s)^{-2} (1+s)^{-1/2}
e^{-\frac{(x - a_j^- (t-s) - a_k^- s)^2}{Mt}} ds. 
\end{aligned}
\end{equation*}
For this last integral we have three cases to consider, 
$a_k^- < 0 < a_j^-$, $a_k^- < a_j^- < 0$, and
$a_j^- < a_k^- < 0$, of which we focus on the second.  In the 
case $|x| \ge |a_k^-| t$, there is no cancellation between 
summands in (\ref{tminussdecomp}), and we have an estimate by 
\begin{equation} \label{chardergf1nl1nocanc1}
\begin{aligned}
C_1 &t^{-5/2} (1+\sqrt{t})^{-1/2} e^{-\frac{(x - a_k^- t)^2}{Lt}}
\int_{\sqrt{t}}^{t/2} e^{-\epsilon \frac{(x - a_j^- (t-s) - a_k^- s)^2}{Mt}} ds \\
&+
C_2 t^{-1/2} (\sqrt{t})^{-2} (1+t)^{-1/2} e^{-\frac{(x - a_k^- t)^2}{Lt}}
\int_{t/2}^{t-\sqrt{t}} e^{-\epsilon \frac{(x - a_j^- (t-s) - a_k^- s)^2}{Mt}} ds \\
&\le
C t^{-3/2} e^{-\frac{(x - a_k^- t)^2}{Lt}},
\end{aligned}
\end{equation} 
while similarly in the case $|x| \le |a_j^-| t$, there is no cancellation
between summands in (\ref{sdecomp}), and we obtain an estimate by 
\begin{equation*}
Ct^{-3/2} e^{-\frac{(x - a_j^- t)^2}{Lt}}.
\end{equation*}
For the case $|a_j^-| t \le |x| \le |a_k^-| t$, we divide the analysis into 
cases $s \in [\sqrt{t}, t/2]$ and $s \in [t/2, t-\sqrt{t}]$ (we take $t$ 
large enough so that $\sqrt{t} < t/2$, proceeding alternatively for $t$
bounded, in which we need not establish $t$ decay).  For 
$s \in [\sqrt{t}, t/2]$, we observe through (\ref{sdecomp}) the 
inequality (\ref{gf1nl4balance1}) with $\gamma = 1/2$.
For the first estimate in (\ref{gf1nl4balance1}), we proceed as in 
(\ref{chardergf1nl1nocanc1}), while for the second we have an estimate 
by 
\begin{equation*}
C_1 t^{-5/2} (1 + |x-a_j^- t|)^{-1/2} \int_{\sqrt{t}}^{t/2}
e^{-\frac{(x - a_j^- (t-s) - a_k^- s)^2}{Mt}} ds
\le
C t^{-3/2} (1+t)^{-1/2} (1 + |x - a_j^- t|)^{-1/2}. 
\end{equation*}
For $s \in [t/2, t-\sqrt{t}]$, we observe through (\ref{tminussdecomp})
the inequality (\ref{gf1nl1balance3}) with $\gamma = 3/2$.
For the second estimate in (\ref{gf1nl1balance3}), we proceed as 
in (\ref{chardergf1nl1nocanc1}), while for the first we obtain 
an estimate by 
\begin{equation*}
\begin{aligned}
C_2 t^{-1/2} &(1+t)^{-1/2} |x - a_k^- t|^{-3/2} 
\int_{t/2}^{t-\sqrt{t}} (t-s)^{-1/2} e^{-\frac{(x - a_j^- (t-s) - a_k^- s)^2}{Mt}} ds \\
&\le C t^{-1/4} (1+t)^{-1/2} |x - a_k^- t|^{-3/2}, 
\end{aligned}
\end{equation*} 
which is better than the required estimate along the $a_k^-$ directions
(since $j \ne k$).
For integration over the third summand on the right-hand side of 
(\ref{charderpartsinsrearranged}), we estimate 
\begin{equation*}
\begin{aligned}
\int_{\sqrt{t}}^{t-\sqrt{t}} &\int_{-\infty}^{+\infty} 
(t-s)^{-3/2} e^{-\frac{(x-\xi)^2}{M(t-s)}}
(1 + s)^{-2} e^{-\frac{(\xi - a_j^- (t-s) - a_k^- s)^2}{Ms}} d\xi ds \\
&\le
Ct^{-1/2} \int_{\sqrt{t}}^{t-\sqrt{t}} (t-s)^{-1} (1+s)^{-3/2}
e^{-\frac{(x - a_j^- (t-s) - a_k^- s)^2}{Mt}} ds, 
\end{aligned}
\end{equation*}
for which we can proceed almost exactly as in our analysis of 
the second summand on the right-hand side of 
(\ref{charderpartsinsrearranged}) (the full analysis is omitted).

{\it Case 2: $l = k \ne j$.}
For the case $l = k \ne j$, the characteristic derivative of $\varphi_k^-$ behaves like a time
derivative of the heat kernel, and we have 
\begin{equation*}
(\partial_t + a_k^- \partial_x) \varphi_k^- (y,s)
=
{\mathbf O}((1+s)^{-3/2}) e^{-\frac{(y - a_k^- s)^2}{Ms}}.
\end{equation*}
In general, our strategy for this case will be to integrate by 
parts, shifting the characteristic derivative onto $\varphi$.
For the first integral in (\ref{threeparts}), after integration 
by parts in $y$, we estimate 
\begin{equation*}
\begin{aligned}
\int_0^{\sqrt{t}} &\int_{-\infty}^{+\infty} 
(t-s)^{-3/2} e^{-\frac{(x - y - a_j^- (t-s))^2}{M(t-s)}}
(1+s)^{-1} e^{-\frac{(y - a_k^- s)^2}{Ms}} dy ds \\
&\le
C t^{-1/2} \int_0^{\sqrt{t}} (t-s)^{-1} (1+s)^{-1/2}
e^{-\frac{(x - a_j^- (t-s) - a_k^- s)^2}{Mt}} ds,
\end{aligned}
\end{equation*}
which is bounded by 
\begin{equation*}
C t^{-1/2} (1+t)^{-3/4} e^{-\frac{(x - a_j^- t)^2}{Lt}},
\end{equation*}
sufficient for $j \ne l$ by the argument of 
(\ref{kerneltopsi1}).  For the third integral in (\ref{threeparts}),
we have 
\begin{equation} \label{threepartsthree}
\begin{aligned}
\int_{t-\sqrt{t}}^t &\int_{-\infty}^{+\infty} 
(\partial_t + a_k^- \partial_x) \tilde{G}^j (x,t;y)
\partial_y (\varphi_k^-(y,s)^2) dy ds \\
&=
-\int_{t-\sqrt{t}}^t \int_{-\infty}^{+\infty} 
(\partial_s + a_k^- \partial_y) \tilde{G}^j (x,t;y)
\partial_y (\varphi_k^-(y,s)^2) dy ds \\
&=
-\int_{-\infty}^{+\infty} \tilde{G}^j (x,0;y) \partial_y (\varphi_k^- (y,t))^2 dy \\
&+\int_{-\infty}^{+\infty} \tilde{G}^j (x,\sqrt{t};y) \partial_y (\varphi_k^- (y,t-\sqrt{t}))^2 dy \\
&-\int_{t-\sqrt{t}}^t \int_{-\infty}^{+\infty} 
\partial_y \tilde{G}^j (x,t-s;y)
(\partial_s + a_k^- \partial_y) ((\varphi_k^-)^2) dy ds.
\end{aligned}
\end{equation}
For the first integral in (\ref{threepartsthree}), we have an 
estimate by 
\begin{equation*}
|\varphi_k^- (x ,t) \partial_x \varphi_k^- (x, t)| 
\le
C (1+t)^{-3/2} e^{-\frac{(x-a_k^- t)^2}{Lt}},
\end{equation*}
while for the second we estimate 
\begin{equation} \label{refonce1}
\begin{aligned}
\int_{-\infty}^{+\infty} &(\sqrt{t})^{-1/2} e^{-\frac{(x - y - a_j^- \sqrt{t})^2}{M\sqrt{t}}}
(1 + (t-\sqrt{t}))^{-3/2} e^{-\frac{(y - a_k^- (t-\sqrt{t}))}{M(t-\sqrt{t})^2}} dy \\
&\le
C t^{-1/2} (1 + (t-\sqrt{t}))^{-1} 
e^{-\frac{(x - a_j^- \sqrt{t} -  a_k^- (t-\sqrt{t}))^2}{M(t-\sqrt{t})}},
\end{aligned}
\end{equation}
which gives an estimate by 
\begin{equation*}
t^{-1/2} (1+t)^{-1} e^{-\frac{(x-a_k^- t)^2}{Lt}}.
\end{equation*}
For the third integral on the right hand side of  (\ref{threepartsthree}), we estimate
\begin{equation*}
\begin{aligned}
\int_{t-\sqrt{t}}^t &\int_{-\infty}^{+\infty} 
(t-s)^{-1} e^{-\frac{(x - y - a_j^- (t-s))^2}{M(t-s)}}
(1+s)^{-2} e^{-\frac{(y - a_k^- s)^2}{Ms}} dy ds \\
&\le
C t^{-1/2} \int_{t-\sqrt{t}}^t (t-s)^{-1/2} (1+s)^{-3/2}
e^{-\frac{(x - a_j^- (t-s) - a_k^- s)^2}{Mt}} ds,
\end{aligned}
\end{equation*}
which gives an estimate by 
\begin{equation*}
t^{-1/4} (1+t)^{-3/2} e^{-\frac{(x-a_k^- t)^2}{Lt}}.
\end{equation*}
For the second integral in (\ref{threeparts}), we first integrate 
by parts, moving the characteristic derivative onto $\varphi^k$, 
and then proceed by a Liu-type characteristic derivative estimate
as described in Section \ref{liu}.
We have
\begin{equation} \label{maincharderest}
\begin{aligned}
\int_{\sqrt{t}}^{t-\sqrt{t}} &\int_{-\infty}^{+\infty} 
(\partial_t + a_k^- \partial_x) \tilde{G}^j (x,t;y)
\partial_y (\varphi^k(y,s)^2) dy ds \\
&=
-\int_{\sqrt{t}}^{t-\sqrt{t}} \int_{-\infty}^{+\infty} 
(\partial_s + a_k^- \partial_y) \tilde{G}^j (x,t;y)
\partial_y ((\varphi^k)^2) dy ds \\
&=
-\int_{-\infty}^{+\infty} \tilde{G}^j (x,\sqrt{t};y) \partial_y (\varphi^k(y,t-\sqrt{t}))^2 dy \\
&-\int_{-\infty}^{+\infty} \tilde{G}^j_y (x,t-\sqrt{t};y) \varphi^k(y,\sqrt{t})^2 dy \\
&+\int_{\sqrt{t}}^{t-\sqrt{t}} \int_{-\infty}^{+\infty} 
\tilde{G}^j (x,t -s ;y)
(\partial_s + a_k^- \partial_y) (\varphi^k(y,s)^2)_y dy ds.
\end{aligned}
\end{equation}
For the first integral on the right hand side of (\ref{maincharderest}),
we proceed exactly as in (\ref{refonce1}) to obtain an estimate by 
\begin{equation*}
t^{-1/2} (1+t)^{-1} e^{-\frac{(x-a_k^- t)^2}{Lt}},
\end{equation*}
while for the second we estimate 
\begin{equation*} 
\begin{aligned}
\int_{-\infty}^{+\infty} &(t-\sqrt{t})^{-1} e^{-\frac{(x - y - a_j^- (t-\sqrt{t}))^2}{M(t-\sqrt{t})}}
(1 + \sqrt{t})^{-1} e^{-\frac{(y - a_k^- \sqrt{t})^2}{M\sqrt{t}}} dy \\
&\le
C t^{-1/2} (t-\sqrt{t})^{-1/2} (1 + \sqrt{t})^{-1} (\sqrt{t})^{1/2} 
e^{-\frac{(x - a_j^- (t-\sqrt{t}) - a_k^- \sqrt{t})^2}{Mt}},
\end{aligned}
\end{equation*}
which gives an estimate by 
\begin{equation*}
t^{-3/4} (1+t)^{-1/2} e^{-\frac{(x-a_j^- t)^2}{Lt}}.
\end{equation*}
An argument similar to that of (\ref{kerneltopsi1}) shows that
this last estimate is sufficient for $j \ne l$, which is the current
setting.  For the third integral on the right hand side of 
(\ref{maincharderest}), we proceed in terms of the non-convecting
variables (\ref{nonconvectingvariables}), for which the integral
can be re-written as 
\begin{equation*}
\int_{\sqrt{t}}^{t-\sqrt{t}} \int_{-\infty}^{+\infty}
g (x, t-s; y + a_j^- (t-s)) (\partial_\tau (\phi (y-a_k^- s, s)^2))_y dy ds,
\end{equation*} 
where we recall that $\partial_\tau$ represents differentiation with respect 
to the second argument of $\phi$.  Setting $\xi = y + a_j^- (t-s)$, 
this last integral becomes 
\begin{equation*}
\int_{\sqrt{t}}^{t-\sqrt{t}} \int_{-\infty}^{+\infty} 
g (x, t-s; \xi) \Big(\partial_\tau (\phi (\xi - a_j^- (t-s) - a_k^- s, s)^2) \Big)_\xi
d\xi ds,
\end{equation*}
for which we write
\begin{equation} \label{partstrick2}
\begin{aligned}
&\Big(g (x, t-s; \xi) \partial_\tau (\phi (\xi - a_j^- (t-s) - a_k^- s, s)^2) \Big)_s \\
&=
-g_\tau (x, t-s; \xi) \partial_\tau (\phi (\xi - a_j^- (t-s) - a_k^- s, s)^2) \\
&+ (a_j^- - a_k^-)
g (x, t-s; \xi) (\partial_\tau (\phi (\xi - a_j^- (t-s) - a_k^- s, s)^2))_\xi \\
&+
g (x, t-s; \xi) \partial_{\tau \tau} (\phi (\xi - a_j^- (t-s) - a_k^- s, s)^2).
\end{aligned}
\end{equation}
For integration over the left hand side of (\ref{partstrick2}), we exchange the
order of integration to obtain 
\begin{equation} \label{partstrick2left}
\begin{aligned}
\int_{-\infty}^{+\infty} &g (x, \sqrt{t}; \xi) 
\partial_\tau (\phi (\xi - a_j^- \sqrt{t} - a_k^- (t-\sqrt{t}), t-\sqrt{t})^2) d\xi \\
&-
\int_{-\infty}^{+\infty} g (x, t-\sqrt{t}; \xi) 
\partial_\tau (\phi (\xi - a_j^- (t-\sqrt{t}) - a_k^- \sqrt{t}, \sqrt{t})^2) d\xi.
\end{aligned}
\end{equation}
For the first expression in (\ref{partstrick2left}), we estimate 
\begin{equation*}
\begin{aligned}
\int_{-\infty}^{+\infty} &(\sqrt{t})^{-1/2} e^{-\frac{(x-\xi)^2}{M\sqrt{t}}}
(1 + (t-\sqrt{t}))^{-2} e^{-\frac{(\xi - a_j^- \sqrt{t} - a_k^- (t-\sqrt{t}))^2}{M(t-\sqrt{t})}}
d \xi \\
&\le
C t^{-1/2} (1 + t)^{-3/2} 
e^{-\frac{(x - a_j^- \sqrt{t} - a_k^-(t- \sqrt{t}))^2}{Mt}},
\end{aligned}
\end{equation*}
which gives an estimate by 
\begin{equation*}
C t^{-1/2} (1 + t)^{-3/2} e^{-\frac{(x - a_k^- t)^2}{Lt}}.
\end{equation*}
For the second integral in (\ref{partstrick2left}), we estimate
\begin{equation*}
\begin{aligned}
\int_{-\infty}^{+\infty} &(t-\sqrt{t})^{-1/2} e^{-\frac{(x-\xi)^2}{M(t-\sqrt{t})}}
(1 + \sqrt{t})^{-2} e^{-\frac{(\xi - a_j^- (t-\sqrt{t}) - a_k^- \sqrt{t})^2}{M\sqrt{t}}}
d \xi \\
&\le
C t^{-1/2} (1 + \sqrt{t})^{-3/2} 
e^{-\frac{(x - a_j^- (t-\sqrt{t}) - a_k^- \sqrt{t})^2}{Mt}},
\end{aligned}
\end{equation*}
which gives an estimate by 
\begin{equation*}
C t^{-1/2} (1 + t)^{-3/4} e^{-\frac{(x - a_j^- t)^2}{Lt}},
\end{equation*}
sufficient for $l \ne j$.  For integration over the first term 
on the right-hand side of (\ref{partstrick2}), we estimate
\begin{equation} \label{partstrick2right1}
\begin{aligned}
\int_{\sqrt{t}}^{t-\sqrt{t}} &\int_{-\infty}^{+\infty} (t-s)^{-3/2}
e^{-\frac{(x-\xi)^2}{M(t-s)}} (1 + s)^{-2} 
e^{-\frac{(\xi - a_j^- (t-s) - a_k^- s)^2}{Ms}} d\xi ds \\
&\le
C t^{-1/2} \int_{\sqrt{t}}^{t-\sqrt{t}} (t-s)^{-1} (1 + s)^{-3/2} 
e^{-\frac{(x - a_j^- (t-s) - a_k^- s)^2}{Mt}} ds. 
\end{aligned}
\end{equation}  
We have three cases to consider here, $a_k^- < 0 < a_j^-$, 
$a_k^- < a_j^- < 0$, and $a_j^- < a_k^- < 0$, of which we
focus on the second.  In the case $|x| \ge |a_k^-| t$, we 
observe that there is no cancellation between summands in 
(\ref{tminussdecomp}), and consequently we obtain an 
estimate by 
\begin{equation} \label{partstrick2right1nocanc}
\begin{aligned}
C_1 &t^{-3/2} e^{-\frac{(x - a_k^- t)^2}{Lt}}
\int_{\sqrt{t}}^{t/2} (1+s)^{-3/2} ds \\
&+
C_2 t^{-1/2} (1+t)^{-3/2} e^{-\frac{(x - a_k^- t)^2}{Lt}}
\int_{t/2}^{t-\sqrt{t}} (t-s)^{-1} ds \\
&\le
Ct^{-1/2} (1+t)^{-1} e^{-\frac{(x - a_k^- t)^2}{Lt}}. 
\end{aligned}
\end{equation}
Similarly, for $|x| \le |a_j^-| t$, there is no cancellation between
summands in (\ref{sdecomp}), and the same calculation gives an 
estimate by 
\begin{equation*}
Ct^{-1/2} (1+t)^{-1} e^{-\frac{(x - a_j^- t)^2}{Lt}}.
\end{equation*}
For $|a_j^-| t \le |x| \le |a_k^-| t$, we divide the analysis into 
cases $s \in [\sqrt{t}, t/2]$ and $s \in [t/2, t-\sqrt{t}]$.  For
$s \in [\sqrt{t}, t/2]$, we observe through (\ref{sdecomp}) the inequality 
(\ref{gf1nl4balance1}) with $\gamma = 3/2$.
For the first estimate in (\ref{gf1nl4balance1}), we proceed 
as in (\ref{partstrick2right1nocanc}), while for the second we obtain 
an estimate by 
\begin{equation*}
C_1 t^{-3/2} (1+|x - a_j^- t|)^{-3/2} 
\int_{\sqrt{t}}^{t/2} e^{-\frac{(x - a_j^- (t-s) - a_k^- s)^2}{Mt}} ds
\le 
C t^{-1/2} (1+t)^{-1/2} (1+|x - a_j^- t|)^{-3/2},
\end{equation*}
which is sufficient.  For $s \in [t/2, t-\sqrt{t}]$, we observe through 
(\ref{tminussdecomp}) the estimate (\ref{gf1nl1balance3}) with $\gamma = 1$.
For the second estimate in (\ref{gf1nl1balance3}), we proceed
as in (\ref{partstrick2right1nocanc}), while for the first we have
an estimate by 
\begin{equation*}
C_2 t^{-1/2} |x - a_k^- t|^{-1} (1+t)^{-3/2}
\int_{t/2}^{t-\sqrt{t}} e^{-\frac{(x - a_j^- (t-s) - a_k^- s)^2}{Mt}} ds
\le
C (1 + t)^{-3/2} |x - a_k^- t|^{-1},  
\end{equation*}
which is sufficient since we are in the case $t \ge |x|/|a_k^-|$.  
For integration over the third term 
on the right-hand side of (\ref{partstrick2}), we estimate
\begin{equation} \label{partstrick2right3}
\begin{aligned}
\int_{\sqrt{t}}^{t-\sqrt{t}} &\int_{-\infty}^{+\infty} (t-s)^{-1/2}
e^{-\frac{(x-\xi)^2}{M(t-s)}} (1 + s)^{-3} 
e^{-\frac{(\xi - a_j^- (t-s) - a_k^- s)^2}{Ms}} d\xi ds \\
&\le
C t^{-1/2} \int_{\sqrt{t}}^{t-\sqrt{t}} (1 + s)^{-5/2} 
e^{-\frac{(x - a_j^- (t-s) - a_k^- s)^2}{Mt}} ds,
\end{aligned}
\end{equation}  
which can be analyzed similarly as in the immediately preceeding case.

{\it Case 3: $l \ne j$, $l \ne k$, $k\ne j$.}
For the case $l \ne j$, $l \ne k$, $k \ne j$, we proceed directly 
in the first and third integrals of (\ref{threeparts}) and through 
the method of Liu described in Section \ref{liu} and elsewhere 
for the second integral in (\ref{threeparts}).  
In each case, we obtain precisely the reduced decay estimate 
for characteristic directions other than $a_l^-$.  This concludes
the proof of Estimate (\ref{Ginteract}(vi)). 


{\bf Proof of (\ref{Ginteract}(v)), nonlinearity $\mathcal{F}$.} The analysis 
of characteristic derivatives in the case of nonlinearity $\mathcal{F}$
constitutes the single most involved section of the paper.  

{\it (\ref{Ginteract}(v)), Nonlinearity} $\varphi (y, s) v(y, s)$.
The critical nonlinearity of $\mathcal{F}$ 
is $\varphi (y, s) v(y, s)$,
for which we consider integrals
\begin{equation} \label{charderphiv}
\int_0^t \int_{-\infty}^{+\infty}
(\partial_t + a_l^\pm \partial_x) \tilde{G}^j (x, t-s; y) (\varphi_k^\pm (y, s) v(y,s))_y dy ds.
\end{equation}

{\it Case 1: $l = j$.} For the case $l=j$, and for $s \in [0, t/2]$, we integrate
by parts in $y$ and employ a supremum norm on $|v|$ to arrive at integrals
\begin{equation} \label{charderphiv1}
\begin{aligned}
\int_0^{t/2} &\int_{-\infty}^{+\infty} (t-s)^{-2} 
e^{-\frac{(x - y - a_j^- (t-s))^2}{M(t-s)}} (1+s)^{-5/4}
e^{-\frac{(y - a_k^- s)^2}{Ms}} dy ds \\
&\le
C t^{-1/2} \int_0^{t/2} (t-s)^{-3/2} (1+s)^{-3/4} 
e^{-\frac{(x - a_j^- (t-s) - a_k^- s)^2}{Mt}} ds.
\end{aligned}
\end{equation}
We have three cases to consider, $a_k^- < 0 < a_j^-$, $a_k^- \le a_j^- < 0$, 
and $a_j^- < a_k^- < 0$, of which we focus on the second.  (We observe that 
for the nonlinearity $\varphi(y,s) v(y,s)$, the case $a_j^- = a_k^-$
arises.)  For the case $a_k^- \le a_j^- < 0$, we first observe that in 
the case $|x| \ge |a_k^-| t$, there is no cancellation between summands in 
(\ref{tminussdecomp}), and proceeding similarly as in (\ref{partstrick2right1nocanc})
we determine an estimate by 
\begin{equation*}
C t^{-2} (1+t)^{1/4} e^{-\frac{(x - a_k^- t)^2}{Lt}}.
\end{equation*}
In the case $|x| \le |a_j^-| t$, there is no cancellation between summands in 
(\ref{sdecomp}) and we obtain an estimate by 
\begin{equation*}
C t^{-2} (1+t)^{1/4} e^{-\frac{(x - a_j^- t)^2}{Lt}}.
\end{equation*}
We note in particular that the case $a_j^- = a_k^-$ has been 
accomodated in this analysis.
For the critical case $|a_j^-| t \le |x| \le |a_k^-| t$ (now $j \ne k$),
we observe through (\ref{sdecomp}) the inequality 
(\ref{gf1nl4balance1}) with $\gamma = 3/4$.
For the first estimate in (\ref{gf1nl4balance1}), we proceed 
similarly as in (\ref{partstrick2right1nocanc}), while for the 
second, we have an estimate by 
\begin{equation*}
C_1 t^{-2} (1 + |x - a_j^- t|)^{-3/4} 
\int_0^{t/2} e^{-\frac{(x - a_j^- (t-s) - a_k^- s)^2}{Mt}} ds
\le 
C t^{-1} (1+t)^{-1/2} (1 + |x - a_j^- t|)^{-3/4},
\end{equation*}  
which is (precisely) sufficient for $t \ge |x|/|a_k^-|$. 
For $s \in [t/2, t-\sqrt{t}]$, we do not integrate by parts in $y$, 
and consequently obtain integrals
\begin{equation} \label{charderphiv2}
\begin{aligned}
\int_{t/2}^{t-\sqrt{t}} &\int_{-\infty}^{+\infty} (t-s)^{-3/2} 
e^{-\frac{(x - y - a_j^- (t-s))^2}{M(t-s)}} s^{-1/2} (1+s)^{-5/4}
e^{-\frac{(y - a_k^- s)^2}{Ms}} dy ds \\
&\le
C t^{-1/2} \int_{t/2}^{t-\sqrt{t}} (t-s)^{-1} (1+s)^{-5/4} 
e^{-\frac{(x - a_j^- (t-s) - a_k^- s)^2}{Mt}} ds,
\end{aligned}
\end{equation}
the last of which can be analyzed similarly as was 
(\ref{charderphiv1}).  We note that the form of the nonlinearity
in (\ref{charderphiv2}) arises from a supremum norm on $v$ (in 
the case that $\varphi$ is differentiated) and from the observation 
that the estimates on $v_y$ that do not decay at rate $s^{-5/4}$
have spatial decay different from that of $\varphi$ so that when 
multiplied by $\varphi$ the combination decays at rate 
$s^{-5/4}$ or better.
In the final case, $s \in [t-\sqrt{t}, t]$,
the expression $(t-s)^{-3/2}$ is not integrable up to $s = t$, and
we must proceed by integrating the characteristic derivative 
by parts.  We have
\begin{equation} \label{charparts2}
\begin{aligned}
\int_{t-\sqrt{t}}^t &\int_{-\infty}^{+\infty} 
(\partial_t + a_l^- \partial_x) \tilde{G}^j (x, t- s; y)
(\varphi_k^- (y, s) v(y,s))_y dy ds \\
&=
\int_{t-\sqrt{t}}^t \int_{-\infty}^{+\infty} 
(- \partial_s - a_l^- \partial_y) \tilde{G}^j (x, t- s; y)
(\varphi_k^- (y, s) v(y,s))_y dy ds \\
&=
\int_{-\infty}^{+\infty} \tilde{G}^j (x, 0; y)
(\varphi_k^- (y, t) v(y,t))_y dy \\
&-
\int_{-\infty}^{+\infty} \tilde{G}^j (x, \sqrt{t}; y)
(\varphi_k^- (y, t-\sqrt{t}) v(y,t-\sqrt{t}))_y dy \\
&-
\int_{t-\sqrt{t}}^t \int_{-\infty}^{+\infty} 
\tilde{G}^j_y (x, t- s; y)
(\partial_s + a_l^- \partial_y) (\varphi_k^- (y, s) v(y,s)) dy ds.  
\end{aligned}
\end{equation}
For the first integral on the right hand side of (\ref{charparts2}), we have the immediate estimates
\begin{equation*}
\begin{aligned}
|\varphi^k_x (x, t) v (x, t)| &\le CE_0 \zeta(t) 
(1 + t)^{-7/4} e^{-\frac{(x - a_k^- t)^2}{Lt}} \\
|\varphi^k (x, t) v_x (x, t)| &\le CE_0 \zeta(t) 
t^{-1/2} (1 + t)^{-5/4} e^{-\frac{(x - a_k^- t)^2}{Lt}}, 
\end{aligned}
\end{equation*}
where the second of these follows from the observation above that 
the combination $\varphi^k (x, t) v_x (x, t)$ decays at a rate faster
than a product of the supremum norms.
For the second integral on the right hand side of (\ref{charparts2}), we estimate 
\begin{equation*}
\begin{aligned}
\int_{-\infty}^{+\infty} &(\sqrt{t})^{-1/2} e^{-\frac{(x - y - a_j^- \sqrt{t})^2}{M\sqrt{t}}}
(t - \sqrt{t})^{-1/2} (1+(t-\sqrt{t}))^{-5/4} 
e^{-\frac{(y - a_k^- (t-\sqrt{t}))^2}{M(t-\sqrt{t})}} dy \\
&\le
C t^{-1/2} (1+t)^{-5/4} e^{-\frac{(x - a_j^- \sqrt{t} - a_k^- (t-\sqrt{t}))^2}{Mt}},
\end{aligned}
\end{equation*}
which gives an estimate by 
\begin{equation*}
C t^{-1/2} (1+t)^{-5/4} e^{-\frac{(x - a_k^- t)^2}{L t}}.
\end{equation*}
For the third integral on the right hand side of (\ref{charparts2}), we have
two nonlinearities to consider.  We begin with integrals  
\begin{equation} \label{charparts2three1}
\int_{t-\sqrt{t}}^t \int_{-\infty}^{+\infty} \tilde{G}^j_y (x, t-s; y) v(y, s)
(\partial_s + a_j^- \partial_y) \varphi^k (y, s) dy ds,
\end{equation} 
for which we estimate
\begin{equation*}
\begin{aligned}
\int_{t-\sqrt{t}}^t &\int_{-\infty}^{+\infty} (t-s)^{-1}
e^{-\frac{(x - y - a_j^- (t-s))^2}{M(t-s)}} (1+s)^{-7/4} e^{-\frac{(y - a_k^- s)^2}{Ms}}
dy ds \\
&\le
C t^{-1/2} \int_{t-\sqrt{t}}^t (t-s)^{-1/2} (1+s)^{-7/4} s^{1/2}
e^{-\frac{(x - a_j^- (t-s) - a_k^- s)^2}{Mt}} ds,
\end{aligned}
\end{equation*}
which gives an estimate by 
\begin{equation*}
C (1+t)^{-3/2} e^{-\frac{(x - a_k^- t)^2}{Lt}}.
\end{equation*}
For the second nonlinearity, we have integrals
\begin{equation} \label{charparts2three2}
\int_{t-\sqrt{t}}^t \int_{-\infty}^{+\infty} \tilde{G}^j_y (x, t-s; y) \varphi^k (y, s)
(\partial_s + a_j^- \partial_y) v(y, s) dy ds.
\end{equation} 
Observing again that the combination $\varphi^k (y, s)
(\partial_s + a_j^- \partial_y) v(y, s)$ decays at rate $s^{-1} (1 + s)^{-3/4}$,
we see that we can proceed as in the previous case.

{\it Case 2: $l = k \ne j$.}  For the case $l = k$, we divide the 
analysis as in (\ref{threeparts}) into integrals
\begin{equation} \label{threeparts2}
\begin{aligned}
\int_0^t &\int_{-\infty}^{+\infty}
(\partial_t + a_l^- \partial_x) \tilde{G}^j (x, t-s; y) (\varphi_k^- (y, s) v(y, s))_y  dy ds \\
&=
\int_0^{\sqrt{t}} \int_{-\infty}^{+\infty}
(\partial_t + a_l^- \partial_x) \tilde{G}^j (x, t-s; y) (\varphi_k^- (y, t) v(y, s))_y  dy ds \\
&+
\int_{\sqrt{t}}^{t-\sqrt{t}} \int_{-\infty}^{+\infty}
(\partial_t + a_l^- \partial_x) \tilde{G}^j (x, t-s; y) (\varphi_k^- (y, t) v(y,s))_y  dy ds \\
&+
\int_{t-\sqrt{t}}^t \int_{-\infty}^{+\infty}
(\partial_t + a_l^- \partial_x) \tilde{G}^j (x, t-s; y) (\varphi_k^- (y, s) v(y, s))_y  dy ds.
\end{aligned}
\end{equation}
For the first integral in (\ref{threeparts2}), upon integration by parts in $y$, 
we estimate
\begin{equation*}
\begin{aligned}
\int_0^{\sqrt{t}} &\int_{-\infty}^{+\infty}
(t-s)^{-3/2} e^{-\frac{(x - y - a_j^- (t-s))^2}{M(t-s)}}
(1+s)^{-5/4} e^{-\frac{(y - a_k^- s)^2}{Ms}} dy ds \\
&\le
C t^{-1/2} \int_0^{\sqrt{t}} (t-s)^{-1} (1+s)^{-5/4} s^{1/2}
e^{-\frac{(x - a_j^- (t-s) - a_k^- s)^2}{Mt}} ds,
\end{aligned}
\end{equation*}
which gives an estimate by 
\begin{equation*}
C t^{-1/2} (1+t)^{-3/4} e^{-\frac{(x - a_j^- t)^2}{Lt}},
\end{equation*}
sufficient for $l \ne j$.  
For the third integral in (\ref{threeparts2}), we do not integrate by 
parts in $y$, and consequently we can estimate
\begin{equation*}
\begin{aligned}
\int_{t-\sqrt{t}}^t &\int_{-\infty}^{+\infty}
(t-s)^{-1} e^{-\frac{(x - y - a_j^- (t-s))^2}{M(t-s)}}
s^{-1/2} (1+s)^{-5/4} e^{-\frac{(y - a_k^- s)^2}{Ms}} dy ds \\
&\le
C t^{-1/2} \int_{t-\sqrt{t}}^t (t-s)^{-1/2} (1+s)^{-5/4}
e^{-\frac{(x - a_j^- (t-s) - a_k^- s)^2}{Mt}} ds,
\end{aligned}
\end{equation*}
which gives an estimate by 
\begin{equation*}
C t^{-1/2} (1+t)^{-1} e^{-\frac{(x - a_k^- t)^2}{Lt}},
\end{equation*}
which is sufficient.  
For the second estimate in (\ref{threeparts2}), we integrate by 
parts as in (\ref{charparts2}), moving the characteristic derivative
onto the nonlinearity, and, when appropriate, moving the $y$ derivative 
onto $\tilde{G}^j$.
We have
\begin{equation} \label{charparts3}
\begin{aligned}
\int_{\sqrt{t}}^{t-\sqrt{t}} &\int_{-\infty}^{+\infty} 
(\partial_t + a_k^- \partial_x) \tilde{G}^j (x, t- s; y)
(\varphi_k^- (y, s) v(y,s))_y dy ds \\
&=
\int_{\sqrt{t}}^{t-\sqrt{t}} \int_{-\infty}^{+\infty} 
(- \partial_s - a_k^- \partial_y) \tilde{G}^j (x, t- s; y)
(\varphi_k^- (y, s) v(y,s))_y dy ds \\
&=
- \int_{-\infty}^{+\infty} \tilde{G}^j (x, \sqrt{t}; y)
(\varphi_k^- (y, t-\sqrt{t}) v(y,t-\sqrt{t}))_y dy \\
&-
\int_{-\infty}^{+\infty} \tilde{G}^j_y (x, t-\sqrt{t}; y)
\varphi_k^- (y, \sqrt{t}) v(y,\sqrt{t}) dy \\
&-
\int_{\sqrt{t}}^{t-\sqrt{t}} \int_{-\infty}^{+\infty} 
\tilde{G}^j_y (x, t- s; y)
(\partial_s + a_k^- \partial_y) (\varphi_k^- (y, s) v(y,s)) dy ds.  
\end{aligned}
\end{equation}
For the first integral on the right hand side of (\ref{charparts3}), we proceed 
precisely as with the second term in (\ref{charparts2}).  For the second
integral on the right hand side of (\ref{charparts3}), we estimate
\begin{equation*}
\begin{aligned}
\int_{-\infty}^{+\infty} &(t-\sqrt{t})^{-1} e^{-\frac{(x - y - a_j^- (t-\sqrt{t}))^2}{M(t-\sqrt{t})}}
(1 + \sqrt{t})^{-5/4} e^{-\frac{(y - a_k^2 \sqrt{t})^2}{M\sqrt{t}}} dy \\
&\le
C t^{-1} (1 + \sqrt{t})^{-3/4} e^{-\frac{(x - a_j^- (t-\sqrt{t}) - a_k^- \sqrt{t})^2}{Mt}}, 
\end{aligned}
\end{equation*}
which gives an estimate by 
\begin{equation*}
C t^{-1} (1 + t)^{-3/8} e^{-\frac{(x - a_j^- t)^2}{Lt}},
\end{equation*}
suffificient for $l \ne j$.
For the third integral on the right hand side of (\ref{charparts3}), we have
two integrals to consider, 
\begin{equation} \label{charparts3three}
\begin{aligned}
&\int_{\sqrt{t}}^{t-\sqrt{t}} \int_{-\infty}^{+\infty} 
\tilde{G}^j_y (x, t- s; y)
v(y, s) (\partial_s + a_k^- \partial_y) \varphi_k^- (y, s) dy ds \\
&+
\int_{\sqrt{t}}^{t-\sqrt{t}} \int_{-\infty}^{+\infty} 
\tilde{G}^j_y (x, t- s; y) \varphi_k^- (y, s)
(\partial_s + a_k^- \partial_y) v(y,s) dy ds.
\end{aligned}
\end{equation}
For the first integral in (\ref{charparts3three}), we employ the supremum norm 
on $|v(y,s)|$ and estimate
\begin{equation} \label{charparts3three1}
\begin{aligned}
\int_{\sqrt{t}}^{t-\sqrt{t}} &\int_{-\infty}^{+\infty} (t-s)^{-1} 
e^{-\frac{(x - y - a_j^- (t-s))^2}{M(t-s)}}
(1+s)^{-9/4} e^{-\frac{(y-a_k^- s)^2}{Ms}} dy ds \\
&\le
C t^{-1/2} \int_{\sqrt{t}}^{t-\sqrt{t}} (t-s)^{-1/2} 
(1+s)^{-7/4} e^{-\frac{(x - a_j^- (t-s) - a_k^- s)^2}{Mt}} ds.
\end{aligned}
\end{equation}
We have three cases to consider for this last integral, 
$a_k^- < 0 < a_j^-$, $a_k^- < a_j^- < 0$, and 
$a_j^- < a_k^- < 0$, of which we focus on the second
(we have already considered the case $j = k$).
For $|x| \ge |a_k^-| t$, there is no cancellation between 
summands in (\ref{tminussdecomp}), and we obtain an estimate by 
\begin{equation*}
C (1 + t)^{-3/2} e^{-\frac{(x - a_k^- t)^2}{Lt}}, 
\end{equation*}
while for $|x| \le |a_j^- t|$, there is no cancellation between 
summands in (\ref{sdecomp}) and we similarly obtain an estimate 
by 
\begin{equation*}
C (1 + t)^{-11/8} e^{-\frac{(x - a_j^- t)^2}{Lt}},
\end{equation*}
sufficient for $l \ne j$.
For the critical case $|a_j^-| t \le |x| \le |a_k^-| t$, we divide
the analysis into cases, $s \in [\sqrt{t}, t/2]$ and $s \in [t/2, t-\sqrt{t}]$.
For $s \in [\sqrt{t}, t/2]$ we observe through (\ref{sdecomp}) the  
inequality (\ref{gf1nl4balance1}) with $\gamma = 3/2$.
For the first estimate in (\ref{gf1nl4balance1}), we proceed as above to obtain 
an estimate by 
\begin{equation*}
C (1 + t)^{-11/8} e^{-\frac{(x - a_j^- t)^2}{Lt}},
\end{equation*}
while for the second we estimate
\begin{equation*}
\begin{aligned}
C &t^{-1} (1+|x - a_j^- t|)^{-3/2} \int_{\sqrt{t}}^{t/2}
(1+s)^{-1/4} e^{-\frac{(x - a_j^- (t-s) - a_k^- s)^2}{Mt}} ds \\
&\le
C (1 + t)^{-5/8} (1 + |x - a_j^- t|)^{-3/2},
\end{aligned}
\end{equation*}
sufficient for $l \ne j$.
For $s \in [t/2, t - \sqrt{t}]$, we observe through (\ref{tminussdecomp}) 
the inequality (\ref{gf1nl1balance3}) with $\gamma = 1/2$.
For the second estimate in (\ref{gf1nl1balance3}),
we immediately obtain an estimate by 
\begin{equation*}
C (1 + t)^{-3/2} e^{-\frac{(x - a_k^- t)^2}{Lt}},
\end{equation*} 
precisely as required, while for the first we estimate 
\begin{equation*}
\begin{aligned}
C &t^{-1/2} (1 + t)^{-7/4} |x - a_k^- t|^{-1/2}
\int_{t/2}^{t-\sqrt{t}} e^{-\frac{(x - a_j^- (t-s) - a_k^- s)^2}{Mt}} ds \\
&\le
C (1+t)^{-7/4} |x - a_k^- t|^{-1/2},
\end{aligned}
\end{equation*}
which is sufficient since $t \ge |x|/|a_k^-|$.   
For the second integral in (\ref{charparts3three}), we have five 
estimates on $|(\partial_s + a_k^- \partial_y) v|$ to consider
(see (\ref{charestonv})),
beginning with integrals
\begin{equation} \label{charparts3three2}
\int_{\sqrt{t}}^{t-\sqrt{t}} \int_{-\infty}^{+\infty} 
(t-s)^{-1} e^{-\frac{(x - y - a_j^- (t-s))^2}{M(t-s)}}
s^{-1} (1+s)^{-1/4} (1 + |y - a_k^- s| + s^{1/2})^{-3/2}
e^{-\frac{(y - a_k^- s)^2}{Ms}} dy ds. 
\end{equation}
Writing 
\begin{equation*}
x - y - a_j^- (t-s) = (x - a_j^- (t-s) - a_k^- s) - (y - a_k^- s),
\end{equation*}
we observe the estimate
\begin{equation} \label{charparts3three2balance1}
\begin{aligned}
&e^{-\frac{(x - y - a_j^- (t-s))^2}{M(t-s)}} (1 + |y - a_k^- s| + s^{1/2})^{-3/2}
e^{-\frac{(y - a_k^- s)^2}{Ms}} \\
&\le
C \Big[ e^{-\epsilon \frac{(x - y - a_j^- (t-s))^2}{M(t-s)}}
e^{-\frac{(x - a_j^- (t-s) - a_k^- s)^2}{\bar{M}(t-s)}}
(1 + |y - a_k^- s| + s^{1/2})^{-3/2}
e^{-\frac{(y - a_k^- s)^2}{Ms}} \\
&+
e^{-\frac{(x - y - a_j^- (t-s))^2}{M(t-s)}} (1 + |x - a_j^- (t-s) - a_k^- s| + s^{1/2})^{-3/2}
e^{-\frac{(x - a_j^- (t-s) - a_k^- s)^2}{\bar{M}s}} \Big].
\end{aligned}
\end{equation}
For the first estimate in (\ref{charparts3three2balance1}), we have
integrals
\begin{equation} \label{charparts3three2balance1first}
\int_{\sqrt{t}}^{t-\sqrt{t}} \int_{-\infty}^{+\infty}
(t-s)^{-1} e^{-\epsilon \frac{(x - y - a_j^- (t-s))^2}{M(t-s)}}
s^{-1} (1+s)^{-1/4}
e^{-\frac{(x - a_j^- (t-s) - a_k^- s)^2}{\bar{M}(t-s)}}
(1 + |y - a_k^- s| + s^{1/2})^{-3/2}
e^{-\frac{(y - a_k^- s)^2}{Ms}} dy ds.
\end{equation}
We have three cases to consider, $a_k^- < 0 < a_j^-$,
$a_k^- < a_j^- < 0$, and $a_j^- < a_k^- < 0$, of which we focus on 
the second.  For $|x| \ge |a_k^-| t$, we observe that there is 
no cancellation between summands in (\ref{tminussdecomp}), and 
consequently that we obtain an estimate by 
\begin{equation*}
\begin{aligned}
C_1& t^{-1} e^{-\frac{(x - a_k^- t)^2}{Lt}} \int_{\sqrt{t}}^{t/2}
s^{-1/2} (1+s)^{-1} e^{-\frac{((a_j^- - a_k^-) (t-s))^2}{\bar{M}(t-s)}} ds \\
&+
C_2 t^{-1} (1+t)^{-1} e^{-\frac{(x - a_k^- t)^2}{Lt}}
\int_{t/2}^{t-\sqrt{t}} (t-s)^{-1/2} e^{-\frac{(x - a_j^- (t-s) - a_k^- s)^2}{\bar{M}(t-s)}} \\
&\le
C (1 + t)^{-3/2} e^{-\frac{(x - a_k^- t)^2}{Lt}}.
\end{aligned}
\end{equation*}
For $|x| \le |a_j^-| t$, we observe that there is no cancellation between summands
in (\ref{tminussdecomp}), and consequently we obtain an estimate by 
\begin{equation*}
\begin{aligned}
C_1& t^{-1} e^{-\frac{(x - a_j^- t)^2}{Lt}} \int_{\sqrt{t}}^{t/2}
e^{-\epsilon \frac{(x - a_j^- (t-s) - a_k^- s)^2}{\bar{M}(t-s)}} 
s^{-1/2} (1+s)^{-1} ds \\
&+
C_2 t^{-1} (1 + t)^{-1} e^{-\frac{(x - a_j^- t)^2}{Lt}} \int_{t/2}^{t-\sqrt{t}}
(t-s)^{-1/2} e^{-\epsilon \frac{(x - a_j^- (t-s) - a_k^- s)^2}{\bar{M}(t-s)}} ds \\
&\le
C t^{-1/2} (1+t)^{-3/4} e^{-\frac{(x - a_j^- t)^2}{Lt}},   
\end{aligned}
\end{equation*}
which is sufficient for $j \ne l$.  For $|a_j^-| t \le |x| \le |a_k^-| t$,
we observe the inequality
\begin{equation} \label{charparts3three2balance2}
\begin{aligned}
&e^{-\frac{(x - a_j^- (t-s) - a_k^- s)^2}{\bar{M}(t-s)}} s^{-1} (1+s)^{-1/2} \\
&\le
C \Big[e^{-\frac{(x - a_j^- t)^2}{Lt}}  e^{-\epsilon \frac{(x - a_j^- (t-s) - a_k^- s)^2}{\bar{M}(t-s)}}
s^{-1} (1 + s)^{-1/2} \\
&+
e^{-\frac{(x - a_j^- (t-s) - a_k^- s)^2}{\bar{M}(t-s)}} |x - a_j^- t|^{-1} (1 + |x - a_j^- t|)^{-1/2}
\Big].
\end{aligned}
\end{equation}
For integration over the first estimate in (\ref{charparts3three2balance2}), 
we have an estimate by 
\begin{equation*}
\begin{aligned}
C_1 &t^{-1} e^{-\frac{(x - a_j^- t)^2}{Lt}} \int_{\sqrt{t}}^{t/2} 
s^{-1/2} (1+s)^{-1} e^{-\epsilon \frac{(x - a_j^- (t-s) - a_k^- s)^2}{\bar{M}(t-s)}} ds \\
&+
C_2 t^{-1} (1+t)^{-1} e^{-\frac{(x - a_j^- t)^2}{Lt}} \int_{t/2}^{t-\sqrt{t}}
(t-s)^{-1/2} e^{-\epsilon \frac{(x - a_j^- (t-s) - a_k^- s)^2}{\bar{M}(t-s)}} ds \\
&\le
C (1+t)^{-5/4} e^{-\frac{(x - a_j^- t)^2}{Lt}},
\end{aligned}
\end{equation*}
which is sufficient for $l \ne j$.  For integration over the second estimate in (\ref{charparts3three2balance2}), 
we have an estimate by 
\begin{equation*}
\begin{aligned}
C_1 &t^{-1} |x - a_j^- t|^{-1} (1 + |x - a_j^- t|)^{-1/2} 
\int_{\sqrt{t}}^{t/2} e^{-\frac{(x - a_j^- (t-s) - a_k^- s)^2}{\bar{M}(t-s)}} ds \\
&+
C_2 (1+t)^{-1/2} |x - a_j^- t|^{-1} (1 + |x - a_j^- t|)^{-1/2} \int_{t/2}^{t-\sqrt{t}}
(t-s)^{-1/2} e^{-\epsilon \frac{(x - a_j^- (t-s) - a_k^- s)^2}{\bar{M}(t-s)}} ds \\
&\le
C t^{-1/2} |x - a_j^- t|^{-1} (1 + |x - a_j^- t|)^{-1/2},
\end{aligned}
\end{equation*}
which, along with an alternative estimate in the case 
$|x - a_j^- t| \le \sqrt{t}$, is sufficient for $l \ne j$.  We 
remark that the critical observation in this calculation was
that since we require less decay along the $j$-characteristic, 
we use the same estimate (\ref{charparts3three2balance2}) for 
both cases $s \in [\sqrt{t}, t/2]$ and $s \in [t/2, t-\sqrt{t}]$.
For the second estimate in (\ref{charparts3three2balance1}), we have
integrals
\begin{equation} \label{charparts3three2balance1second}
\begin{aligned}
\int_{\sqrt{t}}^{t-\sqrt{t}} &\int_{-\infty}^{+\infty}
(t-s)^{-1} e^{-\frac{(x - y - a_j^- (t-s))^2}{M(t-s)}}
s^{-1} (1+s)^{-1/4}
(1 + |x - a_j^- (t-s) - a_k^- s| + s^{1/2})^{-3/2} \\
&\times 
e^{-\frac{(x - a_j^- (t-s) - a_k^- s)^2}{\bar{M}s}}
e^{- \epsilon \frac{(y - a_k^- s)^2}{Ms}} dy ds.
\end{aligned}
\end{equation}
We have three cases to consider, $a_k^- < 0 < a_j^-$,
$a_k^- < a_j^- < 0$, and $a_j^- < a_k^- < 0$, of which we focus on 
the second.  For $|x| \ge |a_k^-| t$, we observe that there is 
no cancellation between summands in (\ref{tminussdecomp}), 
and consequently we obtain an estimate by 
\begin{equation} \label{charparts3three2balance1secondnocanc}
\begin{aligned}
C_1 &t^{-1} (1 + |x - a_k^- t|)^{-3/2} \int_{\sqrt{t}}^{t/2} 
s^{-1/2} (1 + s)^{-1/4} e^{-\frac{(x - a_j^- (t-s) - a_k^- s)^2}{\bar{M}s}} ds \\
&+
C_2 t^{-1} (1+t)^{-1/4} (1 + |x - a_k^- t|)^{-3/2} \int_{t/2}^{t-\sqrt{t}}
(t-s)^{-1/2} ds \\
&\le
C (1 + t)^{-3/4} (1 + |x - a_k^- t|)^{-3/2},
\end{aligned}
\end{equation}
which is the required estimate since $l = k$.  For the case 
$|x| \le |a_j^-| t$, we have no cancellation between summands in 
(\ref{sdecomp}) and proceeding as in 
(\ref{charparts3three2balance1secondnocanc}) we immediately obtain an 
estimate by 
\begin{equation*}
C (1 + t)^{-3/4} (1 + |x - a_j^- t|)^{-3/2}.
\end{equation*}
For $|a_j^-| t \le |x| \le |a_k^-| t$, we divide the analysis into 
cases $s \in [\sqrt{t}, t/2]$ and $s \in [t/2, t-\sqrt{t}]$.  For 
$s \in [\sqrt{t}, t/2]$, we observe the estimate 
(\ref{charparts3three2balance2}).
For integration over the first estimate in (\ref{charparts3three2balance2}), 
we have an estimate by 
\begin{equation*}
\begin{aligned}
C_1 &t^{-1} e^{-\frac{(x - a_j^- t)^2}{Lt}} \int_{\sqrt{t}}^{t/2} 
s^{-1/2} (1+s)^{-1} e^{-\epsilon \frac{(x - a_j^- (t-s) - a_k^- s)^2}{\bar{M}s}}
ds \\
&\le
C (1+t)^{-5/4} e^{-\frac{(x - a_j^- t)^2}{Lt}},
\end{aligned}
\end{equation*}
while for integration over the second estimate in (\ref{charparts3three2balance2}), 
we have an estimate by
\begin{equation*}
\begin{aligned}
C_1 &t^{-1} |x - a_j^- t|^{-1} (1 + |x - a_j^- t|)^{-1/2} \int_{t/2}^{t-\sqrt{t}}
e^{-\frac{(x - a_j^- (t-s) - a_k^- s)^2}{\bar{M}s}} ds \\
&\le
C (1 + t)^{-1/2} |x - a_j^- t|^{-1} (1 + |x - a_j^- t|)^{-1/2},
\end{aligned}
\end{equation*}
which, along with an alternative estimate in the case 
$|x - a_j^- t| \le \sqrt{t}$, is sufficient for $l \ne j$.
For $s \in [t/2, t-\sqrt{t}]$, we compute an estimate directly from 
(\ref{charparts3three2balance1second})
\begin{equation*}
\begin{aligned}
C_2 &t^{-1} (1 + t)^{-1} \int_{t/2}^{t-\sqrt{t}} (t-s)^{-1/2}
e^{-\frac{(x - a_j^- (t-s) - a_k^- s)^2}{\bar{M}s}} ds \\
&\le
C (1 + t)^{-7/4},
\end{aligned}
\end{equation*}
which for $t \ge |x|/|a_k^-|$ gives an estimate by 
\begin{equation*}
C (1 + |x| + t)^{-7/4}.
\end{equation*}
For the remaining estimates in $|(\partial_s + a_l^- \partial_y) v|$, we 
have decay with with a different scaling than in the diffusion 
wave.  For example, we have terms of the form 
\begin{equation*}
\int_{\sqrt{t}}^{t-\sqrt{t}} \int_{-\infty}^{+\infty} 
(t-s)^{-1} e^{-\frac{(x - y - a_j^- (t-s))^2}{M(t-s)}}
s^{-1/2} (1+s)^{-1/2} (1 + |y - a_m^- s| + s^{1/2})^{-3/2}
e^{-\frac{(y - a_k^- s)^2}{Ms}} dy ds,
\end{equation*}
where $m \ne k$.  In such cases, we observe that for $y$ near $a_m^- s$
we have exponential decay in $s$ (which gives exponential decay in $\sqrt{t}$),
while for $y$ away from $a_m^- s$, we have integrals of the form 
\begin{equation*}
\int_{\sqrt{t}}^{t-\sqrt{t}} \int_{-\infty}^{+\infty} 
(t-s)^{-1} e^{-\frac{(x - y - a_j^- (t-s))^2}{M(t-s)}}
s^{-1/2} (1+s)^{-1/2} (1 + s)^{-3/2}
e^{-\frac{(y - a_k^- s)^2}{Ms}} dy ds,
\end{equation*}
which are better than previous cases (see (\ref{charparts3three1})).

{\it Case 3: $l \ne k$, $l \ne j$.}  For the case $l \ne k$, $l \ne j$, we divide the 
analysis into precisely the same three terms as in (\ref{threeparts2}).
For the first integral in (\ref{threeparts2}), upon integration by parts in 
$y$, we estimate
\begin{equation*}
\begin{aligned}
\int_0^{\sqrt{t}} &\int_{-\infty}^{+\infty} (t-s)^{-3/2} 
e^{-\frac{(x - y - a_j^- (t-s))^2}{M(t-s)}} (1+s)^{-5/4} e^{-\frac{(y - a_k^- s)^2}{Ms}}
dy ds \\
&\le
C_1 t^{-1/2} \int_0^{\sqrt{t}} (t-s)^{-1} (1+s)^{-3/4}
e^{-\frac{(x - a_j^- (t-s) - a_k^- s)^2}{Mt}} ds
\le
C t^{-1/2} (1+t)^{-7/8} e^{-\frac{(x - a_j^- t)^2}{Lt}},
\end{aligned}
\end{equation*}
which is sufficient for $l \ne j$.
For the third integral in (\ref{threeparts2}), we estimate
\begin{equation*}
\begin{aligned}
\int_{t - \sqrt{t}}^t &\int_{-\infty}^{+\infty} (t-s)^{-1}
e^{-\frac{(x - y - a_j^- (t-s))^2}{M(t-s)}} 
s^{-1/2} (1+s)^{-5/4} e^{-\frac{(y - a_k^- s)^2}{Ms}} \\
&\le
C t^{-1/2} \int_{t-\sqrt{t}}^t (t-s)^{-1/2} (1+s)^{-5/4}
e^{-\frac{(x - a_j^- (t-s) - a_k^- s)^2}{Mt}} ds
\le
C (1+t)^{-3/2} e^{-\frac{(x - a_k^- t)^2}{Lt}},
\end{aligned}
\end{equation*}
which is sufficient.  For the second integral in (\ref{threeparts2}),
we begin with the case $j = k$, for which we write 
\begin{equation} \label{jequalk}
\begin{aligned}
&\int_{\sqrt{t}}^{t-\sqrt{t}} \int_{-\infty}^{+\infty}
(\partial_t + a_l^- \partial_x) \tilde{G}^j (x, t-s; y) 
\Big(\varphi^k (y, s) v(y, s) \Big)_y dy ds \\
&=
- \int_{\sqrt{t}}^{t/2} \int_{-\infty}^{+\infty}
(\partial_t + a_l^- \partial_x) \tilde{G}^j_y (x, t-s; y) 
\varphi^k (y, s) v(y, s) dy ds \\
&+
\int_{t/2}^{t-\sqrt{t}} \int_{-\infty}^{+\infty}
(\partial_t + a_l^- \partial_x) \tilde{G}^j (x, t-s; y) 
\Big(\varphi^k (y, s) v(y, s) \Big)_y dy ds.
\end{aligned}
\end{equation}
For the first integral on the right hand side of (\ref{jequalk}), 
we estimate 
\begin{equation*}
\begin{aligned}
&\int_{\sqrt{t}}^{t/2} (t-s)^{-3/2} e^{-\frac{(x - y - a_j^- (t-s))^2}{M(t-s)}}
(1+s)^{-5/4} e^{-\frac{(y - a_j^- s)^2}{Ms}} ds \\
&\le
C t^{-1/2} \int_{\sqrt{t}}^{t/2} (t-s)^{-1} (1 + s)^{-5/4} s^{1/2} 
e^{-\frac{(x - a_j^- t)^2}{Mt}} ds \\
&\le
C_1 (1 + t)^{-5/4} e^{-\frac{(x - a_j^- t)^2}{Mt}} ds,
\end{aligned}
\end{equation*}
which is sufficient for $l \ne j$.  Similarly, for 
the second integral on the right hand side of (\ref{jequalk}), 
we estimate 
\begin{equation*}
\begin{aligned}
&\int_{t/2}^{t-\sqrt{t}} (t-s)^{-1} e^{-\frac{(x - y - a_j^- (t-s))^2}{M(t-s)}}
s^{-1/2} (1+s)^{-5/4} e^{-\frac{(y - a_j^- s)^2}{Ms}} ds \\
&\le
C t^{-1/2} \int_{t/2}^{t-\sqrt{t}} (t-s)^{-1/2} (1 + s)^{-5/4} 
e^{-\frac{(x - a_j^- t)^2}{Mt}} ds \\
&\le
C_1 (1 + t)^{-5/4} e^{-\frac{(x - a_j^- t)^2}{Mt}} ds,
\end{aligned}
\end{equation*}
which again is sufficient for $l \ne j$.
In the case $j \ne k$,
we employ the non-convecting variables (\ref{nonconvectingvariables}),
along with 
\begin{equation} \label{charnonconvectingvariable}
g_l (x, t; y + a_j^- (t-s)) := (\partial_t + a_l^- \partial_x)
\tilde{G}^j (x, t-s; y),
\end{equation}
for which we can write the second integral in (\ref{threeparts2})
as
\begin{equation*}
\int_{\sqrt{t}}^{t-\sqrt{t}} \int_{-\infty}^{+\infty} 
g_l (x, t-s; y + a_j^- (t-s)) 
\Big(\phi (y - a_k^- s, s) V (y - a_k^- s, s) \Big)_y dy ds.
\end{equation*}
Setting $\xi = y + a_j^- (t-s)$, this becomes 
\begin{equation} \label{threeparts2redux}
\int_{\sqrt{t}}^{t-\sqrt{t}} \int_{-\infty}^{+\infty} 
g_l (x, t-s; \xi) 
\Big(\phi (\xi - a_j^- (t-s) - a_k^- s, s) V (\xi - a_j^- (t-s) - a_k^- s, s) \Big)_y dy ds.
\end{equation}
Proceeding similarly as in (\ref{partsins}), we write 
\begin{equation} \label{charderpartsinsredux}
\begin{aligned}
&\Big(g_l (x,t-s;\xi) \phi (\xi - a_j^- (t-s) - a_k^- s, s) V (\xi - a_j^- (t-s) - a_k^- s, s) \Big)_s \\
&=
- g_{l \tau} (x,t-s;\xi) \phi (\xi - a_j^- (t-s) - a_k^- s, s) V (\xi - a_j^- (t-s) - a_k^- s, s) \\
&+ (a_j^- - a_k^-) g_l (x,t-s;\xi) 
\Big(\phi (\xi - a_j^- (t-s) - a_k^- s, s) V(\xi - a_j^- (t-s) - a_k^- s, s) \Big)_\xi \\
&+
g_l (x,t-s;\xi) \Big(\phi (\xi - a_j^- (t-s) - a_k^- s, s) V (\xi - a_j^- (t-s) - a_k^- s, s) \Big)_\tau.
\end{aligned}
\end{equation}
For integration over the left hand side of (\ref{charderpartsinsredux}),
we exchange the order of integration to obtain 
\begin{equation} \label{charderpartsinsreduxfirst}
\begin{aligned}
&\int_{-\infty}^{+\infty} 
g_l (x,\sqrt{t};\xi) \phi (\xi - a_j^- \sqrt{t} - a_k^- (t-\sqrt{t}), t-\sqrt{t}) 
V (\xi - a_j^- \sqrt{t} - a_k^- (t-\sqrt{t}), t-\sqrt{t}) d\xi \\
&-
\int_{-\infty}^{+\infty} 
g_l (x,t-\sqrt{t};\xi) \phi (\xi - a_j^- (t-\sqrt{t}) - a_k^- \sqrt{t}, \sqrt{t}) 
V (\xi - a_j^- (t-\sqrt{t}) - a_k^- \sqrt{t}, \sqrt{t}) d\xi.
\end{aligned}
\end{equation}
For the first integral in (\ref{charderpartsinsreduxfirst}), using the 
supremum norm of $V$, we estimate 
\begin{equation*}
\begin{aligned}
\int_{-\infty}^{+\infty} &(\sqrt{t})^{-1} e^{-\frac{(x - \xi)^2}{M\sqrt{t}}}
(1 + (t-\sqrt{t}))^{-5/4} e^{-\frac{(\xi - a_j^- \sqrt{t} - a_k^- (t - \sqrt{t}))^2}{M(t-\sqrt{t})}}
d\xi \\
&\le
C t^{-1/2} (\sqrt{t})^{-1/2} (1 + (t-\sqrt{t}))^{-5/4} (t - \sqrt{t})^{1/2} 
e^{-\frac{(x - a_j^- \sqrt{t} - a_k^- (t - \sqrt{t}))^2}{Mt}} \\
&\le
C (1 + t)^{-3/2} e^{-\frac{(x - a_k^- t)^2}{Lt}}.
\end{aligned}
\end{equation*}
For the second integral in (\ref{charderpartsinsreduxfirst}), again using the 
supremum norm of $V$, we estimate 
\begin{equation*}
\begin{aligned}
\int_{-\infty}^{+\infty} &(t-\sqrt{t})^{-1} e^{-\frac{(x - \xi)^2}{M(t-\sqrt{t})}}
(1 + \sqrt{t})^{-5/4} e^{-\frac{(\xi - a_j^- (t-\sqrt{t}) - a_k^- \sqrt{t})}{M(t-\sqrt{t})}}
d\xi \\
&\le
C t^{-1/2} (t-\sqrt{t})^{-1/2} (1 + \sqrt{t})^{-5/4} (\sqrt{t})^{1/2} 
e^{-\frac{(x - a_j^- (t-\sqrt{t}) - a_k^- \sqrt{t})}{Mt}} \\
&\le
C (1 + t)^{-11/8} e^{-\frac{(x - a_j^- t)}{Lt}},
\end{aligned}
\end{equation*}
which is slightly better than required for $l \ne j$ (we require $t^{-10/8}$).
For integration over the first term on the right hand side of 
(\ref{charderpartsinsredux}), we estimate
\begin{equation} \label{charderpartsinsreduxsecond}
\begin{aligned}
\int_{\sqrt{t}}^{t-\sqrt{t}} &\int_{-\infty}^{+\infty} 
(t-s)^{-2} e^{-\frac{(x - \xi)^2}{M(t-s)}} (1 + s)^{-5/4}
e^{-\frac{(\xi - a_j^- (t-s) - a_k^- s)^2}{Ms}} d\xi ds \\
&\le
C t^{-1/2} \int_{\sqrt{t}}^{t-\sqrt{t}} (t-s)^{-3/2} (1 + s)^{-5/4} s^{1/2}
e^{-\frac{(x - a_j^- (t-s) - a_k^- s)^2}{Mt}} ds.
\end{aligned}
\end{equation}
We have three cases to consider, $a_k^- < 0 < a_j^-$, 
$a_k^- < a_j^- < 0$, and $a_j^- < a_k^- < 0$, of which we 
focus on the second.  For $|x| \ge |a_k^-| t$, there is no 
cancellation between summands in (\ref{tminussdecomp}), and 
we can estimate 
\begin{equation} \label{charderpartsinsreduxsecondnocanc1}
\begin{aligned}
C_1 &t^{-2} e^{-\frac{(x - a_k^- t)^2}{Lt}} \int_{\sqrt{t}}^{t/2}
(1 + s)^{-5/4} s^{1/2} ds
+
C_2 (1 + t)^{-5/4} e^{-\frac{(x - a_k^- t)^2}{Lt}} \int_{t/2}^{t-\sqrt{t}}
(t-s)^{-3/2} ds \\
&\le
C (1 + t)^{-3/2} e^{-\frac{(x - a_k^- t)^2}{Lt}}. 
\end{aligned}
\end{equation}
In the case $|x| \le |a_j^-| t$, we have no cancellation between 
summands in (\ref{sdecomp}) and similarly obtain an estimate by 
\begin{equation*}
C (1 + t)^{-3/2} e^{-\frac{(x - a_k^- t)^2}{Lt}}. 
\end{equation*}
For the critical case $|a_j^-| t \le |x| \le |a_k^-| t$, we divide 
the analysis into subcases $s \in [\sqrt{t}, t/2]$ and $s \in [t/2, t-\sqrt{t}]$.
For $s \in [\sqrt{t}, t/2]$, we observe through (\ref{sdecomp}) 
the inequality (\ref{gf1nl4balance1}) with $\gamma = 3/4$.
For the first estimate in (\ref{gf1nl4balance1}), we proceed
similarly as in (\ref{charderpartsinsreduxsecondnocanc1}), while for the second
we estimate
\begin{equation*}
C_1 t^{-2} (1 + |x - a_j^- t|)^{-3/4} \int_{\sqrt{t}}^{t-\sqrt{t}}
(1+s)^{-1/2} s^{1/2} e^{-\frac{(x - a_j^- (t-s) - a_k^- s)^2}{Mt}} ds 
\le
C t^{-3/2} (1 + |x - a_j^- t|)^{-3/4},
\end{equation*}
which is sufficient for $t \ge |x|/|a_k^-|$.  
For $s \in [t/2, t-\sqrt{t}]$, we observe through (\ref{tminussdecomp})
the inequality (\ref{gf1nl1balance3}) with $\gamma = 3/2$.
For the second estimate in (\ref{gf1nl1balance3}), we proceed
similarly as in (\ref{charderpartsinsreduxsecondnocanc1}), while for the first
we estimate
\begin{equation*}
C_2 (1+t)^{-5/4} |x - a_k^- t|^{-3/2} \int_{t/2}^{t-\sqrt{t}} 
e^{-\frac{(x - a_j^- (t-s) - a_k^- s)^2}{Mt}} ds
\le
C (1 + t)^{-3/4} |x - a_k^- t|^{-3/2}.
\end{equation*}
For the third expression on the right hand side of 
(\ref{charderpartsinsredux}), we estimate 
\begin{equation} \label{charderpartsinsreduxthird}
\begin{aligned}
\int_{\sqrt{t}}^{t - \sqrt{t}} &\int_{-\infty}^{+\infty} 
(t-s)^{-1} e^{-\frac{(x - \xi)^2}{M(t-s)}} s^{-1} (1+s)^{-1} 
e^{-\frac{(\xi - a_j^- (t-s) - a_k^- s)^2}{Ms}} d\xi ds \\
&\le
C t^{-1/2} \int_{\sqrt{t}}^{t - \sqrt{t}} (t-s)^{-1/2} s^{-1/2} 
(1 + s)^{-1} e^{-\frac{(x - a_j^- (t-s) - a_k^- s)^2}{Ms}} ds,
\end{aligned}
\end{equation}
where we have once again observe the increased rate of time 
decay for the combination
\begin{equation*}
\Big(\phi (\xi - a_j^- (t-s) - a_k^- s, s) V (\xi - a_j^- (t-s) - a_k^- s, s) \Big)_\tau.
\end{equation*}
We have three cases to consider, $a_k^- < 0 < a_j^-$, 
$a_k^- < a_j^- < 0$, and $a_j^- < a_k^- < 0$, of which we 
focus on the second.  For $|x| \ge |a_k^-| t$, there is no 
cancellation between summands in (\ref{tminussdecomp}), and 
we can estimate 
\begin{equation} \label{charderpartsinsreduxthirdnocanc1}
\begin{aligned}
C_1 &t^{-1} e^{-\frac{(x - a_k^- t)^2}{Lt}} \int_{\sqrt{t}}^{t/2}
s^{-1/2} (1+s)^{-1} ds 
+ C_2 t^{-1} (1 + t)^{-1} e^{-\frac{(x - a_k^- t)^2}{Lt}} 
\int_{t/2}^{t-\sqrt{t}} (t-s)^{-1/2} ds \\
&\le
C (1+t)^{-5/4} e^{-\frac{(x - a_k^- t)^2}{Lt}},
\end{aligned}
\end{equation}
which is sufficient for $k \ne k$.  For $|x| \le |a_j^-| t$, 
there is no cancellation between summands in 
(\ref{sdecomp}) and we obtain an estimate by 
\begin{equation*}
C (1+t)^{-5/4} e^{-\frac{(x - a_j^- t)^2}{Lt}},
\end{equation*}
For $|a_j^-| t \le |x| \le |a_k^-| t$, we divide the analysis
into cases, $s \in [\sqrt{t}, t/2]$ and $s \in [t/2, \sqrt{t}]$.
For $s \in [\sqrt{t}, t/2]$, we observe the inequality
(\ref{gf1nl4balance1}) with $\gamma = 3/2$.  For the first
estimate in (\ref{gf1nl4balance1}), we proceed as in 
(\ref{charderpartsinsreduxthirdnocanc1}), while for the 
second we estimate
\begin{equation*}
\begin{aligned}
C_1 t^{-1} (1 + |x - a_j^- t|)^{-3/2} \int_{\sqrt{t}}^{t/2}
s^{-1/2} (1 + s)^{1/2} e^{-\frac{(x - a_j^- (t-s) - a_k^- s)^2}{Mt}} ds
\le
C (1 + t)^{-1/2} (1 + |x - a_j^- t|)^{-3/2},
\end{aligned}
\end{equation*}
which is precisely enough for $l \ne j$.  For $s \in [t/2, t-\sqrt{t}]$,
we observe through (\ref{tminussdecomp}) the inequality 
(\ref{gf1nl1balance3}) with $\gamma = 1/2$.  For the second estimate in 
(\ref{gf1nl1balance3}), we proceed similarly as in 
(\ref{charderpartsinsreduxthirdnocanc1}), while for the first 
we estimate
\begin{equation*}
\begin{aligned}
C_2 t^{-1} (1+t)^{-1} |x - a_k^- t|^{-1/2} 
\int_{t/2}^{t-\sqrt{t}} e^{-\frac{(x - a_j^- (t-s) - a_k^- s)^2}{Mt}} ds
\le
C (1 + t)^{-3/2} |x - a_k^- t|^{-1/2},
\end{aligned}
\end{equation*}
which is sufficient for $l \ne k$.  This ends the proof of 
(\ref{Ginteract}(vi)) for the leading order 
convection kernels $\tilde{G}^j (x, t; y)$.


{\it (\ref{Ginteract}(v)), Nonlinearity} $(v(y, s)^2)_y$.
We next consider integrals
\begin{equation} \label{chardervsquared}
\begin{aligned}
\int_0^t &\int_{-\infty}^{+\infty} 
(\partial_t + a_l^\pm \partial_x) \tilde{G}^j (x,t-s;y) (v(y,s)^2)_y dy ds \\
&= 
\int_0^{t/2} \int_{-\infty}^{+\infty} 
(\partial_t + a_l^\pm \partial_x) \tilde{G}^j (x,t-s;y) (v(y,s)^2)_y dy ds \\
&+
\int_{t/2}^t \int_{-\infty}^{+\infty} 
(\partial_t + a_l^\pm \partial_x) \tilde{G}^j (x,t-s;y) (v(y,s)^2)_y dy ds, 
\end{aligned}
\end{equation}
for which we have two cases to consider, $l \ne j$ and $l = j$.  

{\it Case 1: $l \ne j$.}  For the case $l \ne j$, we integrate
the first integral on the right hand side of (\ref{chardervsquared})
by parts in $y$ to obtain integrals
\begin{equation*}
\int_0^{t/2} \int_{-\infty}^{+\infty} 
\partial_y (\partial_t + a_l^\pm \partial_x) \tilde{G}^j (x,t-s;y) v(y,s)^2 dy ds.
\end{equation*}
Taking supremum norm on one $v(y, s)$, we have two cases two consider, one
for each estimate on $v$.  For the first, we have
\begin{equation} \label{chardervsquaredfirst}
\int_0^{t/2} \int_{-\infty}^{+\infty} 
(t-s)^{-3/2} e^{-\frac{(x - y - a_j^- (t-s))^2}{M(t-s)}} 
(1 + s)^{-3/4} (1 + |y - a_k^- s| + s^{1/2})^{-3/2} dy ds,
\end{equation}
where $a_k^- < 0$.  We observe through (\ref{maindecomp}) the 
inequality (\ref{gf1nl1balance1}) with $\gamma = 3/2$.  For 
the first estimate in (\ref{gf1nl1balance1}) we have integrals
\begin{equation}
\int_0^{t/2} \int_{-\infty}^{+\infty} 
(t-s)^{-3/2} e^{-\epsilon \frac{(x - y - a_j^- (t-s))^2}{M(t-s)}}
e^{-\frac{(x - a_j^- (t-s) - a_k^- s)^2}{\bar{M}(t-s)}} 
(1 + s)^{-3/4} (1 + |y - a_k^- s| + s^{1/2})^{-3/2} dy ds.
\end{equation}   
We have three cases to consider $a_k^- < 0 < a_j^- s$, 
$a_k^- \le a_j^- < 0$, and $a_j^- < a_k^- < 0$, of which 
we focus on the second.  In the case $s \in [0, \sqrt{t}]$,
we have 
\begin{equation*}
e^{-\frac{(x - a_j^- (t-s) - a_k^- s)^2}{\bar{M}(t-s)}}
\le
C e^{-\frac{(x - a_j^- t)^2}{Lt}},
\end{equation*}
from which we obtain an estimate by 
\begin{equation*}
C_1 t^{-3/2} e^{-\frac{(x - a_j^- t)^2}{Lt}} \int_0^{\sqrt{t}}
(1 + s)^{-1} ds 
\le
C (1 + t)^{-3/2} \ln (1 + t) e^{-\frac{(x - a_j^- t)^2}{Lt}},
\end{equation*}
sufficient for $l \ne j$.  For $s \in [\sqrt{t}, t/2]$, 
we begin with the case $|x| \ge |a_k^-| t$, for which 
there is no cancellation between summands in (\ref{tminussdecomp}),
and we obtain an estimate by 
\begin{equation} \label{chardervsquaredfirstnocanc1}
C_1 t^{-3/2} e^{-\frac{(x - a_k^- t)^2}{Lt}} \int_{\sqrt{t}}^{t/2}
(1 + s)^{-1} e^{-\frac{(x - a_j^- (t-s) - a_k^- s)^2}{\bar{M}(t-s)}} ds.
\end{equation}
In the event that $j = k$, we obtain an estimate by 
\begin{equation*}
C (1 + t)^{-3/2} \ln (1 + t) e^{-\frac{(x - a_k^- t)^2}{Lt}},
\end{equation*}
which is sufficient for $k = j \ne l$, whereas in the event that 
$j \ne k$, we obtain an estimate by 
\begin{equation*}
C (1 + t)^{-3/2} e^{-\frac{(x - a_k^- t)^2}{Lt}}.
\end{equation*} 
Proceeding similarly in the case $|x| \le |a_j^-| t$, we obtain 
an estimate by 
\begin{equation*}
C (1 + t)^{-3/2} \ln (1 + t) e^{-\frac{(x - a_j^- t)^2}{Lt}},
\end{equation*}
which is sufficient for $l \ne j$.  For $|a_j^-| t \le |x| \le |a_k^-| t$,
we observe through (\ref{sdecomp}) the inequality (\ref{gf1nl1balance2})
with $\gamma = 3/4$.
For the first estimate in (\ref{gf1nl1balance2}), we proceed as in 
(\ref{chardervsquaredfirstnocanc1}), while for the second 
we estimate
\begin{equation*}
C_1 t^{-3/2} (1 + |x - a_j^- t|)^{-1} \int_{\sqrt{t}}^{t/2} 
e^{-\frac{(x - a_j^- (t-s) - a_k^- s)^2}{\bar{M}(t-s)}} ds
\le
C t^{-1} (1 + |x - a_j^- t|)^{-1},
\end{equation*}
which is sufficient for $t \ge |x|/|a_k^-|$ and $l \ne j$.  For the 
second estimate in (\ref{gf1nl1balance1}) we have integrals
\begin{equation*}
\int_{\sqrt{t}}^{t/2} \int_{-\infty}^{+\infty} 
(t-s)^{-3/2} e^{-\frac{(x - y - a_j^- (t-s))^2}{M(t-s)}} 
(1 + s)^{-3/4} (1 + |y - a_k^- s| + |x - a_j^- (t-s) - a_k^- s| + s^{1/2})^{-3/2} dy ds,
\end{equation*}
for which we observe through (\ref{sdecomp}) the inequality 
(\ref{gf1nl1balance4}) with $\gamma = 3/4$.  Integrating the 
exponential kernel, for the first estimate in (\ref{gf1nl1balance4})
we obtain an estimate by 
\begin{equation*}
C_1 t^{-1} (1 + |x - a_j^- t|)^{-3/2} \int_{\sqrt{t}}^{t/2} 
(1 + s)^{-3/4} ds
\le
C t^{-3/4} (1 + |x - a_j^- t|)^{-3/2},
\end{equation*}
while for the second estimate in (\ref{gf1nl1balance4}), 
we obtain an estimate by 
\begin{equation*}
C_1 t^{-1} (1 + |x - a_j^- t|)^{-3/4} \int_{\sqrt{t}}^{t/2}
(1 + |x - a_j^- t| + s^{1/2})^{-3/2} ds
\le
C t^{-1} (1 + |x - a_j^- t|)^{-1}.
\end{equation*}
For the second estimate on $v(y, s)$, (\ref{chardervsquaredfirst})
is replaced by 
\begin{equation} \label{chardervsquaredsecond}
\int_0^{t/2} \int_{-|a_1^-| s}^{0} (t-s)^{-3/2}
e^{-\frac{(x - y - a_k^- (t-s))^2}{M(t-s)}} 
(1 + |y|)^{-1/2} (1 + |y| + s)^{-5/4} (1 + |y| + s^{1/2})^{-1/2}
dy ds, 
\end{equation}
for which we observe through (\ref{sdecomp}) the inequality
(\ref{gf1nl2balance1}).  For the first estimate in (\ref{gf1nl2balance1}),
we have integrals
\begin{equation*}
\int_0^{t/2} \int_{-|a_1^-| s}^{0} (t-s)^{-3/2}
e^{- \epsilon \frac{(x - y - a_j^- (t-s))^2}{M(t-s)}}
e^{-\frac{(x - a_j^- (t-s))^2}{\bar{M} (t-s)}} 
(1 + |y|)^{-1/2} (1 + |y| + s)^{-5/4} (1 + |y| + s^{1/2})^{-1/2}
dy ds,
\end{equation*}
for which we have two cases to consider, $a_j^- < 0$ and $a_j^- > 0$.
For the case $a_j^- > 0$, there is no cancellation between 
$x$ and $a_j^- (t-s)$ and the claimed estimate can be deduced in 
straightforward fashion.  For the case $a_j^- < 0$, we first
consider the subcase $|x| \ge |a_j^-| t$, for which there is no
cancellation between summands on the right hand side of 
\begin{equation*}
x - a_j^- (t-s) = (x - a_j^- t) + a_j^- s,
\end{equation*}
and we immediately obtain an estimate by 
\begin{equation} \label{chardervsquaredsecondnocanc}
C_1 t^{-3/2} e^{-\frac{(x - a_j^- t)^2}{Lt}}
\int_0^{t/2} (1 + s)^{-1} ds
\le
C (1 + t)^{-3/2} \ln (e+t) e^{-\frac{(x - a_j^- t)^2}{Lt}},
\end{equation}
which is sufficient for $l \ne j$.  For $|x| \le |a_j^-| t$,
we observe the inequality (\ref{gf1nl2balance2}).    For the 
first estimate in (\ref{gf1nl2balance2}), we proceed similarly 
as in (\ref{chardervsquaredsecondnocanc}), while for the second
we estimate
\begin{equation*}
C_1 t^{-3/2} (1 + |x - a_j^- t|)^{-3/2}
\int_0^{t/2} (1 + s)^{1/2} e^{-\frac{(x - a_j^- (t-s))^2}{\bar{M} (t-s)}} 
\le
C (1 + t)^{-1/2} (1 + |x - a_j^- t|)^{-3/2},
\end{equation*}
which is sufficient for $l \ne j$.  For the 
second estimate in (\ref{gf1nl2balance1}), we have integrals
\begin{equation*}
\begin{aligned}
\int_0^{t/2} &\int_{-|a_1^-| s}^{0} (t-s)^{-3/2}
e^{-\frac{(x - y - a_j^- (t-s))^2}{M(t-s)}}
(1 + |y| + |x - a_j^- (t-s)|)^{-1/2} 
(1 + |y| + |x - a_j^- (t-s)| + s)^{-5/4} \\
&\times (1 + |y| + |x - a_j^- (t-s)| + s^{1/2})^{-1/2}
dy ds,
\end{aligned}
\end{equation*}
for which we have two cases to consider, $a_j^- < 0$ and $a_j^- > 0$,
and as before we need focus only on the former.  
For $a_j^- > 0$ and $|x| \ge |a_j^-| t$, there is no cancellation between $x - a_j^- t$
and $a_j^- s$, and we immediately obtain an estimate by 
\begin{equation} \label{chardervsquaredsecondnocanc2}
C_1 t^{-1} (1 + |x - a_j^- t|)^{-3/2} 
\int_0^{t/2} (1 + |x - a_j^- (t-s)|)^{-3/4} ds
\le
C (1+t)^{-3/4} (1 + |x - a_j^- t|)^{-3/2}.  
\end{equation}
For $|x| \le |a_j^-| t$,
we observe the inequality (\ref{gf1nl2balance4}).  For the first
estimate in (\ref{gf1nl2balance4}), we proceed similarly as in 
(\ref{chardervsquaredsecondnocanc2}), while for the second we 
obtain an estimate by 
\begin{equation*}
C_1 t^{-1} (1 + |x - a_j^- t|)^{-3/2} \int_0^{t/2} 
(1 + |x - a_j^- (t-s)|)^{-1/2} ds
\le
C (1+t)^{-1/2} (1 + |x - a_j^- t|)^{-3/2},
\end{equation*} 
which is sufficient for $l \ne j$.  For the case $s \in [t/2, t]$,
we need to shift the characteristic derivative onto the 
nonlinearity.  We accomplish this precisely as in (\ref{charparts}), 
computing
\begin{equation} \label{charpartsvsquared}
\begin{aligned}
\int_{t/2}^t &\int_{-\infty}^{+\infty} 
(\partial_t + a_l^- \partial_x) \tilde{G}^j (x, t- s; y)
(v (y, s)^2)_y dy ds \\
&=
\int_{t/2}^t \int_{-\infty}^{+\infty} 
(- \partial_s - a_l^- \partial_y) \tilde{G}^j (x, t- s; y)
(v (y, s)^2)_y dy ds \\
&=
- \int_{-\infty}^{+\infty} \tilde{G}^j_y (x, t/2; y)
v (y, t/2)^2 dy \\
&-
\int_{-\infty}^{+\infty} \tilde{G}^j (x, 0; y)
(v (y, t)^2)_y dy \\
&-
\int_{t/2}^t \int_{-\infty}^{+\infty} 
\tilde{G}^j_y (x, t- s; y)
(\partial_s + a_l^- \partial_y) v (y, s)^2 dy ds,  
\end{aligned}
\end{equation}
where in the first and last integrals on the right hand side 
we have additionally integrated by parts 
in $y$.  For the first integral on the right hand side of 
(\ref{charpartsvsquared}), and for the first estimate on $v (y, s)$,
we have integrals
\begin{equation} \label{charpartsvsquaredfirst}
\begin{aligned}
\int_{-\infty}^{+\infty} &t^{-1} e^{-\frac{(x - y - a_k^- (t/2))^2}{M(t/2)}}
(1 + |y - a_k^- (t/2)| + (t/2)^{1/2})^{-3/2} (1 + (t/2))^{-3/4} dy \\
&\le
C t^{-1/2} (1 + t)^{-3/2},
\end{aligned}
\end{equation}
which is sufficient for $t \ge c|x|$, for any constant $c > 0$.  In the
case $|x| \ge t/c$, we observe the decomposition 
\begin{equation*}
x - y - \frac{1}{2} a_j^- t
=
(x - \frac{1}{2} a_j^- t - \frac{1}{2} a_k^- t) - (y - \frac{1}{2} a_k^- t),
\end{equation*}
through which we observe the inequality 
\begin{equation} \label{charpartsvsquaredfirstbalance1}
\begin{aligned}
&e^{-\frac{(x - y - a_j^- (t/2))^2}{M(t/2)}} (1 + |y - \frac{1}{2}a_k^- t| + (t/2)^{1/2})^{-3/2} \\
&\le
C\Big[ e^{-\epsilon \frac{(x - y - a_j^- (t/2))^2}{M(t/2)}}
e^{-\frac{(x - a_j^- (t/2) - a_k^- (t/2))^2}{Lt}}
(1 + |y - \frac{1}{2} a_k^- t| + (t/2)^{1/2})^{-3/2} \\
&+
e^{-\frac{(x - y - a_j^- (t/2))^2}{M(t/2)}} 
(1 + |y - \frac{1}{2}a_k^- t| + |x - a_j^- (t/2) - a_k^- (t/2)| + (t/2)^{1/2})^{-3/2}.
\end{aligned}
\end{equation}
In the event that $|x| \ge t/c$, for $c$ sufficiently small, we have 
exponential decay in both $|x|$ and $t$ for the first estimate in 
(\ref{charpartsvsquaredfirstbalance1}), while for the second we have 
an estimate by 
\begin{equation*}
C_1 t^{-1/2} (1 + t)^{-3/4} (1 + |x| + t)^{-3/2},
\end{equation*} 
which again is sufficient.  For the second estimate on $v(y, s)$
we have integrals
\begin{equation} \label{charpartsvsquaredfirstnl2}
\begin{aligned}
\int_{-|a_1^-| \frac{t}{2}}^{0} &t^{-1} e^{-\frac{(x - y - a_k^- (t/2))^2}{M(t/2)}}
(1 + |y|)^{-1/2} (1 + |y| + (t/2))^{-1/2}
(1 + |y| + (t/2)^{1/2})^{-1/2} (1 + (t/2))^{-3/4} dy \\
&\le
C t^{-1/2} (1 + t)^{-3/2},
\end{aligned}
\end{equation}
which is sufficient for $t \ge c|x|$, for any constant $c > 0$.  The 
case $|x| \ge t/c$ can be analyzed similarly as in the immediately preceding case.
For the second integral on the right hand side of 
(\ref{charpartsvsquared}), we observe that $G^j (x,0;y)$ is a delta
function with mass at $x = y$, and consequently, we have an estimate 
by 
\begin{equation*}
C |v (x, t) v_x (x, t)| \le
C (1 + t)^{-3/4} \Big[t^{-1/2} (\psi_1 (x, t) + \psi_2 (x,t)) + \psi_3 (x, t) + \psi_4 (x, t) \Big],
\end{equation*}
which is sufficient.  For the third integral on the right hand side of 
(\ref{charpartsvsquared}), and for the first estimate on  
$(\partial_s + a_l^- \partial_y) v (y, s)$,
we have integrals
\begin{equation} \label{charpartsvsquaredthird}
\begin{aligned}
&\int_{t/2}^t \int_{-\infty}^{+\infty} (t-s)^{-1} e^{-\frac{(x - y - a_j^- (t-s))^2}{M(t-s)}}
s^{-1} (1+s)^{-1/2} (1 + |y - a_l^- s)| + s^{1/2})^{-3/2} dy ds \\
&\le
C t^{-1/2} (1 + t)^{-1/4}
\int_{t/2}^t \int_{-\infty}^{+\infty} (t-s)^{-1} e^{-\frac{(x - y - a_j^- (t-s))^2}{M(t-s)}}
s^{-1/2} (1+s)^{-1/4} (1 + |y - a_l^- s)| + s^{1/2})^{-3/2} dy ds.
\end{aligned}
\end{equation}
This last integral has already been considered in the analysis of 
(\ref{gf1nl1}), and we obtain an estimate by 
\begin{equation*}
C t^{-1/2} (1 + t)^{-1/4} \psi_1 (x, t).
\end{equation*}
We proceed in precisely the same manner for the second and third estimates 
on $(\partial_s + a_l^- \partial_y) v (y, s)$, obtaining an estimate by 
\begin{equation*}
C \Big[t^{-1/2} (\bar{\psi}_1^{l,-} (x, t) + \psi_2 (x, t)) + \psi_3 (x, t) + \psi_4 (x, t) \Big].
\end{equation*}

  For the fourth estimate 
on $(\partial_s + a_l^- \partial_y) v (y, s)$, we have integrals 
\begin{equation} \label{charpartsvsquaredthird2}
\begin{aligned}
&\int_{t/2}^t \int_{-\infty}^{+\infty} (t-s)^{-1} e^{-\frac{(x - y - a_j^- (t-s))^2}{M(t-s)}}
(1+s)^{-3/4} (1 + |y| + s)^{-1} (1 + |y|)^{-1} dy ds. 
\end{aligned}
\end{equation}   
In this case, we observe the inequality 
\begin{equation} \label{charpartsvsquaredthird2balance1}
\begin{aligned}
&e^{-\frac{(x - y - a_j^- (t-s))^2}{M(t-s)}} (1 + |y|)^{-1} \\
&\le
C\Big[ e^{- \epsilon \frac{(x - y - a_j^- (t-s))^2}{M(t-s)}}
e^{-\frac{(x - a_j^- (t-s))^2}{\bar{M}(t-s)}} (1 + |y|)^{-1}
+
e^{-\frac{(x - y - a_j^- (t-s))^2}{M(t-s)}} (1 + |y| + |x - a_j^- (t-s)|)^{-1}
\Big].
\end{aligned}
\end{equation}
For the first estimate in (\ref{charpartsvsquaredthird2balance1}), we 
have integrals 
\begin{equation} \label{charpartsvsquaredthird2balance1first}
\begin{aligned}
&\int_{t/2}^t \int_{-|a_1^-| s}^{0} (t-s)^{-1} 
e^{- \epsilon \frac{(x - y - a_j^- (t-s))^2}{M(t-s)}}
e^{-\frac{(x - a_j^- (t-s))^2}{\bar{M}(t-s)}}
(1+s)^{-3/4} (1 + |y| + s)^{-1} (1 + |y|)^{-1} dy ds,
\end{aligned}
\end{equation}   
for which we have two cases to consider, $a_j^- < 0$ and $a_j^- > 0$,
and we focus on the former.  For $a_j^- < 0$ and additionally
$|x| \ge |a_j^-| t$, we have no cancellation between $x - a_j^- t$
and $a_j^- s$, and we obtain an estimate by 
\begin{equation} \label{charpartsvsquaredthird2balance1firstnocanc1}
\begin{aligned}
C_2 (1 + t)^{-7/4} \ln(e+t) e^{-\frac{(x - a_j^- t)^2}{Lt}} 
&\int_{t/2}^{t-1} (t-s)^{-1} ds
+ 
C_2 (1 + t)^{-7/4} \ln(e+t) e^{-\frac{(x - a_j^- t)^2}{Lt}} 
\int_{t-1}^t (t-s)^{-1/2} ds \\
&\le
C (1+t)^{-7/4} [\ln (e+t)]^2 e^{-\frac{(x - a_j^- t)^2}{Lt}}, 
\end{aligned}
\end{equation}
which is sufficient.  For $|x| \le |a_j^-| t$, we observe the 
inequality (\ref{gf1nl2balance3}) with $\gamma = 1/2$, for which we estimate
\begin{equation}
\begin{aligned}
C (1 + t)^{-7/4} \ln(e+t) (1 + |x|)^{-1/2} 
\int_{t/2}^t (t-s)^{-1/2} e^{-\frac{(x - a_j^- (t-s))^2}{\bar{M}(t-s)}} ds
\le
C (1 + t)^{-7/4} \ln(e+t) (1 + |x|)^{-1/2},
\end{aligned}
\end{equation}
which is sufficient for $t \ge |x|/|a_j^-|$.  
For the second estimate in (\ref{charpartsvsquaredthird2balance1}), we 
have integrals 
\begin{equation*}
\begin{aligned}
&\int_{t/2}^t \int_{-|a_1^-| s}^{0} (t-s)^{-1} 
e^{-\frac{(x - y - a_j^- (t-s))^2}{M(t-s)}} (1 + |y| + |x - a_j^- (t-s)|)^{-1}
(1 + s)^{-7/4} dy ds \\
&\le
C (1+t)^{-7/4} \int_{t/2}^t (t-s)^{-1/2} (1 + |x - a_j^- (t-s)|)^{-1} ds
\le
C (1+t)^{-7/4},
\end{aligned}
\end{equation*} 
which is sufficient for $|x| \le Kt$, some constant $K$.  For $|x| \ge Kt$,
we estimate  
\begin{equation*}
\begin{aligned}
&\int_{t/2}^t \int_{-|a_1^-| s}^{0} (t-s)^{-1} 
e^{-\frac{(x - y - a_j^- (t-s))^2}{M(t-s)}} (1 + |y| + |x - a_j^- (t-s)|)^{-1}
(1 + s)^{-7/4} dy ds \\
&\le
C (1+t)^{-7/4} (1 + |x| + t)^{-1} \int_{t/2}^t (t-s)^{-1/2} ds
\le
C (1+t)^{-5/4} (1 + |x| + t)^{-1}.
\end{aligned}
\end{equation*} 
For the final estimate on $(\partial_s + a_l^- \partial_y) v(y,s)$, 
we have integrals
\begin{equation}
\begin{aligned}
&\int_{t/2}^t \int_{-|a_1^-| s}^{0} (t-s)^{-1} e^{-\frac{(x - y - a_j^- (t-s))^2}{M(t-s)}}
(1+s)^{-3/4} (1 + |y| + s)^{-7/4} dy ds \\
&\le
C (1 + t)^{-5/2}
\int_{t/2}^t (t-s)^{-1/2} ds
\le 
C (1 + t)^{-2},
\end{aligned}
\end{equation}
which is sufficient for $|x| \le Kt$.  Since $|y|$ is bounded by $|a_1^-| s$, 
for $|x| \ge Kt$ and $K$ sufficiently large, we have exponential decay in 
both $|x|$ and $t$.

{\it Case 2: $l = j$.}  For the case $l = j$, we have additional 
decay at rate $(t-s)^{-1/2}$, from which we immediately recover the claimed
estimates.


{\it (\ref{Ginteract}(v)), Nonlinearity} 
$|\frac{\partial \bar{u}^\delta}{\partial \delta} \delta|^2$.
We next consider integrals
\begin{equation*}
\int_0^t \int_{-\infty}^{+\infty} 
(\partial_t + a_l^\pm \partial_x) \tilde{G}^j_y (x,t-s;y) (1+s)^{-1} e^{-\eta |y|} dy ds,
\end{equation*}
where in this case additional $y$-derivatives on the nonlinearity give no 
additional decay.

{\it Case 1: $l \ne j$.}  For the case $l \ne j$, we have integrals
\begin{equation} \label{charderlastnonlin}
\begin{aligned}
&\int_0^{t-1} \int_{-\infty}^{+\infty} (t-s)^{-3/2} 
e^{-\frac{(x - y - a_j^- (t-s))^2}{M(t-s)}} (1 + s)^{-1} e^{-\eta |y|} dy ds \\
&+
\int_{t-1}^t \int_{-\infty}^{+\infty} (t-s)^{-1} 
e^{-\frac{(x - y - a_j^- (t-s))^2}{M(t-s)}} (1 + s)^{-1} e^{-\eta |y|} dy ds.
\end{aligned}
\end{equation}
In either case, we observe the inequality (\ref{gf1nl3balance1}).  For the first
estimate in (\ref{gf1nl3balance1}), we have, upon integration of $e^{-\eta |y|}$,
integrals
\begin{equation} \label{charderlastnonlinfirst}
\begin{aligned}
&\int_0^{t-1} (t-s)^{-3/2} 
e^{-\frac{(x - a_j^- (t-s))^2}{M(t-s)}} (1 + s)^{-1} ds \\
&+
\int_{t-1}^t (t-s)^{-1/2} 
e^{-\frac{(x - a_j^- (t-s))^2}{M(t-s)}} (1 + s)^{-1} ds.
\end{aligned}
\end{equation}
Focusing as in previous cases on the subcase $a_j^- < 0$, we
first observe that for $|x| \ge |a_j^-| t$, there is no cancellation 
between $x - a_j^- t$ and $a_j^- s$, and we consequently have 
an estimate by
\begin{equation} \label{charderlastnonlinfirstnocanc}
\begin{aligned}
C t^{-3/2} e^{-\frac{(x - a_j^- t)^2}{Lt}} &\int_0^{t/2}  
(1 + s)^{-1} ds 
+
C_2 (1+t)^{-1} e^{-\frac{(x - a_j^- t)^2}{Lt}}
\int_{t/2}^{t} (t-s)^{-1/2} 
e^{-\frac{(a_j^- s)^2}{M(t-s)}} ds \\
&\le
C t^{-3/2} [\ln (e+t)] e^{-\frac{(x - a_j^- t)^2}{Lt}},
\end{aligned}
\end{equation} 
where in this last inequality we have observed that for $s \in [t/2, t]$,
we have
\begin{equation*}
e^{-\frac{(a_j^- s)^2}{M(t-s)}} \le e^{-\eta_1 t},
\end{equation*}
for $\eta_1 > 0$.  For $|x| \le |a_j^-| t$, we divide the analysis into cases, 
$s \in [0, t/2]$, $s \in [t/2, t-1]$, and $s \in [t-1, t]$.  For $s \in [0, t/2]$,
we observe the inequality (\ref{gf1nl3balance2}).  For the second estimate
in (\ref{gf1nl3balance2}), we obtain an estimate by 
\begin{equation*}
C_1 t^{-3/2} (1+|x - a_j^- t|)^{-1} 
\int_0^{t/2} e^{-\frac{(x - a_j^- (t-s))^2}{M(t-s)}} ds
\le
C t^{-1} (1+|x - a_j^- t|)^{-1},
\end{equation*}
which is sufficient for $t \ge |x|/|a_j^-|$ and $l \ne j$.  For
the first estimate in (\ref{gf1nl3balance2}), we proceed as in 
(\ref{charderlastnonlinfirstnocanc}).  For $s \in [t/2, t]$, we 
observe the inequality (\ref{gf1nl2balance3}) with $\gamma = 1$,
for which we have an estimate by 
\begin{equation*}
\begin{aligned}
C_2 (1+t)^{-1} |x|^{-1} &\int_{t/2}^{t-1} (t-s)^{-1/2}
e^{-\frac{(x - a_j^- (t-s))^2}{2 \bar{M}(t-s)}} ds
+
C_3 (1+t)^{-1} |x|^{-1} \int_{t-1}^t ds \\
&\le
C (1 + t)^{-1} |x|^{-1}.
\end{aligned}
\end{equation*}

{\it Case 2: $l = j$.}  For the case $l = j$, we have additional 
decay at rate $(t-s)^{-1/2}$, from which we immediately recover the claimed
estimates.  

This ends the proof of estimate (\ref{Ginteract}(v)) for the leading 
order convection kernel $\tilde{G}^j (x, t; y)$.

{\it (\ref{Ginteract}(v)--(vi)), remainder estimates.}
In our proofs of (\ref{Ginteract}(v)--(vi)), we have considered only the 
leading order convection kernel
\begin{equation*}
\tilde{G}^j (x, t; y) = ct^{-1/2} e^{-\frac{(x - y - a_j^- t)^2}{4 \beta_j t}}. 
\end{equation*}
We must also consider the remaining three scattering estimates (terms in 
$S (x, t; y)$) and additionally the remainder estimates $R(x, t; y)$.
Beginning with the remainder estimate, we focus our attention on the 
nonlinearity $(\varphi^k (y, s)^2)_y$ (analysis of the remaining 
nonlinearities is similar).  We have the decomposition
\begin{equation} \label{residualthreeparts}
\begin{aligned}
\int_0^t &\int_{-\infty}^{+\infty} (\partial_t + a_l^- \partial_x) R(x, t-s; y) 
(\varphi^k (y, s)^2)_y dy ds \\
&=
\int_0^{\sqrt{t}} \int_{-\infty}^{+\infty} (\partial_t + a_l^- \partial_x) R(x, t-s; y) 
(\varphi^k (y, s)^2)_y dy ds \\
&+
\int_{\sqrt{t}}^{t-\sqrt{t}} \int_{-\infty}^{+\infty} (\partial_t + a_l^- \partial_x) R(x, t-s; y) 
(\varphi^k (y, s)^2)_y dy ds \\
&+
\int_{t-\sqrt{t}}^t \int_{-\infty}^{+\infty} (\partial_t + a_l^- \partial_x) R(x, t-s; y) 
(\varphi^k (y, s)^2)_y dy ds,
\end{aligned}
\end{equation} 
where for $y \le 0$ and $a_l^- < 0$,
\begin{equation} \label{charderresidual}
\begin{aligned}
(\partial_t + a_l^- \partial_x) &R (x,t;y) = 
\sum_{j=1}^J \mathbf{O}(e^{-\eta t})\delta_{x-\bar
a_j^* t}(-y) +
\mathbf{O}(e^{-\eta(|x-y|+t)})\\
&+
\mathbf {O} \left( (t+1)^{-3/2} e^{-\eta x^+}
+e^{-\eta|x|} \right)
t^{-1} (t+1)^{1/2}
e^{-(x-y-a_l^{-} t)^2/Mt} \\
&+\sum_{k \ne l} \mathbf {O} \left( (t+1)^{-1} e^{-\eta x^+}
+e^{-\eta|x|} \right)
t^{-1} (t+1)^{1/2}
e^{-(x-y-a_k^{-} t)^2/Mt} \\
&+
\chi_{\{ |a_k^{-} t|\ge |y|
\}} \mathbf{O}((t+1)^{-1/2} t^{-1})
e^{-(x-a_l^{-}(t-|y/a_k^-|))^2/Mt}
e^{-\eta x^+} \\
&+ \sum_{a_k^{-} > 0, \, a_j^{-} < 0, \, j\ne l} \chi_{\{ |a_k^{-} t|\ge |y|
\}} \mathbf{O}((t+1)^{-1/2} t^{-1})
e^{-(x-a_j^{-}(t-|y/a_k^-|))^2/Mt}
e^{-\eta x^+} \\
&+ \sum_{a_k^{-} > 0, \, a_j^{+} > 0} \chi_{\{ |a_k^{-} t|\ge |y|
\}} \mathbf{O}((t+1)^{-1/2} t^{-1})
e^{-(x-a_j^{+}(t-|y/a_k^{-}|))^2/Mt}
e^{-\eta x^-}. 
\end{aligned}
\end{equation}
In each estimate with decay $t^{-3/2}$ or $t^{-2}$, we have time decay 
better than that of characteristic derivatives of $\tilde{G}^j (x, t; y)$,
and we can proceed as in the above analyses.  For terms 
$\sum_{j=1}^J \mathbf{O}(e^{-\eta t})\delta_{x-\bar a_j^* t}(-y)$, we proceed
as in the proof of (\ref{Hnonlinear}), in which the interactions arising from 
our relaxation of strict parabolicity are considered.  
The only genuinely new term is the exponentially decaying 
contribution, which has reduced decay in $t$. 
Focusing on this term, we integrate by parts in $y$,
observing that in the Lax and overcompressive cases 
differentiation with respect to $y$ improves $t$ decay by 
a factor $t^{-1/2}$ (this is a fundamental point of difference
between the Lax and overcompressive cases considered here, and 
the undercompressive case).  We have integrals
\begin{equation} \label{cleanup}
\begin{aligned}
e^{-\eta |x|} \int_0^{\sqrt{t}} &\int_{-\infty}^{+\infty} 
(t-s)^{-1} e^{-\frac{(x - y - a_j^- (t-s))^2}{M(t-s)}} 
(1+s)^{-1} e^{-\frac{(y - a_k^- s)^2}{Ms}} dy ds \\
&\le
C e^{-\eta |x|} t^{-1} \int_0^{\sqrt{t}}
(1+s)^{-1} s^{1/2} e^{-\frac{(x - a_j^- (t-s) - a_k^- s)^2}{Mt}} ds \\
&\le
C e^{-\eta |x|} t^{-1} e^{-\frac{(x - a_j^- t)^2}{Mt}}
\int_0^{\sqrt{t}} (1+s)^{-1} s^{1/2} ds \\
& \le
C (1+t)^{-3/4} e^{-\eta |x|} e^{-\frac{(x - a_j^- t)^2}{Mt}},
\end{aligned}
\end{equation}
for which we observe that 
\begin{equation*}
e^{-\eta |x|} e^{-\frac{(x - a_j^- t)^2}{Mt}}
\le
C e^{-\eta_1 |x|} e^{-\eta_2 t}.
\end{equation*}
For the third integral in (\ref{residualthreeparts}), we 
estimate
\begin{equation*}
\begin{aligned}
e^{-\eta |x|} \int_{t-\sqrt{t}}^t &\int_{-\infty}^{+\infty} 
(t-s)^{-1/2} e^{-\frac{(x - y - a_j^- (t-s))^2}{M(t-s)}} 
(1+s)^{-3/2} e^{-\frac{(y - a_k^- s)^2}{Ms}} dy ds \\
&\le
C e^{-\eta |x|} t^{-1/2} \int_{t-\sqrt{t}}^t
(1+s)^{-3/2} s^{1/2} e^{-\frac{(x - a_j^- (t-s) - a_k^- s)^2}{Mt}} ds \\
&\le
C e^{-\eta |x|} t^{-1/2} (1+t)^{-1} e^{-\frac{(x - a_k^- t)^2}{Lt}}
\int_{t-\sqrt{t}}^t e^{- \epsilon \frac{(x - a_j^- (t-s) - a_k^- s)^2}{Mt}}  ds \\
& \le
C (1+t)^{-1} e^{-\eta |x|} e^{-\frac{(x - a_k^- t)^2}{Mt}},
\end{aligned}
\end{equation*}
which is sufficient, precisely as above.  
For the second integral in (\ref{residualthreeparts}), and for 
$t$ large enough so that $\sqrt{t} < t/2$, we estimate
(integrating by parts in $y$ for $s \in [\sqrt{t}, t/2]$)  
\begin{equation} \label{residualthreepartssecond}
\begin{aligned}
e^{-\eta |x|} \int_{\sqrt{t}}^{t/2} &\int_{-\infty}^{+\infty} 
(t-s)^{-1} e^{-\frac{(x - y - a_j^- (t-s))^2}{M(t-s)}} 
(1+s)^{-1} e^{-\frac{(y - a_k^- s)^2}{Ms}} dy ds \\
&+
e^{-\eta |x|} \int_{t/2}^{t-\sqrt{t}} \int_{-\infty}^{+\infty} 
(t-s)^{-1/2} e^{-\frac{(x - y - a_j^- (t-s))^2}{M(t-s)}} 
(1+s)^{-3/2} e^{-\frac{(y - a_k^- s)^2}{Ms}} dy ds \\
&\le
C_2 e^{-\eta |x|} t^{-1/2} \int_{t/2}^{t-\sqrt{t}} (t-s)^{-1/2} 
(1+s)^{-1} s^{1/2} e^{-\frac{(x - a_j^- (t-s) - a_k^- s)^2}{Mt}} ds \\
&+
C_2 e^{-\eta |x|} t^{-1/2} \int_{t/2}^{t-\sqrt{t}} 
(1+s)^{-3/2} s^{1/2} e^{-\frac{(x - a_j^- (t-s) - a_k^- s)^2}{Mt}} ds.
\end{aligned}
\end{equation}
For $s \in [\sqrt{t}, t/2]$, we first observe that for the case 
$j = k$, we have an estimate by 
\begin{equation} \label{residualthreepartssecondnocanc1}
C (1 + t)^{-1/2} e^{-\eta |x|} e^{-\frac{(x - a_j^- t)^2}{Mt}},
\end{equation}
which is sufficient, as observed above.  For the case $j \ne k$,
we recall the inequality (\ref{gf1nl4balance1}) with $\gamma = 1/2$.
For the first estimate in (\ref{gf1nl4balance1}), we proceed 
as in (\ref{residualthreepartssecondnocanc1}), while for the second
we estimate
\begin{equation*}
\begin{aligned}
C_1 &e^{-\eta |x|} t^{-1} (1 + |x - a_j^- t|)^{-1/2} 
\int_{\sqrt{t}}^{t/2} e^{-\frac{(x - a_j^- (t-s) - a_k^- s)^2}{Mt}} ds \\
&\le
C (1+t)^{-1/2} e^{-\eta |x|} (1 + |x - a_j^- t|)^{-1/2}.
\end{aligned}
\end{equation*}
Observing the estimate 
\begin{equation*}
e^{-\eta |x|} (1 + |x - a_j^- t|)^{-1/2} \le
C e^{-\eta_1 |x|} (1 + t)^{-1/2},
\end{equation*}
we observe that this estimate is sufficient.
For $s \in [t/2, t-\sqrt{t}]$, we first observe that 
in the event that $j = k$, this last integral 
in (\ref{residualthreepartssecond}) provides an estimate by 
\begin{equation*}
C t^{-1/2} e^{-\eta |x|} e^{-\frac{(x - a_j^- t)^2}{Mt}},
\end{equation*}
while in the event that $j \ne k$, we have an estimate by
\begin{equation*}
C (1+t)^{-1} e^{-\eta |x|},
\end{equation*}
either of which is sufficient.

{\it (\ref{Ginteract}(v)--(vi)), full scattering estimates.}
In the scattering estimates $S(x, t; y)$ of Proposition \ref{greenbounds},
we have three corrections to $\tilde{G}^j (x, t; y)$, respectively from 
convection, reflection, and transmission.  For convection, the term is 
\begin{equation*}
\tilde{G}_c^j (x, t; y) = ct^{-1/2} e^{-\frac{(x - y - a_j^- t)^2}{4 \beta_j t}}
\Big{(}\frac{e^{-x}}{e^x + e^{-x}} \Big{)},
\end{equation*}  
$a_j^- > 0$,
which satisfies 
\begin{equation*}
|(\partial_t + a_j^- \partial_x) \tilde{G}_c^j (x, t; y)| \le
C \Big{[}t^{-3/2} + t^{-1/2} e^{-\eta |x|} \Big{]} e^{-\frac{(x - y - a_j^- t)^2}{4 \beta_j t}}.
\end{equation*}
The estimate with $t^{-3/2}$ decay is precisely as in the case of 
$\tilde{G}^j$, and can be analyzed similarly.  For the estimate with 
exponential decay in $|x|$, we can proceed precisely as with the 
exponentially decaying term arising in the analysis of $R(x, t; y)$.
The remaining corrections, from the reflection and transmission terms in 
$S(x, t; y)$ can be analyzed similarly.

In our analysis of $\tilde{G}^j (x, t; y)$, we employed the relation 
\begin{equation*}
\partial_x \tilde{G}^j (x, t; y) = - \partial_y \tilde{G}^j (x, t; y)
\end{equation*}
(see the argument of (\ref{charparts}), in which the characteristic derivative
$(\partial_t + a_j^- \partial_x) \tilde{G}^j (x, t-s; y)$ is converted into 
a characteristic derivative in the variables of integration 
$-(\partial_s + a_j^- \partial_y) \tilde{G}^j (x, t-s; y)$).  Designating
the scattering term arising from reflection as 
\begin{equation*}
\tilde{G}_R^{j, k} (x, t; y) 
=
c_R t^{-1/2} e^{-\frac{(x - (a_j^-/a_k^-)y - a_j^- t)^2}{4 \bar{\beta}_{jk}^- t}}
\Big{(}\frac{e^{-x}}{e^x + e^{-x}}\Big{)},
\end{equation*}
we observe that the analogous estimate is  
\begin{equation*}
\partial_x \tilde{G}_R^{j, k} (x, t; y)
=
-\frac{a_k^-}{a_j^-} \partial_y \tilde{G}_R^{j, k} (x, t; y) 
+
{\mathbf O} (t^{-1/2} e^{-\eta |x|}) 
e^{-\frac{(x - (a_j^-/a_k^-)y - a_j^- t)^2}{4 \bar{\beta}_{jk}^- t}},
\end{equation*}
from which we see that the conversion from a characteristic derivative in 
variables $x$ and $t$ to a characteristic derivative in variables $y$ and $s$ takes the form 
\begin{equation*}
\begin{aligned}
(\partial_t + a_j^- \partial_x) \tilde{G}_R^{j, k} (x, t-s; y)
& =
-(\partial_s + a_k^- \partial_y) \tilde{G}_R^{j, k} (x, t-s; y) \\
& +
{\mathbf O} ((t-s)^{-1/2} e^{-\eta |x|}) 
e^{-\frac{(x - (a_j^-/a_k^-)y - a_j^- (t-s))^2}{4 \bar{\beta}_{jk}^- (t-s)}}.
\end{aligned}
\end{equation*}
In this way, we again have precisely the improved decay required by 
the analysis, namely 
\begin{equation*}
(\partial_s + a_k^- \partial_y) \tilde{G}_R^{j, k} (x, t-s; y) 
=
{\mathbf O} ((t-s)^{-3/2}) 
e^{-\frac{(x - (a_j^-/a_k^-)y - a_j^- (t-s))^2}{4 \bar{\beta}_{jk}^- (t-s)}}.
\end{equation*}
The argument involving transmission terms is entirely similar.

Finally, in each case in which a characteristic derivative in $x$ and $t$ 
is shifted to one in $y$ and $s$ an exponentially decaying error term 
arises, such as (from reflection)
\begin{equation*}
{\mathbf O} (t^{-1/2} e^{-\eta |x|}) 
e^{-\frac{(x - (a_j^-/a_k^-)y - a_j^- t)^2}{4 \bar{\beta}_{jk}^- t}}.
\end{equation*}
In all cases, these terms can be analyzed as were the the exponentially 
decaying corrections in $(\partial_t + a_j^- \partial_x) R(x, t; y)$
(see (\ref{cleanup})).

This completes the proof of Lemma \ref{nonlinearintegralestimates}. 
\hfill $\square$

\section{Nonlinear integral Estimates II}\label{integralestimatesNLII}

Finally, we complete the paper by establishing
the nonlinear integral estimates of Proposition \ref{Hnonlinear}.

\noindent {\bf Proof of Lemma \ref{Hnonlinear}.} To show (\ref{Hinteract})(i)
we need to estimate 
\begin{equation}\label{quant}
\begin{aligned}
\Big| \int_0^t\int_{-\infty}^{+\infty}\mathcal{R}_j^*(x) 
&\mathcal{O}(e^{-\eta_0 (t-s)}) \delta_{x-\bar a_j^*
(t-s)}(-y) \mathcal{L}_j^{*t}(y)\Upsilon(y,s) dyds\Big|\le \\
&\qquad \qquad 
\qquad \qquad 
C\int_0^t (e^{-\eta_0 (t-s)})|\Upsilon(-x+\bar a_j^*
(t-s),s)|ds
\end{aligned}
\end{equation}
for the various sources $\Upsilon(y,s)$ arising in the bounds
for $\fcal$, $\Phi$.

For example, the typical term 
$$
|v||v_{yy}|(y,s)\le C(1+s)^{-1/2} (\psi_1+\psi_2)(y,s)
$$
arising in the bounds for $\fcal$ leads to sources
\begin{equation}
\begin{aligned}
\Upsilon_1(y,s)&= (1+s)^{-1/2} \psi_1(y,s)\\
&=
(1+s)^{-1/2} (1+|y-a_i^\pm t|)^{-3/2}
\end{aligned}
\end{equation}
and 
\begin{equation}
\begin{aligned}
\Upsilon_2(y,s)&=(1+s)^{-1/2} \psi_2(y,s)\\
&\le (1+s)^{-1/2} (1+|y|)^{-1/2}(1+|y|+t)^{-1/2}(1+|y|+t^{1/2})^{-1/2}.
\end{aligned}
\end{equation}

Substituting $\Upsilon_1$ in (\ref{quant}), we obtain
\begin{equation}
\begin{aligned}
\int_0^t e^{-\eta_0(t-s)} &(1+s)^{-\frac12}
(1+|x-\bar a_j^* (t-s)- a_i^\pm s |)^{-\frac32} ds =\\
&\int_0^t e^{-\eta_0(t-s)} (1+s)^{-\frac12}
(1+|x- a_i^-t -(\bar a_j^*-a_i^\pm) (t-s)|)^{-\frac32} ds,
\end{aligned}
\end{equation}
which , by (\ref{abineq}), is smaller than
$$
\int_0^t e^{-\eta_0(t-s)} (1+s)^{-\frac12}
(1+|x- a_i^-t |)^{-\frac32}(|1+(\bar a_j^*-a_i^-) (t-s)|)^{\frac32} ds,
$$
which in turn (absorbing $t-s$ powers in $e^{\frac{-\eta_0 (t-s)}{2}}$) 
is smaller than
\begin{equation}
\begin{aligned}
(1+|x- a_i^-t |)^{-\frac32} 
\int_0^t e^{-\frac{\eta_0}2(t-s)} (1+s)^{-\frac12}ds
&\le C(1+|x- a_i^-t |)^{-\frac32}(1+t)^{-1/2}\\
&=C(1+t)^{-1/2}\psi_1(x,t).
\end{aligned}
\end{equation}

Similarly, substituting $\Upsilon_2$ in (\ref{quant}), 
observing that
$$
(1+|y|+t^{1/2})^{-1/2}\sim
\min\{(1+|y|)^{-1/2}, (1+s)^{-1/4}\},
$$
and following the same procedure, we obtain an estimate of
\begin{equation}
\begin{aligned}
C(1+t)^{-1/2}(1+|x|)^{-1/2}(1+|x|+t)^{-1/2}
&\min\{ (1+|x|)^{-1/2}, (1+t)^{-1/4} \} \sim \\
&\quad C(1+t)^{-1/2}(1+|x|)^{-1/2}(1+|x|+t)^{-1/2} (1+|x|+t^{1/2})^{-1/2}\\
&\le C(1+t)^{-1/2} (\psi_1+\psi_2)(x,t),
\end{aligned}
\end{equation}
where, in the final step, we have used the fact that,
for $\chi(x,t)=0$,
$$
(1+|x|)^{-1/2}(1+|x|+t)^{-1/2} (1+|x|+t^{1/2})^{-1/2}\sim
(1+|x|+t)^{-3/2}\le C\psi_1(x,t).$$

Bounds for other cases follow similarly.
  \hfill $\square$

{\it Acknowledgements.} The authors were partially supported by the National Science Foundation,
under grants DMS--0500988 (Howard), 
and grants DMS--0070765 and DMS--0300487 (Raoofi and Zumbrun).

\noindent Peter HOWARD \\
\noindent Department of Mathematics \\
\noindent Texas A\&M University \\
\noindent College Station, TX 77843 \\
\noindent phoward@math.tamu.edu \\
\noindent and \\
\noindent Mohammadreza RAOOFI \\
\noindent Max Planck Institute for Mathematics in the Sciences\\
\noindent Inselstraße 22-26\\
\noindent D-04103 Leipzig, Germany\\
\noindent mraoofi@indiana.edu \\
\noindent and \\
\noindent Kevin ZUMBRUN \\
\noindent Department of Mathematics \\
\noindent Indiana University \\
\noindent Bloomington, IN  47405-4301\\
\noindent kzumbrun@indiana.edu \\

\end{document}